# Modulated equations of Hamiltonian PDEs and dispersive shocks

Sylvie Benzoni-Gavage[*]   Colin Mietka[†]   L. Miguel Rodrigues[‡]

November 7, 2020

**Abstract**

Motivated by the ongoing study of dispersive shock waves in non integrable systems, we propose and analyze a set of wave parameters for periodic waves of a large class of Hamiltonian partial differential systems — including the generalized Korteweg–de Vries equations and the Euler–Korteweg systems — that are well-behaved in both the small amplitude and large wavelength limits. We use this parametrization to determine fine asymptotic properties of the associated modulation systems, including detailed descriptions of eigenmodes. As a consequence, in the solitary wave limit we prove that modulational instability is decided by the sign of the second derivative — with respect to speed, fixing the endstate — of the Boussinesq moment of instability; and, in the harmonic limit, we identify an explicit modulational instability index, of Benjamin–Feir type.

## Contents

[*]This work has been supported by the French National Research Agency projects NABUCO, grant ANR-17-CE40-0025, BoND, grant ANR-13-BS01-0009-01, and by the LABEX MILYON (ANR-10-LABX-0070) of Université de Lyon, within the program 'Investissements d'Avenir' (ANR-11-IDEX-0007) operated by the French National Research Agency (ANR).
[†]The doctoral scholarship of CM was directly supported by ANR-13-BS01-0009-01.
[‡]Research of LMR has received funding from the ANR project BoND (ANR-13-BS01-0009-01) and the city of Rennes.

# 1 Introduction

Motivated by the study of dispersive shock waves, we investigate some of the universal properties of modulated equations, for a large class of Hamiltonian systems of partial differential equations (PDE) that contains several models of mathematical physics and in particular generalized versions of the Korteweg–de Vries (KdV) equation and dispersive modifications of the Euler equations for compressible fluids - among which we find the fluid formulation via Madelung's transform of the nonlinear Schrödinger (NLS) equations. To place our results in context, we first recall what are modulation systems, dispersive shock waves and their expected role in dispersive regularization. Large parts of this preliminary discussion are exploratory and conjectural since we are still lacking a rigorous analysis of dispersive shock waves and vanishing dispersive limits at the level of generality considered here. Indeed the present analysis is precisely designed as a first step towards a general mathematically rigorous theory, still to come, and some elements of the preliminary discussion might be thought of as a roadmap for this ultimate goal.

## Hamiltonian systems of Korteweg type

As in [BGNR13, BGNR14, BGMR16, BGMR20], we consider some abstract systems of the form

$$\partial_t \mathbf{U} = \partial_x (\mathbf{B}\, \delta\mathcal{H}[\mathbf{U}])\,. \tag{1}$$

where the unknown $\mathbf{U}$ takes values in $\mathbb{R}^N$ with $N \in \{1, 2\}$, $\mathbf{B}$ is a *symmetric* and *nonsingular* matrix so that $\mathbf{B}\partial_x$ is a skew-symmetric differential operator, and $\delta\mathcal{H}[\mathbf{U}]$ denotes the variational derivative at $\mathbf{U}$ of $\mathcal{H} = \mathcal{H}(\mathbf{U}, \mathbf{U}_x)$. We specialize to classes of systems satisfying more precise conditions described in Assumption 1, sufficiently large to include quasilinear versions of the Korteweg–de Vries equation and the abovementioned Euler–Korteweg systems, hence also hydrodynamic formulations of the nonlinear Schrödinger equations. In the former and henceforth spatial derivatives are denoted either as $\partial_x$ or as $_x$.

Associated with the invariance of $\mathcal{H}$ by time and spatial translations comes the fact that smooth solutions of (1) also satisfy the local conservations laws

$$\partial_t(\mathcal{H}(\mathbf{U}, \mathbf{U}_x)) = \partial_x \left( \tfrac{1}{2}\delta\mathcal{H}[\mathbf{U}] \cdot \mathbf{B}\delta\mathcal{H}[\mathbf{U}] + \nabla_{\mathbf{U}_x}\mathcal{H}(\mathbf{U}, \mathbf{U}_x) \cdot \partial_x(\mathbf{B}\,\delta\mathcal{H}[\mathbf{U}]) \right) \tag{2}$$

$$\partial_t(\mathcal{Q}(\mathbf{U})) = \partial_x \left( \delta\mathcal{Q}[\mathbf{U}] \cdot \mathbf{B}\delta\mathcal{H}[\mathbf{U}] + \nabla_{\mathbf{U}_x}\mathcal{H}(\mathbf{U}, \mathbf{U}_x) \cdot \partial_x(\mathbf{B}\,\delta\mathcal{Q}[\mathbf{U}]) - \mathcal{H}(\mathbf{U}, \mathbf{U}_x) \right) \tag{3}$$

for the *Hamiltonian* density $\mathcal{H}$, generating the time evolution, and the *impulse* $\mathcal{Q}$, given by $\mathcal{Q}(\mathbf{U}) := \tfrac{1}{2}\mathbf{U} \cdot \mathbf{B}^{-1}\mathbf{U}$, generating spatial translations. See Section 2 for details.

## Modulated equations

Modulated equations are expected to govern the evolution of modulated periodic wavetrains (also called weakly deformed soliton lattices by Dubrovin and Novikov [DN89, DN93]). Starting from a system of PDEs, such as (1), admitting families of periodic traveling waves, one may derive modulation equations through an averaging procedure, which yields PDEs on large space-time scales for the local parameters of the waves. The corresponding *ansatz*, expected to approximate solutions to the original PDEs, looks like one periodic wave train on small scales but have variable wave parameters on larger scales, hence exhibit varying amplitude and wavelength on these large scales.

A robust way to obtain them is to consider a two-scale formal asymptotic expansion combining slow arbitrary variables and single-phase fast oscillations,

$$\mathbf{U}(x,t) = \mathbf{U}_0\left(\varepsilon x, \varepsilon t, \frac{\phi_{(\varepsilon)}(\varepsilon x, \varepsilon t)}{\varepsilon}\right) + \varepsilon \mathbf{U}_1\left(\varepsilon x, \varepsilon t, \frac{\phi_{(\varepsilon)}(\varepsilon x, \varepsilon x)}{\varepsilon}\right) + \text{h.o.t}, \qquad (4)$$

with $0 < \varepsilon \ll 1$,

$$\phi_{(\varepsilon)}(X,T) = \phi_0(X,T) + \varepsilon \phi_1(X,T) + \text{h.o.t},$$
$$\mathbf{U}_j(X,T,\theta) \text{ one-periodic in } \theta, \quad j = 0, 1, \cdots,$$

where $X = \varepsilon x$ and $T = \varepsilon t$ denote some rescaled spatial and time variables respectively. A leading-order identification shows that for all $(X,T)$,

$$\underline{\mathbf{U}}(\xi) := \mathbf{U}_0(X, T; \xi k(X, T))$$

must be the profile of a periodic, traveling wave solution to the original system, here (1), of (spatial) period $\Xi(X,T) = 1/k(X,T)$ and speed $c(X,T) = -\omega(X,T)/k(X,T)$, with $k = \partial_X \phi_0$ and $\omega = \partial_T \phi_0$. This already leaves as a constraint $\partial_T k = \partial_X \omega$, an equation known as the *conservation of waves* equation.

The missing part of the time evolution is obtained from constraints for the resolution of the next step of the identification. Indeed this step yields an affine equation for $\mathbf{U}_1$, with

- linear part essentially[1] given by the linearization about $\underline{\mathbf{U}}$ of the original system, in the frame moving with speed $c$, acting on functions of $\xi$ with the same period as $\underline{\mathbf{U}}$;

- source terms depending only on $\mathbf{U}_0$, $\phi_0$ and $\phi_1$.

The possibility to solve this step is then equivalent to the orthogonality (for the $L^2$-scalar product in the $\xi$-variable) of source terms to the kernel of the adjoint of the linear operator, a constraint automatically satisfied by the $\phi_1$-part of source terms. Now it turns out that elements of the latter kernel are in correspondence with local conservations laws for the original systems. In the present case the conservative nature of (1), (2), (3) is directly linked to the presence in the kernel of the adjoint of the linearization of respectively constant functions, $\delta \mathcal{Q}[\underline{\mathbf{U}}]$ and $\delta \mathcal{H}[\underline{\mathbf{U}}]$. Orthogonality to those then yields time evolution equations for the averages of the quantities involved in local conservation laws. Note however that for traveling waves such as $\underline{\mathbf{U}}$, $\delta \mathcal{H}[\underline{\mathbf{U}}]$ is already a linear combination of $\delta \mathcal{Q}[\underline{\mathbf{U}}]$ and constants so that the local conservation law for the averaged of $\mathcal{H}$ does not generate a new independent equation but an *entropy* for the modulated system.

The upshot of the detailed process (for which we refer to [NR13, BGNR14]) is the modulated system

$$\partial_T k - \partial_X \omega = 0 \qquad (5a)$$
$$\partial_T \langle \underline{\mathbf{U}} \rangle - \partial_X (\langle \mathbf{B} \, \delta \mathcal{H}[\underline{\mathbf{U}}] \rangle) = 0 \qquad (5b)$$
$$\partial_T \langle \mathcal{Q}(\underline{\mathbf{U}}) \rangle - \partial_X \big(\langle \delta \mathcal{Q}[\underline{\mathbf{U}}] \cdot \mathbf{B} \delta \mathcal{H}[\underline{\mathbf{U}}] + \nabla_{\mathbf{U}_x} \mathcal{H}(\underline{\mathbf{U}}, \underline{\mathbf{U}}_\xi) \cdot \partial_\xi (\mathbf{B} \, \delta \mathcal{Q}[\underline{\mathbf{U}}]) - \mathcal{H}(\underline{\mathbf{U}}, \underline{\mathbf{U}}_\xi) \rangle\big) = 0 \qquad (5c)$$

where brackets $\langle \cdot \rangle$ stand for mean values over the period $\Xi(X,T)$. A few comments on dimensions are worth stating (and we refer to [Rod13] for a more detailed discussion). The original Hamiltonian system (1) counts $N$ equations whereas (5a)-(5b)-(5c) involves $N + 2$ equations. From the way the modulated system has been derived it should be clear that this $N + 2$ breaks into 1 for the number of wavenumber variables[2] and $N+1$ for the dimension of the span of $\delta\mathcal{F}[\underline{\mathbf{U}}]$, with $\mathcal{F}$ ranging along functionals that satisfy a conservation law along the flow of (1), $\underline{\mathbf{U}}$ being a periodic traveling wave. In the identification process, the number of independent averaged conservation laws, $N+1$ here, arises as related to the dimension of the kernel of the adjoint of the linearization in a moving frame, restricted to functions of the same period. This dimension must also be the dimension of the kernel of the adjoint of the linearized operator itself, hence the dimension of the family of periodic waves with fixed period and speed (associated with the abovementioned wavenumber). Thus the formal argument also carries the fact that the dimension of the modulated system ($N+2$ here) differs from the dimension of the family of periodic traveling waves ($N+3$ here) by the number of wavenumbers (1 here), hence agrees with the dimension of periodic wave profiles identified when coinciding up to translation (again $N+2$). The missing piece of information, about phase shifts, may be recovered a posteriori[3] by solving $\partial_T \phi_0 = \omega$, $\partial_X \phi_0 = k$. As proved for instance in [BGNR14, Appendix B], to a large extent, the present informal discussion on dimensions may be turned into sound mathematical arguments.

---

[1] Up to a rescaling from $\xi$ to $\theta$, normalizing period to 1.

As already pointed out in Whitham's seminal work [Whi99] for KdV and NLS, it is possible to use the variational structure of systems such as (1) to derive (5a)-(5b)-(5c) from least action considerations, instead of the geometrical optics expansion (4). For recent accounts of this variational derivation the reader is referred to [Kam00, Bri17]. As for the class of systems considered here, the corresponding form in terms of an action integral along periodic wave profiles was explicited in [BGNR13] and subsequently crucially used in [BGMR16, BGMR20].

A simple situation where one expects that the dynamics of (1) is well-described by a slow modulation *ansatz* similar to (4), hence obeys at leading-order a suitable version of (5a)-(5b)-(5c) is in the large-time regime starting from a smooth and localized perturbation of a single periodic traveling wave of (1), which should correspond to a nearly constant solution to (5a)-(5b)-(5c). Yet, though it is arguably the simplest relevant regime, a rigorous validation of the foregoing scenario has been obtained for none of the equations considered here. See however [BGNR14] for a spectral validation on the full class (1), [Rod18] for a linear validation on KdV, and [JNRZ13, JNRZ14] for full validations but on parabolic systems.

## Small dispersion limit

Our present contribution is rather focused in regimes involving solutions to (5a)-(5b)-(5c) covering the full range of possible amplitudes, and known as *dispersive shock waves*. These are typically expected to play a key role in the regularization of shocks by weak dispersion.

This regime may be described by introducing a small wavenumber parameter $\varepsilon > 0$ and moving to rescaled variables $(X, T) = (\varepsilon\, x, \varepsilon\, t)$. Looking first for a non oscillatory slow expansion (instead of (4))

$$\mathbf{U}(x,t) = \mathbf{U}_0\left(\varepsilon\, x, \varepsilon\, t\right) + \varepsilon \mathbf{U}_1\left(\varepsilon\, x, \varepsilon\, t\right) + \text{h.o.t}\,, \qquad (6)$$

suggests that $\mathbf{U}_0$ should satisfy

$$\partial_T \mathbf{U}_0 \,=\, \partial_X(\mathbf{B}\nabla_{\mathbf{U}}\mathcal{H}(\mathbf{U}_0, 0))\,, \qquad (7)$$

a first-order system of conservation laws. To make the discussion more concrete, let us temporarily focus on the KdV case, where (1) becomes in slow variables

$$\partial_T v \,+\, \partial_X(\tfrac{1}{2}v^2) \,+\, \varepsilon^2 \partial_X^3 v \,=\, 0$$

and (7) reduces to the Hopf — or inviscid Burgers — equation

$$\partial_T v \,+\, \partial_X(\tfrac{1}{2}v^2) \,=\, 0\,.$$

In the KdV case, for nonnegative initial data and as long as the Hopf equation does not develop a shock, Lax and Levermore proved in [LL83] that as $\varepsilon \to 0$ the solutions to the above scaled KdV equation starting with the same initial datum converge strongly in $L^2$ to the solution of the Hopf equation. They also proved that after the shock formation a weak limit still exists but it does not solve the Hopf equation almost everywhere anymore. Instead there coexist some zones where the weak limit of $v^2$ coincide with the square of the limit of $v$ and the latter solves the Hopf equation and other zones where this fails and the weak limit of $v$ does not satisfy an uncoupled scalar PDE but is the component $\langle v \rangle$ of a solution to the system of three equations given by (5a)-(5b)-(5c), specialized to KdV.

---

[2]The nonlinear Schrödinger equations (in original formulations) form typically a case with a two-dimensional group of symmetries. Their reduction to hydrodynamic form lowers the symmetry dimension by 1 but adds a conservation law.

[3]This requires some extra study though. See for instance the dedicated analysis in [JNRZ13, JNRZ14, Rod18].

Since the seminal [LL83] there has been a lot of attention devoted to weak dispersion limits and their link to modulated equations (including some multi-phase versions), but for the moment all mathematically complete analyses are restricted to the consideration of completely integrable PDEs, such as KdV, the modified Korteweg–de Vries equation, the cubic Schrödinger equations, the Benjamin-Ono equation and equations in their hierarchies. We refer the interested reader to [Ven85b, Ven85a, Ven87, Ven90, Wri93, DVZ97, JLM99, Gra02, Mil02, TVZ04, TVZ06, PT07, MX11, MX12, JM14, Jen15, MW16] for original papers and to [EJLM03, EH16, Mil16] for a detailed account of progresses made so far in this direction. Note however that, whatever the precise definition we use for completely integrable PDEs, these correspond to specific nonlinearities and in fact few models are completely integrable.

To unravel some of the reasons why the completely integrable case is significantly simpler to analyze — but still tremendously involved ! —, let us step back a little and draw some analogies with the *vanishing viscosity limit*. The natural parabolic counterpart to the KdV equation above is the viscous Burgers equation

$$\partial_T v \,+\, \partial_X(\tfrac{1}{2}v^2) \,=\, \varepsilon\, \partial_X^2 v\,.$$

In this case, by using the Hopf-Cole [Hop50, Col51] transform, precisely introduced for this purpose, it is relatively easy to analyze the limit $\varepsilon \to 0^+$ and check that solutions to the Burgers equation converge to the weak solution of the Hopf equation given by the Lax-Oleĭnik formula.

The convergence towards a weak solution of the inviscid equation and the characterization of the vanishing viscosity limits by an entropy criterion has been extended even beyond the scalar case [Ole57] to solutions to systems starting from initial data with small total variation [BB05]. At the heart of these general treatments lies an understanding not only of the limiting slow behavior encoded by the solution of the inviscid equation but also of the fast part essentially supported near discontinuities of the inviscid solution. To give a heuristic flavor of the latter, let us focus again on the scalar case and zoom from $(T,X)$ to $(T,\widetilde{x}) = (T,(X-\psi_{(\varepsilon)}(T))/\varepsilon)$, with $\psi_0(\cdot)$ describing the position of a discontinuity of the limiting solution $v_0$ and $X$ living in a neighborhood of the latter discontinuity. Then the identification of powers of $\varepsilon$ suggests that the fast local structure, at time $T$ near the discontinuity located at $\psi_0(T)$ should be described by a front of the Burgers equation (in fast variables) traveling with velocity $\partial_T \psi_0(T)$ (satisfying the Rankine-Hugoniot condition) and joining $v_0(\psi_0(T)^-,T)$ to $v_0(\psi_0(T)^+,T)$ (where $\pm$ denote limits from above or from below). The existence of such viscous fronts plays a deep role in both the heuristic and rigorous treatment of the vanishing viscosity limit. In particular, even the rigorous masterpiece by Bianchini and Bressan [BB05] proceeds through such a local traveling-wave decomposition of solutions.

In the dispersive case the possibility of this global-slow/local-fast scenario fails already in general by the non existence of the required traveling fronts. Indeed, elementary considerations show that whereas in the scalar diffusive case the set of pairs of values that may be joined by a nondegenerate front is an open subset of $\mathbb{R}^2$ and the selection of the speed coincides with the one given by the Rankine-Hugoniot condition, in the scalar KdV-like case[4] this set is a submanifold of dimension 1. At this stage it should be clear that the understanding of what are the possible fast structures replacing viscous fronts, also called *viscous shock waves*, to join both sides of a shock should already provide a good grasp on the weak dispersion limit.

## The Gurevich–Pitaevskii problem

Leaving aside the possibility that the fast part of solutions could be given by well-localized elementary blocks, steady in the frame moving with the speed of the shock they regularize and interacting with

---

[4]Note moreover that in the classical KdV case this set is empty since fronts require a potential with at least two local maxima.

the slow part only through their limiting values at infinity, one is naturally led to the consideration as elementary fast blocks of unsteady patterns as in (4), mixing slow and fast scales but whose limits at $\pm\infty$ are purely slow, that is, constant with respect to the fast variable. Modulated periodic wave trains may reach these limiting constant states in two ways

- by letting their amplitude go to zero, they reach a constant state by asymptoting a harmonic periodic wavetrain oscillating about the reference constant;

- by letting their wavelength go to infinity, they converge to a solitary wave connected to the reference constant by its limiting trail.

From the foregoing considerations arises the question of determining when two given constant states may be joined by a rarefaction wave of (5a)-(5b)-(5c) in the sense that limiting values of the rarefaction wave are parameters corresponding to either harmonic or solitary waves and the limiting values of the average part $\langle \underline{\mathbf{U}} \rangle$ fits the prescribed constants. The corresponding unsteady, oscillatory patterns, recovered through (4), are referred to as *dispersive shock waves*. Note that the question differs from the investigation of classical rarefaction waves of hyperbolic systems in at least two ways. On one hand, both harmonic wavetrains with a prescribed limiting value and solitary waves with a prescribed endstate form one-dimensional families (when identified up to translation) and the knowledge of through which harmonic train or solitary wave given constant states may be joined is an important part of the unknown elements to determine. On the other hand the modulated system (5a)-(5b)-(5c) is a priori not defined at the limiting values and yet the hope to match dispersive shock waves with solutions to (7) heavily relies on the expectation that in both limits — solitary or harmonic — (5b) uncouples from the rest of the system and converges to (7).

It is worth stressing that even though one expects to obtain, afterwards, a multi-scale pattern through (4) the foregoing problem is $\varepsilon$-free. It is a dispersive analogous to the determination of conditions under which two constants may be joined by a viscous front. Note that in the viscous version of the problem such considerations are then included in the definition of admissibility of weak solutions to (7), and expected to determine reachability by vanishing viscosity limits. Notably, in the classical Riemann problem, one considers how to solve (7) starting from an initial datum taking one value up to some point then another value, by gluing constants, rarefaction waves, admissible shocks and contact discontinuities. Solutions to the Riemann problems may then be used themselves as elementary blocks to solve the general Cauchy problem for (7) (designed from vanishing viscosity considerations). See for instance [Ser99, Bre00] for background on the classical Riemann problem. Likewise, in the zero dispersion limit, the Gurevich–Pitaevskii problem consists in joining two given constants on two complementary half-lines with constant sectors, rarefaction waves of (7) and rarefaction waves of (5a)-(5b)-(5c), the junction between solutions to (7) and solutions to (5a)-(5b)-(5c) being understood in the sense mentioned herein-above. This approach was introduced for KdV by Gurevich and Pitaevskii in [GP73] and has been referred to as the Gurevich–Pitaevskii problem since then, or sometimes the dispersive Riemann problem.

It must be stressed that as for Riemann problems in the weak dissipation limit, the Gurevich–Pitaevskii problems are expected to carry a wealth of information on the weak dispersion regime. We already pointed out that a fully rigorous treatment of the weak dispersion limit is for the moment restricted to some equations, associated with Lax pairs including a scalar Schrödinger operator and completely integrable through inverse scattering transforms. Unfortunately the same is true for the associated Gurevich–Pitaevskii problems. Indeed modulated systems of those particular systems inherit from the Lax pair representation of the original system, a family of strong Riemann invariants, given by edges of Lax spectral bands; see [DN74, DMN76, FFM80, DN83, FL86, Pav87] for original papers pointing this connection. The latter observation was certainly the main motivation for the introduction and the study of the classes of hyperbolic systems possessing a complete set of strong Riemann invariants, a class coined as rich by Serre [Ser00, Chapter 12] and as semi-Hamiltonian by Tsarev [Tsa85, Tsa90, Tsa00]. Along a rarefaction wave of a rich system all but one Riemann invariant are constants. Moreover in both the solitary wave limit and the harmonic limit of PDEs associated with such Lax pairs one of the Lax spectral gaps closes so that the Riemann invariant varying along a rarefaction wave of such a rich modulated system connecting two harmonic/solitary limits is actually merging in both limits with one of the steady Riemann invariants. This makes the rarefaction wave problem considerably simpler to solve, at least as far as determining which states may be connected and what are the trail and edge speeds of the pattern.

Given its particular importance for some classes of applications, there have been a few attempts to propose solutions to the Gurevich–Pitaevskii not relying on strong Riemann invariants. One of the most remarkable attempt is due to El and the reader is referred to [EH16] for details on the method and to [El05, Hoe14] for two instances of application. The elegant method of El provides answers consistent with rigorous analyses of integrable cases and displays reasonably good agreement with numerical experiments. Yet unfortunately, so far, it still lacks strong[5] theoretical support. Elucidating the mathematical validity of the approach of El may be considered as one of the key problems of the field.

To conclude the exploratory part of the paper, we point out that even from a formal point of view there are at least two important features of the weak dispersion limit that we have left aside and on which we comment now.

**Remark 1.** System (5a)-(5b)-(5c) is itself a — hopefully[6] hyperbolic — first-order system so that it may be expected to develop shocks in finite time and the expansion in (4) to suffer from a finite-time validity (in the slow variables) in the same way as the relevance of (6) stops when the corresponding solution to (7) forms a shock. Yet the formal process itself hints at an $\varepsilon$-dispersive correction to (5a)-(5b)-(5c) — see for instance [Rod18] —, so that the phenomena may be thought itself as a weak dispersion limit in the presence of wave-breaking at the level of (single-phase) modulation equations, suggesting the presence of oscillations at this level, and resulting in a two-phase oscillation pattern at the original level. For KdV a compatible scenario (with arbitrary number of phases) was proposed in [FFM80] within the terminology of integrability by inverse scattering ; it was subsequently recast in terms of averaged modulation equations in [EKV01] and proved to hold in [GT02, Gra04]. Note that the prediction includes a description of where 0-phase, 1-phase and 2-phases patterns live in the space-time diagram.

**Remark 2.** There has been considerable effort devoted to the description of what is seen on a zoom in a neighborhood of a wave-breaking point. This results in a different asymptotic regime and a suitable scenario was first proposed by Dubrovin [Dub06] and then proved for various integrable PDEs in [CG09, CG12, BT13].

## Structure of general modulated systems

As far as the formal description of dispersive shocks by means of modulated equations is concerned, multi-scale regions are connected to single-scale ones by either one of the two asymptotic regimes corresponding to the small amplitude limit - when the amplitude of the waves goes to zero - and the solitary wave limit - when the wavelength of the waves goes to infinity. The understanding of both these regimes is a crucial step towards the actual construction of dispersive shocks. In particular, to analyze rarefaction waves of modulated equations connecting such asymptotic regimes, one needs to elucidate the hyperbolic nature of its eigenfields in both distinguished limits. Indeed, as for the classical Riemann problem, the resolution of the Gurevich–Pitaevskii problem crucially relies on the hyperbolicity and the structure of the eigenfields of modulated equations. This requires a study not only of averaged quantities involved in the conservative formulation but also of their derivatives, as appearing in the expanded, quasilinear form.

With this respect, note that unfortunately, the formulation of modulated equations in terms of what is arguably the most natural set of wave parameters blows up in the solitary wave limit. This issue has been partially resolved by El [El05] by replacing one of the parameters with the so-called *conjugate*

---

[5]The method implicitly extrapolates rarefaction curves in the bulk $k > 0$, $\alpha > 0$ from computations purely carried out in limiting hyperplanes $k = 0$ and $\alpha = 0$, a process that seems hard to justify without assuming the modulation system under consideration to be rich/semi-Hamiltonian, whereas this property is known to be satisfied for almost no hyperbolic system besides modulation systems of integrable equations.

[6]Since otherwise the modulation system is ill-posed, (some of) the underlying periodic waves are unstable and the scenario in (4) fails instantly.

*wavenumber*. However, this new parametrization is in general limited to the large wavelength regime[7].

One of our main contributions here is to provide a global set of parameters. For the latter, we prove in great generality that

- it may serve as a parametrization of periodic wave profiles (identified up to spatial translation) exactly when the original averaged quantities $(k, \langle \underline{\mathbf{U}} \rangle, \langle \mathcal{Q}(\underline{\mathbf{U}}) \rangle)$ can, that is, as proved in [BGNR14, Appendix B] and [BGMR16, Theorem 1], exactly when wave profiles with fixed period form an $(N+1)$-dimensional manifold (when identified up to spatial translation);

- in these variables the modulated system possesses an Hamiltonian formulation, with Hamiltonian function the original averaged Hamiltonian energy (see Theorem 1);

- these variables may be extended to solitary-wave and harmonic limits in such a way that the modulated system admits regular limits even in quasilinear form (see Theorem 4).

The proposed system of coordinates already appeared for the Euler–Korteweg system in mass Lagrangian coordinates in [GS95] (also see [BG13] for an account of Gavrilyuk–Serre's result with our notation), but its significance remained unclear. As we show hereafter, it turns out to apply to our more general framework, and to give new insight on modulated systems. The point is to replace the conserved variable $\langle \mathcal{Q}(\underline{\mathbf{U}}) \rangle$ in Equations (5a)-(5b)-(5c) by another one, denoted merely by $\alpha$ hereafter. This new variable tends to zero when the amplitude of the wave tends to zero — as the amplitude squared, as we shall see later on —, and has a finite limit when $k$ goes to zero, that is, in the solitary-wave limit. Remarkably enough, this quantity can be defined as simply as

$$\alpha := \frac{1}{k} \left( \langle \mathcal{Q}(\underline{\mathbf{U}}) \rangle - \mathcal{Q}(\langle \underline{\mathbf{U}} \rangle) \right) . \tag{8}$$

It turns out that, as far as smooth solutions are concerned, the modulation equations take the alternative form

$$\partial_T \begin{pmatrix} k \\ \alpha \\ \mathbf{M} \end{pmatrix} = \mathbb{B} \, \partial_X \left( \nabla_{k, \alpha, \mathbf{M}} \mathrm{H} \right) , \tag{9}$$

where

$$\mathbb{B} := \left( \begin{array}{cc|c} 0 & 1 & \\ 1 & 0 & \\ \hline & & \mathbf{B} \end{array} \right)$$

is symmetric and nonsingular and

$$\mathbf{M} := \langle \underline{\mathbf{U}} \rangle , \qquad \mathrm{H} := \langle \mathcal{H}[\underline{\mathbf{U}}] \rangle .$$

The Hamiltonian structure of System (9) provides a form of symmetry in the spirit of Godunov's theory of hyperbolic systems. Nevertheless, this form does not automatically provide energy estimates nor imply hyperbolicity because, as our expansions show (see Remark 14), the associated potential, natural symmetrizer is not definite in either one of the limits.

System (9) has also an appealing symmetric form with respect to the distinguished limits, $k \to 0$ corresponding to the long wavelength limit and $\alpha \to 0$ to the small amplitude limit. Yet another upshot of our analysis is a strong, somewhat surprising asymmetry as regards the asymptotic nature of the eigenfields. In the solitary wave regime, the hyperbolicity of the modulated equations is equivalent to its weak hyperbolicity and may persist even at the limit in the presence of the solitary wave speed as a double root. We see this striking property as a consequence of the strong separation of scales displayed in asymptotic expansions of large wavelength profiles (see Remark 17). By contrast, in the harmonic wave regime, in general the hyperbolicity of modulated equations is lost at the limit, the characteristic speed corresponding to the group velocity being non semi-simple — associated with a Jordan block of height two. See Theorems 5 and 6.

---

[7]To see this on a concrete example, consider the defocusing modified KdV equation. In this case potentials are coercive quartic polynomials, hence in general possess two wells. Except in the very exceptional case when the two wells have exactly the same depth, the conjugated wavenumber parametrization necessarily breaks down before the small amplitude limit for one of the two families of periodic waves associated with each potential well.

We stress that many asymptotic properties of the modulation systems are much easier to study once a limiting system has been identified. This is precisely where we benefit from our new set of parameters. In particular, once Theorems 1 and 4 are known it is relatively easy to derive the most basic properties of the modulation systems, both mentioned and used in the preceding formal discussion of dispersive shocks. For instance Corollary 5 contains that at both limits the modulation system split into a block coinciding with the original dispersionless system, System (7), and a $2 \times 2$ block with double characteristic given by either the solitary wave speed at the long wavelength limit or with the harmonic linear group velocity at the small amplitude limit. Yet, even with good variables in hands, some further properties require finer details of higher-order expansions.

Our new set of parameters enables us to show how the eigenfields of modulated equations degenerate in the small amplitude and the large wavelength limits (see Theorems 7 and 8) but it relies on a more involve analysis. The main upshots are that

- Near the harmonic limit, we derive explicit conditions determining modulational instability (see Theorem 7 and Appendix A), those being known in some cases as the Benjamin–Feir criteria.

- Near the soliton limit, we prove that modulational instability is determined by exactly the same condition ruling stability of solitary waves and, as proved in [BGMR20], co-periodic stability of nearby periodic waves, that is, it is decided by the sign of the second derivative — with respect to speed, fixing the endstate — of the Boussinesq moment of instability.

For the conclusions mentioned here it should be emphasized that it is relatively easy to support wrong deceptive claims when arguing formally; see for instance Remarks 15 and A.i. Another somewhat surprising, but not unrelated (see Remark 17), discrepancy between both limits is that the convergence of eigenvalues splitting from the double root is exponentially fast in the solitary limit.

We insist on the fact that, surprising or not, the properties of the modulation systems discussed in the present contribution are *proved* here for the first time for a class of systems not restricted to integrable systems. Nevertheless, for the satisfaction of readers interested in knowing which part of our conclusions could be intuited from sound formal arguments, we have singled out in Section 4.2 conclusions that could be inferred from the knowledge that both variables $(k, \alpha, \mathbf{M})$ and averaged Hamiltonian H extend with sufficient smoothness to both limits. We stress however that in the solitary wave limit these extensions are highly nontrivial to prove.

We conclude this introductory section with a few words on the nature of proofs contained in the rest of the paper. The most elementary ones are purely algebraic manipulations. For the other ones we rely on asymptotic expansions of the abbreviated action integral of the profile ODE, and of its derivatives up to second order. These were derived in detail in [BGMR20] and used there to deduce some consequences on the stability of periodic waves with respect to co-periodic perturbations. As we show in the present paper, that asymptotic behavior gives insight on the modulated equations as well.

The essence of some of the results — providing, for a large class of Hamiltonian partial differential equations, connections between various formulations and stability indices mixing many different coefficients — and of their proofs — involving two kinds of asymptotic expansions for an $N + 2$-dimensional matrix, mixing four different orders in the soliton regime, — leads to an unavoidable but sometimes overwhelming notational inflation. To counterbalance this and help the reader finding his/her way in the notational maze, we have included in Appendix C a symbolic index and gathered some of the explicit formulas scattered through the text, specialized to the generalized Korteweg–de Vries equation.

We hope that this also provides sufficient support for the reader interested in derived results and take-home messages but not in their lengthy proofs so as to jump through the text from key result to key result and consult the table of symbols when needed. For such a reader, here is a possible roadmap:

- Theorem 3, p. 27, proves that variables $(k, \alpha, \mathbf{M})$ extend up to asymptotic boundaries without loss of dimension;

- Theorem 4, p. 29, proves that the averaged Hamiltonian H extends up to asymptotic boundaries as a $\mathcal{C}^2$ function;

- Theorem 7, p. 38, provides expansions of eigenvalues and eigenvectors of the modulated system in the harmonic limit;

- Theorem 8, p. 40, provides expansions of eigenvalues and eigenvectors of the modulated system in the soliton limit ;

- Appendix A details explicit formula for the small-amplitude instability index in both the scalar ((58), p. 43) and the system ((59), p. 48) cases.

The general setting and various formulations of modulated equations are presented in Section 2. Asymptotic properties of the alternate parametrization are established in Section 3. Limits of the modulated system are derived in Section 4. Eigenfields are studied in Section 5. Appendix A contains explicit modulational instability criteria for the harmonic limit. Appendix B provides some insights on the splitting of double roots. Appendix C is the symbolic index.

**Acknowledgement.** The first and third author would like to express their gratitude to Gennady El, Sergey Gavrilyuk, Mark Hoefer and Michael Schearer for enlightening discussions during the preparation of the present paper. They also thank an anonymous referee for the suggestion to work out the possible connection with the notion of wave action, that has resulted in Remark 3.

**Matrix notation.** Along the text, in matrices, 0 may denote scalar, vector or matrix-valued zeroes. Moreover empty entries denote zeroes and $*$ entries denote values that are irrelevant for the discussion and may vary from line to line.

## 2 Various formulations of modulated equations

### 2.1 General framework

As announced in the introduction, we consider abstract systems of the form

$$\partial_t \mathbf{U} = \partial_x (\mathbf{B}\, \delta\mathcal{H}[\mathbf{U}]). \tag{10}$$

where the unknown $\mathbf{U}$ takes values in $\mathbb{R}^N$, $\mathbf{B}$ is a symmetric and nonsingular matrix, and $\delta\mathcal{H}[\mathbf{U}]$ denotes the variational derivative at $\mathbf{U}$ of $\mathcal{H} = \mathcal{H}(\mathbf{U}, \mathbf{U}_x)$. For the sake of clarity, here, we shall mostly stick to Assumption 1 below, all the more so when we are to apply results from [BGMR20].

**Assumption 1.** *There are smooth functions $f$, $\kappa$ and $\tau$ with $\kappa$ and $\tau$ taking only positive values, and a nonzero real number $b$ such that*

- *either $N = 1$, $\mathbf{U} = v$, $\mathcal{H} = e(v, v_x)$, and $\mathbf{B} = b$,*

- *or $N = 2$, $\mathbf{U} = (v, u)^\mathsf{T}$,*

$$\mathcal{H} = \frac{1}{2}\tau(v)\, u^2 + e(v, v_x), \qquad \mathbf{B} = \begin{pmatrix} 0 & b \\ b & 0 \end{pmatrix},$$

*with*

$$e(v, v_x) = \frac{1}{2}\kappa(v)\, v_x^2 + f(v)$$

*in both cases.*

By definition we have

- in the case $N = 1$, $\delta\mathcal{H}[\mathbf{U}] = \delta e[v] := \partial_v e(v, v_x) - \partial_x(\partial_{v_x} e(v, v_x))$,
- in the case $N = 2$,

$$\delta\mathcal{H}[\mathbf{U}] = \begin{pmatrix} \frac{1}{2}\tau'(v)\, u^2 + \delta e[v] \\ \tau(v)\, u \end{pmatrix}.$$

The *impulse*

$$\mathcal{Q}(\mathbf{U}) := \frac{1}{2}\mathbf{U} \cdot \mathbf{B}^{-1}\mathbf{U},$$

generates $x$-translations in that

$$\partial_x \mathbf{U} = \partial_x(\mathbf{B}\,\delta\mathcal{Q}[\mathbf{U}]).$$

From the invariance of $\mathcal{H}(\mathbf{U}, \mathbf{U}_x)$ with respect to $x$-translations, that reads in differentiated form

$$\partial_x(\mathcal{H}(\mathbf{U}, \mathbf{U}_x)) = \delta\mathcal{H}[\mathbf{U}] \cdot \partial_x \mathbf{U} + \partial_x(\mathbf{U}_x \cdot \nabla_{\mathbf{U}_x}\mathcal{H}(\mathbf{U}, \mathbf{U}_x)),$$

stems the local conservation law

$$\partial_t \mathcal{Q}(\mathbf{U}) = \partial_x(\mathbf{U} \cdot \delta\mathcal{H}[\mathbf{U}] + \mathsf{L}\mathcal{H}[\mathbf{U}]), \tag{11}$$

along smooth solutions of (10), where

$$\mathsf{L}\mathcal{H}[\mathbf{U}] := \mathbf{U}_x \cdot \nabla_{\mathbf{U}_x}\mathcal{H}(\mathbf{U}, \mathbf{U}_x) - \mathcal{H}(\mathbf{U}, \mathbf{U}_x) = v_x\,\partial_{v_x} e(v, v_x) - \mathcal{H}(\mathbf{U}, \mathbf{U}_x).$$

The modulated system (5a)-(5b)-(5c) is also written as

$$\partial_T k - \partial_X \omega = 0, \tag{12}$$
$$\partial_T \langle \underline{\mathbf{U}} \rangle - \partial_X \langle \mathbf{B}\,\delta\mathcal{H}[\underline{\mathbf{U}}] \rangle = 0, \tag{13}$$
$$\partial_T \langle \mathcal{Q}(\underline{\mathbf{U}}) \rangle - \partial_X \langle \underline{\mathbf{U}} \cdot \delta\mathcal{H}[\underline{\mathbf{U}}] + \mathsf{L}\mathcal{H}[\underline{\mathbf{U}}] \rangle = 0, \tag{14}$$

where $\xi \mapsto \underline{\mathbf{U}}(X, T; \xi)$ is the profile of a periodic, traveling wave solution to (10) of (spatial) period $\Xi(X, T) = 1/k(X, T)$ and speed $c(X, T) = -\omega(X, T)/k(X, T)$, and brackets $\langle \cdot \rangle$ stand for mean values over the period $\Xi(X, T)$. Again we refer to [BGNR14] for a formal derivation of the system from a geometrical optics expansion.

As proved in [BGNR14, Appendix B] and [BGMR16, Theorem 1], the fact that System (5a)-(5b)-(5c) is a closed system, of evolution type, for initial data under consideration, is equivalent to the fact that, for each fixed period, periodic wave profiles under consideration form a non-degenerate manifold of dimension $N + 1$ when identified up to translation. In this case wave profiles may be smoothly parametrized by $(k, \langle \underline{\mathbf{U}} \rangle, \langle \mathcal{Q}(\underline{\mathbf{U}}) \rangle)$. As mentioned in the introduction, from the point of view of modulation theory, the range of validity of the latter parametrization is optimal. Yet these coordinates come with at least three serious drawbacks:

- they are not very explicit so that within this set of coordinates the modulation system is hard to manipulate;
- they are degenerate in the solitary-wave limit, losing two dimensions instead of one dimension;
- they do not provide a clear variational form for the modulation system.

We first recall how the first and third issues may be fixed by choosing a parametrization involving constants of integration of the wave profile ODEs.

**Remark 3.** Concerning the third point, let us make explicit that we focus mostly on variational structures of Hamiltonian type — or more precisely on Poisson structures. It is equally common to study Lagrangian formalisms. For the reader willing to work out Lagrangian counterparts we note that the starting point is that (1) is the Euler–Lagrange equation for the Lagrangian density

$$\mathfrak{L}[\boldsymbol{\psi}] := \frac{1}{2}\boldsymbol{\psi}_x \cdot \mathbf{B}^{-1}\boldsymbol{\psi}_t - \mathcal{H}[\boldsymbol{\psi}_x]$$

with $\mathbf{U} = \boldsymbol{\psi}_x$. Incidentally we stress that

$$\mathcal{Q}(\boldsymbol{\psi}_x) = \nabla_{\boldsymbol{\psi}_t}\mathfrak{L}[\boldsymbol{\psi}] \cdot \boldsymbol{\psi}_x$$

so that $\langle \mathcal{Q}(\mathbf{U}) \rangle$ may be identified with the *wave action*, a quantity that plays a prominent role in the classical theory of wave motion ; see [Hay70] and [Whi99, Section 11.7]. To go a bit further we point out the traveling-wave form

$$\boldsymbol{\psi} : \ (t,x) \mapsto \underline{\boldsymbol{\psi}}(k\,x + \omega\,t)$$

with $\underline{\mathbf{U}} = k\underline{\boldsymbol{\psi}}'(k\,\cdot)$, which yields the following formula for the averaged Lagrangian on wave profiles

$$\mathfrak{L}^{\mathrm{av}} = \omega\,\langle \mathcal{Q}(\underline{\mathbf{U}}) \rangle - \langle \mathcal{H}(\underline{\mathbf{U}}, k\underline{\mathbf{U}}_x) \rangle\,.$$

In particular when a choice of wave parameters including $(\omega, k)$ is done, say $(\omega, k, \mathbf{A})$ for some $\mathbf{A}$, by using the profile equation one may readily check that

$$\partial_\omega \mathfrak{L}^{\mathrm{av}} = \langle \mathcal{Q}(\underline{\mathbf{U}}) \rangle\,, \qquad \partial_k \mathfrak{L}^{\mathrm{av}} = -k\,\Theta\,, \qquad \nabla_{\mathbf{A}} \mathfrak{L}^{\mathrm{av}} = 0\,,$$

where $\Theta$ denotes the abbreviated action integral introduced below.

## 2.2 Modulated equations in terms of constants of integration

For a traveling wave $\mathbf{U} = \underline{\mathbf{U}}(x - ct)$ of speed $c$ to be solution to (10), there must exist a $\boldsymbol{\lambda} \in \mathbb{R}^N$ such that

$$\delta(\mathcal{H} + c\mathcal{Q})[\underline{\mathbf{U}}] + \boldsymbol{\lambda} = 0\,, \tag{15}$$

which can be viewed as the *Euler–Lagrange equation* $\delta\mathcal{L}[\underline{\mathbf{U}}] = 0$ associated with the *Lagrangian*

$$\mathcal{L} = \mathcal{L}(\mathbf{U}, \mathbf{U}_x; c, \boldsymbol{\lambda}) := \mathcal{H}(\mathbf{U}, \mathbf{U}_x) + c\mathcal{Q}(\mathbf{U}) + \boldsymbol{\lambda} \cdot \mathbf{U}\,.$$

Exactly as (3) was derived from (1) and the invariance by translation of $\mathcal{H}(\mathbf{U}, \mathbf{U}_x)$, it stems from the translational invariance of $\mathcal{L}(\mathbf{U}, \mathbf{U}_x)$ that solutions to (15) are such that for some $\mu \in \mathbb{R}$

$$\mathsf{L}\mathcal{L}[\underline{\mathbf{U}}] = \mu\,, \tag{16}$$

where

$$\mathsf{L}\mathcal{L}[\mathbf{U}] := \mathbf{U}_x \cdot \nabla_{\mathbf{U}_x}\mathcal{L}(\mathbf{U}, \mathbf{U}_x) - \mathcal{L}(\mathbf{U}, \mathbf{U}_x)\,.$$

This more concrete point of view introduces as wave parameters their speed $c$ and the integration constants $\boldsymbol{\lambda} \in \mathbb{R}^N$ and $\mu \in \mathbb{R}$. As discussed with more details in [BGNR13, BGMR16] — see in particular the proof of [BGMR16, Theorem 1] —, once these parameters are fixed, non-constant periodic wave profiles form a discrete set and the corresponding profiles perturb smoothly with respect to parameters. Henceforth we will often omit to specify that one of the branches of periodic waves have been chosen. By doing so, we obtain (a discrete number of) natural parametrizations of wave profiles.

As already pointed out in [BGNR13, BGMR16, BGMR20], many key properties of periodic traveling waves are more neatly stated in terms of the wave-speed and constants of integration by introducing the abbreviated action integral

$$\Theta(\mu, c, \boldsymbol{\lambda}) := \int_0^{\Xi} (\mathcal{H}[\underline{\mathbf{U}}] + c\mathcal{Q}(\underline{\mathbf{U}}) + \boldsymbol{\lambda} \cdot \underline{\mathbf{U}} + \mu)\,\mathrm{d}\xi\,, \tag{17}$$

whose definition involves a periodic profile $\underline{\mathbf{U}}$ of fundamental period $\Xi$ corresponding to parameter values $(\mu, c, \boldsymbol{\lambda})$. The action provides a nice closed form of the modulated equations in (12)-(13)-(14) and it encodes the duality between constants of integration and averaged quantities. Indeed, let us recall from [BGNR13] the following.

**Proposition 1.** *Under assumption 1, consider $\Omega$ an open subset of $\mathbb{R}^{N+2}$ and*

$$(\mu, c, \boldsymbol{\lambda}) \in \Omega \mapsto (\underline{\mathbf{U}}, \Xi) \in L^\infty(\mathbb{R}) \times (0, +\infty)$$

*a smooth mapping[8] such that for each value of the parameters $(\mu, c, \boldsymbol{\lambda})$, the function $\underline{\mathbf{U}} = (\underline{v}, \underline{u})$ is a smooth, non-constant periodic solution to (15)-(16), and $\Xi$ is the fundamental period of $\underline{\mathbf{U}}$.*

*Then the function $\Theta$ defined in (17) is also smooth on $\Omega$, and such that*

$$\partial_\mu \Theta = \Xi, \quad \partial_c \Theta = \int_0^\Xi \mathcal{Q}(\underline{\mathbf{U}}) \, dx, \quad \nabla_{\boldsymbol{\lambda}} \Theta = \int_0^\Xi \underline{\mathbf{U}} \, dx. \tag{18}$$

**Corollary 1.** *In the framework of Proposition 1,*

1. *the system in (12)-(13)-(14) equivalently reads, as far as smooth solutions are concerned,*

$$(\partial_T + c \, \partial_X)(\nabla_{\mu,c,\boldsymbol{\lambda}} \Theta) = \partial_\mu \Theta \begin{pmatrix} \begin{array}{cc|c} 0 & 1 & \\ 1 & 0 & \\ \hline & & -\mathbf{B} \end{array} \end{pmatrix} \partial_X \begin{pmatrix} \mu \\ c \\ \boldsymbol{\lambda} \end{pmatrix}; \tag{19}$$

2. *the mapping*

$$(\mu, c, \boldsymbol{\lambda}) \in \Omega \mapsto \left( k = 1/\Xi, \langle \underline{\mathbf{U}} \rangle = \tfrac{1}{\Xi} \int_0^\Xi \underline{\mathbf{U}} \, dx, \langle \mathcal{Q}(\underline{\mathbf{U}}) \rangle = \tfrac{1}{\Xi} \int_0^\Xi \mathcal{Q}(\underline{\mathbf{U}}) \, dx \right)$$

*is a local diffeomorphism if and only if*

$$\det \left( \nabla^2_{\mu,c,\boldsymbol{\lambda}} \Theta(\mu, c, \boldsymbol{\lambda}) \right) \neq 0, \quad \forall (\mu, c, \boldsymbol{\lambda}) \in \Omega.$$

**Remark 4.** The simple, closed form in (19) of modulated equations is well-known for the KdV equation. It is for instance pointed out by Kamchatnov [Kam00, eq. (3.135)], who says that 'despite the simple appearance of these equations, they are not very useful in practice'.

Exactly as pointed out in the introduction for System (9), the symmetric form of (19) does not readily imply that this system is hyperbolic. This would be the case if the matrix $\nabla^2_{\mu,c,\boldsymbol{\lambda}} \Theta$ were definite, which is not the case in general. As was shown in [BGMR20, Corollaries 1 and 2], in non-degenerate cases, $\nabla^2_{\mu,c,\boldsymbol{\lambda}} \Theta$ has a negative signature — or Morse index — equal to $N$ for small amplitude waves and equal either to $N$ or to $N+1$ for those of large wavelength. In addition, as follows from [BGMR16], for $N = 1$ a definite Hessian matrix $\nabla^2_{\mu,c,\boldsymbol{\lambda}} \Theta$ is incompatible with the spectral stability of the underlying periodic wave.

An important drawback of the formulation of modulated equations in the variables $(\mu, c, \boldsymbol{\lambda})$ is that all the quantities appearing in the time derivatives in (19) blow up in the solitary wave limit. Indeed, $\partial_\mu \Theta = \Xi = 1/k$ goes to infinity when $k$ goes to zero, as well as $\nabla_{\boldsymbol{\lambda}} \Theta = \Xi \langle \underline{\mathbf{U}} \rangle$ and $\partial_c \Theta = \Xi \langle \mathcal{Q}(\underline{\mathbf{U}}) \rangle$.

## 2.3 An important averaged variable

We claim that, despite their complicated and implicit form, Equations in (12)-(13)-(14) admit an equivalent form in a system of coordinates that is rather well suited for the study of dispersive shocks, in that it allows to take both the small amplitude limit and the solitary wave limit ($k \to 0$), in a most symmetric manner. We achieve this goal by replacing the conserved variable $\langle \mathcal{Q}(\underline{\mathbf{U}})\rangle$ with another one, which we merely denote by $\alpha$, and that is given by

$$\alpha := \frac{1}{k}\left(\langle \mathcal{Q}(\underline{\mathbf{U}})\rangle - \mathcal{Q}(\langle \underline{\mathbf{U}}\rangle)\right) \;=\; \frac{1}{k}\langle \mathcal{Q}(\underline{\mathbf{U}} - \langle \underline{\mathbf{U}}\rangle)\rangle\,.$$

---
[8]That is, we choose one branch of waves.

As already pointed out above and proved below, this new variable tends to zero when the amplitude of $\underline{\mathbf{U}}$ tends to zero and has a finite limit when $k$ goes to zero, that is when $\underline{\mathbf{U}}$ eventually becomes a solitary wave profile. Indeed, first we observe that $\alpha$ goes indeed to zero when the amplitude of $\underline{\mathbf{U}}$ goes to zero, because $k$ goes to the nonzero *harmonic* wavenumber and $\langle \mathcal{Q}(\underline{\mathbf{U}})\rangle$ and $\mathcal{Q}(\langle \underline{\mathbf{U}}\rangle)$ both go to the value of $\mathcal{Q}$ at the constant limiting state of the small amplitude wave. As to the limit when $k$ goes to zero, we see that

$$\begin{aligned}
\alpha &= \int_{-\Xi/2}^{\Xi/2} \left(\mathcal{Q}(\underline{\mathbf{U}}(\xi)) - \mathcal{Q}(\langle \underline{\mathbf{U}}\rangle)\right) \mathrm{d}\xi \\
&= \int_{-\Xi/2}^{\Xi/2} \left(\mathcal{Q}(\underline{\mathbf{U}}(\xi)) - \mathcal{Q}(\langle \underline{\mathbf{U}}\rangle) - \nabla_{\mathbf{U}}\mathcal{Q}(\langle \underline{\mathbf{U}}\rangle)\cdot(\underline{\mathbf{U}}(\xi) - \langle \underline{\mathbf{U}}\rangle)\right) \mathrm{d}\xi \\
&\to \int_{-\infty}^{+\infty} \left(\mathcal{Q}(\underline{\mathbf{U}}^s(\xi)) - \mathcal{Q}(\mathbf{U}_s) - \nabla_{\mathbf{U}}\mathcal{Q}(\mathbf{U}_s)\cdot(\underline{\mathbf{U}}^s(\xi) - \mathbf{U}_s)\right) \mathrm{d}\xi
\end{aligned}$$

when $\Xi$ goes to infinity, where $\underline{\mathbf{U}}^s$ denotes the limiting, solitary wave profile, homoclinic to $\mathbf{U}_s$, the limit of $\langle \underline{\mathbf{U}}\rangle$. The asymptotic behavior of $\alpha$ in these limits is proved in more details in Subsection 3.2.

Another remarkable property of $\alpha$ is that, at least for our main models of interest, scalar equations ($N=1$) and Euler–Korteweg systems, one may determine the sign of $\alpha$ in terms of parameters governing the traveling profiles.

**Proposition 2.** *Under Assumption 1,*

- *if $N = 1$ then $\alpha$ has the sign of $b$;*
- *if $N = 2$ and $\tau = \mathrm{Id}$, then $\alpha$ has the sign of $b\lambda_2$;*
- *if $N = 2$ and $\tau \equiv 1$, then $\alpha$ has the sign of $-c$.*

*Proof.* All signs are ultimately deduced from the Cauchy–Schwarz inequality: when $F$ and $G$ are linearly independent, $\langle F\,G\rangle^2 < \langle F^2\rangle\,\langle G^2\rangle$. The simplest case is for scalar equations, for which $\mathcal{Q}(v) = v^2/(2b)$, so that

$$2\,b\,\alpha = \langle \underline{v}^2\rangle \times \langle 1\rangle - \langle \underline{v}\rangle^2 > 0$$

since $\underline{v}$ is not constant.

If $N=2$ then $\mathcal{Q}(v,u) = v\,u/b$. Yet, when $\tau \equiv 1$, it follows from (15) that $\underline{u} + (c\underline{v})/b$ is constant so that

$$\langle \mathcal{Q}(\underline{v},\underline{u})\rangle - \mathcal{Q}(\langle \underline{v}\rangle, \langle \underline{u}\rangle) = -\frac{c}{b^2}(\langle \underline{v}^2\rangle \times \langle 1\rangle - \langle \underline{v}\rangle^2)$$

is of the same sign as $-c$. Indeed from the relation between $\underline{u}$ and $\underline{v}$ it follows that $\underline{v}$ is not constant since $\underline{\mathbf{U}}$ is not constant. Likewise when $\tau = \mathrm{Id}$, $\underline{u} + \lambda_2/\underline{v}$ is constant and

$$\langle \mathcal{Q}(\underline{v},\underline{u})\rangle - \mathcal{Q}(\langle \underline{v}\rangle, \langle \underline{u}\rangle) = \frac{\lambda_2}{b}\left(\langle \underline{v}\rangle\,\left\langle \frac{1}{\underline{v}}\right\rangle - \langle 1\rangle^2\right)$$

is of the sign of $\lambda_2/b$ since, as $\underline{v}$ is non-constant, $\sqrt{\underline{v}}$ and $1/\sqrt{\underline{v}}$ are independent. $\square$

**Remark 5.** The case when $N=2$ and $\tau = \mathrm{Id}$ includes Eulerian formulations of the Euler–Korteweg systems, whereas the case when $N=2$ and $\tau \equiv 1$ encompasses mass Lagrangian formulations of such systems. Each element of the latter class is conjugated to an element of the former and vice versa. As pointed out in [BG13], correspondences respect traveling wave types, and, as proved in [BGNR14], they also respect details of (the nonzero part of) the spectrum of linearizations about periodic waves. Obviously the foregoing proposition is consistent with corresponding conjugacies. Indeed, denoting with subscripts $_E$ and $_L$ quantities corresponding to each formulation, it follows from [BGNR14] that $b_E = -1$, $b_L = 1$, $(\lambda_2)_E = c_L/b_L$ and

$$\frac{\alpha_E}{k_E} = \frac{\alpha_L}{k_L}.$$

## 2.4 Alternative form of modulated equations

Returning to our general framework, we claim that variables $(k, \alpha, \mathbf{M} := \langle \underline{\mathbf{U}} \rangle)$ may be used exactly when $(k, \mathbf{M}, \langle \mathcal{Q}(\underline{\mathbf{U}}) \rangle)$ may be used and that using the former yields an alternate formulation of the modulated equations (19) that still has a nice symmetric-looking structure, and is now well-suited for both the small amplitude limit ($\alpha \to 0$) and the solitary wave limit ($k \to 0$).

To begin with, note that the vector $(k, \alpha, \mathbf{M})$ is deduced from $(k, \mathbf{M}, \langle \mathcal{Q}(\underline{\mathbf{U}}) \rangle)$ through the map $(k, \mathbf{M}, \mathbf{P}) \mapsto (k, (\mathbf{P} - \mathcal{Q}(\mathbf{M}))/k, \mathbf{M})$, which is obviously a (local) diffeomorphism so that parametrizations are indeed equivalent. In particular, Corollary 1 provides a characterization of when parametrization by $(k, \alpha, \mathbf{M})$ is possible.

Now we provide counterparts to Proposition 1 and Corollary 1 for variables $(k, \alpha, \mathbf{M})$. Here, the role of $\Theta$ in (19) is to some extent played by the averaged Hamiltonian

$$\mathrm{H} := \langle \mathcal{H}[\underline{\mathbf{U}}] \rangle.$$

**Remark 6.** Remarkably the action integral $\Theta$ and the averaged Hamiltonian $\mathrm{H}$ are closely related. It follows indeed from the definition of $\Theta$ in (17) and the expression of its derivatives in (18) that

$$\Theta = \Xi\,\mathrm{H} + c\,\partial_c \Theta + \boldsymbol{\lambda} \cdot \nabla_{\boldsymbol{\lambda}} \Theta + \mu\,\partial_\mu \Theta.$$

Would $\Theta$ be strictly convex, we would recognize

$$-\Xi\,\mathrm{H} = c\,\partial_c \Theta + \boldsymbol{\lambda} \cdot \nabla_{\boldsymbol{\lambda}} \Theta + \mu\,\partial_\mu \Theta - \Theta \tag{20}$$

as being the *conjugate* function of $\Theta$.

**Theorem 1.** *In the framework of Proposition 1, assume that the action $\Theta$ defined in (17) has a nonsingular Hessian at all points in $\Omega$. Then the mapping $(\mu, c, \boldsymbol{\lambda}) \mapsto (k, \alpha, \mathbf{M})$ defined by*

$$k = \frac{1}{\Xi}, \qquad \alpha = \frac{1}{k}\left( \langle \mathcal{Q}(\underline{\mathbf{U}}) \rangle - \mathcal{Q}(\langle \underline{\mathbf{U}} \rangle) \right), \qquad \mathbf{M} = \langle \underline{\mathbf{U}} \rangle,$$

*is a local diffeomorphism. Moreover*

1. *as a function of $(k, \alpha, \mathbf{M})$ the averaged Hamiltonian*

$$\mathrm{H} := \langle \mathcal{H}[\underline{\mathbf{U}}] \rangle$$

   *is such that,*
$$\partial_k \mathrm{H} = \Theta - \alpha c, \qquad \partial_\alpha \mathrm{H} = -kc, \qquad \nabla_{\mathbf{M}} \mathrm{H} = \langle \delta \mathcal{H}[\underline{\mathbf{U}}] \rangle; \tag{21}$$

2. *the modulated equations in (12)-(13)-(14) — or equivalently (19) — have a closed form in the variables $(k, \alpha, \mathbf{M})$, which reads*

$$\partial_T \begin{pmatrix} k \\ \alpha \\ \mathbf{M} \end{pmatrix} = \left( \begin{array}{cc|c} 0 & 1 & \\ 1 & 0 & \\ \hline & & \mathbf{B} \end{array} \right) \partial_X \left( \nabla_{k, \alpha, \mathbf{M}} \mathrm{H} \right). \tag{22}$$

*Proof.* We have already established the first assertion. However, for later use let us point out more precisely that the relations in (18) — linked to the fact that $\Theta$ is indeed an abbreviated action integral — imply that

$$k = \frac{1}{\partial_\mu \Theta}, \qquad \mathbf{M} = \frac{\nabla_\lambda \Theta}{\partial_\mu \Theta}, \qquad \alpha = \partial_c \Theta - (\partial_\mu \Theta)\, \mathcal{Q}\left(\frac{\nabla_\lambda \Theta}{\partial_\mu \Theta}\right), \qquad (23)$$

$$\partial_\mu \Theta = \frac{1}{k}, \qquad \nabla_\lambda \Theta = \frac{\mathbf{M}}{k}, \qquad \partial_c \Theta = \alpha + \frac{1}{k}\, \mathcal{Q}(\mathbf{M}). \qquad (24)$$

In order to compute the partial derivatives of H in the variables $(k, \alpha, \mathbf{M})$, it is expedient to use (20). Indeed, in this way, by combining classical cancellation of derivatives of conjugate functions with relations (24), we derive

$$\partial_k \mathrm{H} = \frac{\mathrm{H}}{k} + c\, \frac{\mathcal{Q}(\mathbf{M})}{k} + \boldsymbol{\lambda} \cdot \frac{\mathbf{M}}{k} + \frac{\mu}{k} = \Theta - \alpha c,$$

$$\partial_\alpha \mathrm{H} = -c\, k,$$

$$\nabla_\mathbf{M} \mathrm{H} = -c\, \nabla \mathcal{Q}(\mathbf{M}) - \boldsymbol{\lambda} = \langle \delta \mathcal{H}[\underline{\mathbf{U}}] \rangle,$$

where the last relation is obtained by averaging (15).

It follows at once that equations (12)-(13), which are also the first and last lines of (19), are equally written as

$$\partial_T k - \partial_x(\partial_\alpha \mathrm{H}) = 0, \qquad (25)$$

$$\partial_T \mathbf{M} - \partial_X(\mathbf{B}\, \nabla_\mathbf{M} \mathrm{H}) = 0. \qquad (26)$$

So the only remaining task is to manipulate (19) to obtain an equation for $\alpha$. By using (23), (19) and the symmetry of $\mathbf{B}$, one derives

$$(\partial_T + c\partial_X)\alpha = (\partial_T + c\partial_X)\left(\partial_c \Theta - \frac{\mathcal{Q}(\nabla_\lambda \Theta)}{\partial_\mu \Theta}\right)$$

$$= \partial_\mu \Theta\, \partial_X \mu + \nabla_\lambda \Theta\, \partial_X \boldsymbol{\lambda} + \frac{\mathcal{Q}(\nabla_\lambda \Theta)}{\partial_\mu \Theta}\, \partial_X c$$

$$= \partial_X \Theta - \alpha \partial_X c$$

thus

$$\partial_T \alpha = \partial_X(\partial_k \mathrm{H}). \qquad (27)$$

□

**Remark 7.** The '*symmetric*' form of (22) readily implies that H is a mathematical *entropy* for this system. Indeed, along smooth solutions of (22) we have

$$\partial_T \mathrm{H} = \partial_X \left((\partial_k \mathrm{H})(\partial_\alpha \mathrm{H}) + \frac{1}{2} (\nabla_\mathbf{M} \mathrm{H}) \cdot \mathbf{B} \nabla_\mathbf{M} \mathrm{H}\right)$$

by the symmetry of $\mathbf{B}$. For the sake of consistency, we now check that this conservation law for H coincides with the averaging of (2) — the original conservation law for $\mathcal{H}$ — along wave profiles. On one hand from (15), (16) and (20) stems

$$\langle \nabla_{\underline{\mathbf{U}}_x} \mathcal{H}[\underline{\mathbf{U}}] \cdot \partial_x(\mathbf{B}\delta\mathcal{H}[\underline{\mathbf{U}}])\rangle = -c\langle \nabla_{\underline{\mathbf{U}}_x} \mathcal{H}[\underline{\mathbf{U}}] \cdot \partial_x \underline{\mathbf{U}}\rangle$$

$$= -c\, (\mu + c\langle \mathcal{Q}(\underline{\mathbf{U}})\rangle + \boldsymbol{\lambda} \cdot \mathbf{M} + \mathrm{H})$$

$$= -c\, k\, \Theta.$$

On the other hand from (15) and the symmetry of **B** follows

$$\langle \tfrac{1}{2}\delta\mathcal{H}[\underline{\mathbf{U}}] \cdot \mathbf{B}\delta\mathcal{H}[\underline{\mathbf{U}}]\rangle = \tfrac{1}{2}\langle\delta\mathcal{H}[\underline{\mathbf{U}}]\rangle \cdot \mathbf{B}\langle\delta\mathcal{H}[\underline{\mathbf{U}}]\rangle + c^2\langle\mathcal{Q}(\underline{\mathbf{U}}-\mathbf{M})\rangle$$
$$= \tfrac{1}{2}\langle\delta\mathcal{H}[\underline{\mathbf{U}}]\rangle \cdot \mathbf{B}\langle\delta\mathcal{H}[\underline{\mathbf{U}}]\rangle + c^2\left(\langle\mathcal{Q}(\underline{\mathbf{U}})\rangle - \mathcal{Q}(\mathbf{M})\right)$$
$$= \tfrac{1}{2}\langle\delta\mathcal{H}[\underline{\mathbf{U}}]\rangle \cdot \mathbf{B}\langle\delta\mathcal{H}[\underline{\mathbf{U}}]\rangle + c^2\,k\,\alpha\,.$$

Combining the foregoing with (21) proves the claim.

**Remark 8.** The conservation law (14) itself also admits a nice formulation in terms of H. It equivalently reads

$$\partial_T \langle \mathcal{Q}(\underline{\mathbf{U}}) \rangle = \partial_X \mathrm{H}^*\,,$$

where

$$\mathrm{H}^* := k\partial_k \mathrm{H} + \alpha \partial_\alpha \mathrm{H} + \mathbf{M}\cdot \nabla_\mathbf{M} \mathrm{H} - \mathrm{H}\,,$$

($\mathrm{H}^*$ would be the *conjugate* function of H if this function were strictly convex). Indeed it already follows from previous computations that

$$\langle \mathsf{L}\mathcal{H}[\underline{\mathbf{U}}]\rangle \;=\; k\Theta - \mathrm{H}\,.$$

Moreover from (15) and the symmetry of **B** we deduce

$$\langle \underline{\mathbf{U}} \cdot \delta\mathcal{H}[\underline{\mathbf{U}}]\rangle = \mathbf{M}\cdot\langle\delta\mathcal{H}[\underline{\mathbf{U}}]\rangle - 2\,c\,\langle\mathcal{Q}(\underline{\mathbf{U}}-\mathbf{M})\rangle$$
$$= \mathbf{M}\cdot\langle\delta\mathcal{H}[\underline{\mathbf{U}}]\rangle - 2\,c\,k\,\alpha$$

so that the claim follows from (21).

The quasilinear form of (22) reads

$$\partial_T \begin{pmatrix} k \\ \alpha \\ \mathbf{M} \end{pmatrix} \;=\; \mathbb{B}\,\nabla^2_{k,\alpha,\mathbf{M}}\mathrm{H}\,\partial_X \begin{pmatrix} k \\ \alpha \\ \mathbf{M} \end{pmatrix}\,, \qquad (28)$$

where

$$\mathbb{B} := \left(\begin{array}{cc|c} 0 & 1 & \\ 1 & 0 & \\ \hline & & \mathbf{B} \end{array}\right),$$

so that (22) is hyperbolic at points in the state space where the matrix $\mathbb{B}\,\nabla^2_{k,\alpha,\mathbf{M}}\mathrm{H}$ is diagonalizable with real eigenvalues. A first, natural approach to check hyperbolicity is to try and use the symmetry of the matrices $\nabla^2_{k,\alpha,\mathbf{M}}\mathrm{H}$ and $\mathbb{B}$.

**Corollary 2.** *In the framework of Theorem 1, if* H *is a* strictly convex *function of* $(k,\alpha,\mathbf{M})$, *then the modulated system* (22) *is hyperbolic.*

*Proof.* This follows from the fact that the Hessian $\nabla^2_{k,\alpha,\mathbf{M}}\mathrm{H}$ of H is a symmetrizer for (22) whenever $\nabla^2_{k,\alpha,\mathbf{M}}\mathrm{H}$ is positive definite. Indeed, as soon as $\nabla^2_{k,\alpha,\mathbf{M}}\mathrm{H}$ is nonsingular the quasilinear form (28) of (22) is equivalent to

$$\nabla^2_{k,\alpha,\mathbf{M}}\mathrm{H}\,\partial_T \begin{pmatrix} k \\ \alpha \\ \mathbf{M} \end{pmatrix} = \nabla^2_{k,\alpha,\mathbf{M}}\mathrm{H}\,\mathbb{B}\,\nabla^2_{k,\alpha,\mathbf{M}}\mathrm{H}\,\partial_X \begin{pmatrix} k \\ \alpha \\ \mathbf{M} \end{pmatrix}\,.$$

Since the matrix $\nabla^2_{k,\alpha,\mathbf{M}}\mathrm{H}\,\mathbb{B}\,\nabla^2_{k,\alpha,\mathbf{M}}\mathrm{H}$ is symmetric, if in addition $\nabla^2_{k,\alpha,\mathbf{M}}\mathrm{H}$ is positive definite then (22) is necessarily hyperbolic by a standard observation in the theory of hyperbolic systems (see for instance [Ser99, Theorem 3.1.6]). □

However, our numerical experiments tend to show that $\nabla^2_{k,\alpha,\mathbf{M}}\mathrm{H}$ is hardly ever definite positive [Mie17]. Moreover, as made explicit in Remark 14, our analysis proves that $\nabla^2_{k,\alpha,\mathbf{M}}\mathrm{H}$ is not definite positive in either one of the small amplitude limit and the large wavelength limit. Indeed, the upper diagonal block in the limits of $\nabla^2_{k,\alpha,\mathbf{M}}\mathrm{H}$ found in Theorems 5 and 6 has signature $(1,1)$, and therefore $\nabla^2_{k,\alpha,\mathbf{M}}\mathrm{H}$ cannot be definite.

The main purpose of subsequent sections is to draw rigorous conclusions on the modulated system in quasilinear form (28), in the small amplitude and soliton limits, when either $\alpha \to 0$ or $k \to 0$. Required expansions are derived from expansions of $\nabla^2_{\mu,c,\boldsymbol{\lambda}}\Theta$ obtained in [BGMR20]. Thus, before going to the most technical part of the present paper, we need to point out the explicit connection between the Hessian of the averaged Hamiltonian H as as function of $(k,\alpha,\mathbf{M})$ and the Hessian of the abbreviated action $\Theta$ as a function of parameters $(\mu,c,\boldsymbol{\lambda})$.

**Proposition 3.** *In the framework of Theorem 1,*

$$\nabla^2_{k,\alpha,\mathbf{M}}\mathrm{H} = -\frac{1}{k}\mathbb{A}^\mathsf{T}(\nabla^2_{\mu,c,\boldsymbol{\lambda}}\Theta)^{-1}\mathbb{A} \;-\; c\,\mathbb{B}^{-1}\,, \tag{29}$$

*with*

$$\mathbb{A} = \mathbb{A}(k,\mathbf{M}) := \begin{pmatrix} -\frac{1}{k} & 0 & 0 \\ -\frac{\mathcal{Q}(\mathbf{M})}{k} & k & \nabla_\mathbf{U}\mathcal{Q}(\mathbf{M})^\mathsf{T} \\ -\frac{\mathbf{M}}{k} & 0 & \mathbf{I}_N \end{pmatrix}.$$

*Proof.* Along the proof we find it convenient to use first and second differentials, denoted with d and $\mathrm{d}^2$, rather than gradients and Hessians. We proceed by differentiating at points $(k,\alpha,\mathbf{M})$ (left implicit) in an arbitrary direction $(\dot{k},\dot{\alpha},\dot{\mathbf{M}})$ (made explicit). In the present proof all functions are implicitly considered as functions of $(k,\alpha,\mathbf{M})$.

The starting point is the differentiation of (20), already used in the proof of Theorem 1,

$$\mathrm{d}\mathrm{H}(\dot{k},\dot{\alpha},\dot{\mathbf{M}}) = \frac{\dot{k}}{k}\mathrm{H} - k\begin{pmatrix}\mu\\c\\\boldsymbol{\lambda}\end{pmatrix}\cdot \mathrm{d}(\nabla_{\mu,c,\boldsymbol{\lambda}}\Theta)(\dot{k},\dot{\alpha},\dot{\mathbf{M}})$$

that we differentiate once more to derive

$$\mathrm{d}^2\mathrm{H}((\dot{k},\dot{\alpha},\dot{\mathbf{M}}),(\dot{k},\dot{\alpha},\dot{\mathbf{M}})) \tag{30}$$

$$= -k\,\mathrm{d}\begin{pmatrix}\mu\\c\\\boldsymbol{\lambda}\end{pmatrix}(\dot{k},\dot{\alpha},\dot{\mathbf{M}})\cdot \nabla^2_{\mu,c,\boldsymbol{\lambda}}\Theta\,\mathrm{d}\begin{pmatrix}\mu\\c\\\boldsymbol{\lambda}\end{pmatrix}(\dot{k},\dot{\alpha},\dot{\mathbf{M}})$$

$$-\frac{1}{k}\begin{pmatrix}\mu\\c\\\boldsymbol{\lambda}\end{pmatrix}\cdot \mathrm{d}\Big[k^2\mathrm{d}(\nabla_{\mu,c,\boldsymbol{\lambda}}\Theta)(\dot{k},\dot{\alpha},\dot{\mathbf{M}})\Big](\dot{k},\dot{\alpha},\dot{\mathbf{M}})\,.$$

Now differentiating (24) yields

$$\mathrm{d}(\nabla_{\mu,c,\boldsymbol{\lambda}}\Theta)(\dot{k},\dot{\alpha},\dot{\mathbf{M}}) = \frac{1}{k}\mathbb{A}\begin{pmatrix}\dot{k}\\\dot{\alpha}\\\dot{\mathbf{M}}\end{pmatrix}$$

which also implies

$$\mathrm{d}\begin{pmatrix}\mu\\c\\\boldsymbol{\lambda}\end{pmatrix}(\dot{k},\dot{\alpha},\dot{\mathbf{M}}) = \frac{1}{k}(\nabla^2_{\mu,c,\boldsymbol{\lambda}}\Theta)^{-1}\mathbb{A}\begin{pmatrix}\dot{k}\\\dot{\alpha}\\\dot{\mathbf{M}}\end{pmatrix}. \tag{31}$$

In turn
$$\mathrm{d}(k\mathbb{A})(\dot{k},\dot{\alpha},\dot{\mathbf{M}}) = \begin{pmatrix} 0 & 0 & 0 \\ -\nabla_{\mathbf{U}}\mathcal{Q}(\mathbf{M})\cdot\dot{\mathbf{M}} & 2k\dot{k} & \dot{k}\nabla_{\mathbf{U}}\mathcal{Q}(\mathbf{M})^{\mathsf{T}} + k\nabla_{\mathbf{U}}\mathcal{Q}(\dot{\mathbf{M}})^{\mathsf{T}} \\ -\dot{\mathbf{M}} & 0 & \dot{k}\,\mathbf{I}_N \end{pmatrix}$$

so that

$$\mathrm{d}\Big[k^2 \mathrm{d}(\nabla_{\mu,c,\boldsymbol{\lambda}}\Theta)(\dot{k},\dot{\alpha},\dot{\mathbf{M}})\Big](\dot{k},\dot{\alpha},\dot{\mathbf{M}}) = \begin{pmatrix} 0 \\ 2k\,\dot{k}\,\dot{\alpha} + 2k\,\mathcal{Q}(\dot{\mathbf{M}}) \\ 0 \end{pmatrix}. \tag{32}$$

Inserting (31) and (32) in (30) achieves the proof by identification of relevant symmetric matrices with corresponding quadratic forms since

$$2\dot{k}\,\dot{\alpha} + 2\,\mathcal{Q}(\dot{\mathbf{M}}) = \begin{pmatrix} \dot{k} \\ \dot{\alpha} \\ \dot{\mathbf{M}} \end{pmatrix} \cdot \mathbb{B}^{-1} \begin{pmatrix} \dot{k} \\ \dot{\alpha} \\ \dot{\mathbf{M}} \end{pmatrix}.$$

□

**Remark 9.** Relation (29) leaves the possibility for $\nabla^2_{k,\alpha,\mathbf{M}}H$ to be definite without $\nabla^2_{\mu,c,\boldsymbol{\lambda}}\Theta$ being so, and vice and versa. This could be of importance since any of those yields hyperbolicity of the modulated system and it was shown in [BGMR16] that the negative signature of $\nabla^2_{\mu,c,\boldsymbol{\lambda}}\Theta$ must be equal to $N$ modulo an even number for the underlying wave to be spectrally stable. Yet as already mentioned, in practice, this is hardly ever the case; see [Mie17] for numerical experiments, and [BGMR20] and Remark 14 for the analysis of signatures in either one of the extreme regimes, small amplitude or large wavelength.

As seen on the quasilinear form (28), the characteristic matrix of System (22) reads

$$\mathbb{W} := -\mathbb{B}\,\nabla^2_{k,\alpha,\mathbf{M}}H\,. \tag{33}$$

We refer to $\mathbb{W}$ in the sequel as the Whitham matrix of (22). It can be rewritten using Equation (29) as

$$\mathbb{W} = \frac{1}{k}\,\mathbb{B}\,\mathbb{A}^{\mathsf{T}}\,(\nabla^2_{\mu,c,\boldsymbol{\lambda}}\Theta)^{-1}\,\mathbb{A}\, + \, c\,\mathbf{I}_{N+2}\,.$$

**Remark 10.** For comparison, the characteristic matrix of System (19) in variables $(\mu,c,\boldsymbol{\lambda})$ is

$$\frac{1}{k}\,(\nabla^2_{\mu,c,\boldsymbol{\lambda}}\Theta)^{-1}\,\mathbb{S}\, + \, c\,\mathbf{I}_{N+2}$$

with

$$\mathbb{S} := \left(\begin{array}{cc|c} 0 & -1 & \\ -1 & 0 & \\ \hline & & \mathbf{B} \end{array}\right). \tag{34}$$

As follows from (31), $(\nabla^2_{\mu,c,\boldsymbol{\lambda}}\Theta)^{-1}\,\mathbb{A}$ provides a change of basis between characteristic matrices. This may be checked directly thanks to the identity

$$\mathbb{S} = \mathbb{A}\,\mathbb{B}\,\mathbb{A}^{\mathsf{T}}.$$

# 3 Asymptotic expansions of parameters

## 3.1 Expansions of action derivatives

Our study of extreme regimes hinges on asymptotic expansions of the action and its derivatives, obtained in [BGMR20] and that we partially recall here.

To conveniently write some of the coefficients of the expansions, we first make more explicit the profile equations (15)-(16). As in [BGMR20] we introduce the potential $\mathcal{W}(v; c, \boldsymbol{\lambda})$ defined in the case $N = 1$ by

$$\mathcal{W}(v; c, \lambda) := - f(v) - \frac{1}{2}\frac{c}{b} v^2 - \lambda v,$$

and in the case $N = 2$ by

$$\mathcal{W}(v; c, \boldsymbol{\lambda}) := - f(v) - \frac{1}{2}\tau(v)\, g(v; c, \lambda_2)^2 - \frac{c}{b}\, v\, g(v; c, \lambda_2) - \lambda_1\, v - \lambda_2\, g(v; c, \lambda_2)$$

with

$$g(v; c, \lambda) := -\frac{1}{\tau(v)}\left(\frac{c}{b} v + \lambda\right).$$

The point is that (15)-(16) is equivalently written

$$\kappa(\underline{v})\,\underline{v}_{xx} + \tfrac{1}{2}\kappa'(\underline{v})\,\underline{v}_x^2 + \mathcal{W}'(\underline{v}; c, \boldsymbol{\lambda}) = 0, \qquad \tfrac{1}{2}\kappa(\underline{v})\,\underline{v}_x^2 + \mathcal{W}(\underline{v}; c, \boldsymbol{\lambda}) = \mu,$$

completed, in the case $N = 2$, with

$$\underline{u} = g(\underline{v}; c, \lambda_2).$$

We only consider non-degenerate limits. The nature of the non-degeneracy is made precise in the following set-up.

**Assumption 2.** <u>Harmonic limit</u> *Fix $(\underline{\mu}_0, \underline{c}_0, \boldsymbol{\underline{\lambda}}_0) \in \overline{\Omega}$ such that there exists $\underline{v}_0 > 0$ such that*

$$\underline{\mu}_0 = \mathcal{W}(\underline{v}_0; \underline{c}_0, \boldsymbol{\underline{\lambda}}_0), \qquad \partial_v \mathcal{W}(\underline{v}_0; \underline{c}_0, \boldsymbol{\underline{\lambda}}_0) = 0, \qquad \partial_v^2 \mathcal{W}(\underline{v}_0; \underline{c}_0, \boldsymbol{\underline{\lambda}}_0) > 0.$$

*Then there exists $\Lambda$ a connected open neighborhood of $(\underline{c}_0, \boldsymbol{\underline{\lambda}}_0)$ and smooth functions $v_0 : \Lambda \to (0, \infty)$ and $\mu_0 : \Lambda \to \mathbb{R}$ such that $(v_0, \mu_0)(\underline{c}_0, \boldsymbol{\underline{\lambda}}_0) = (\underline{v}_0, \underline{\mu}_0)$ and for any $(c, \boldsymbol{\lambda}) \in \Lambda$*

$$\mu_0(c, \boldsymbol{\lambda}) = \mathcal{W}(v_0(c, \boldsymbol{\lambda}); c, \boldsymbol{\lambda}), \qquad \partial_v \mathcal{W}(v_0(c, \boldsymbol{\lambda}); c, \boldsymbol{\lambda}) = 0, \qquad \partial_v^2 \mathcal{W}(v_0(c, \boldsymbol{\lambda}); c, \boldsymbol{\lambda}) > 0.$$

*Moreover one may ensure[9] that for some $r_0 > 0$*

$$\Omega_0^{r_0} := \bigcup_{(c, \boldsymbol{\lambda}) \in \Lambda} (\mu_0(c, \boldsymbol{\lambda}), \mu_0(c, \boldsymbol{\lambda}) + r_0) \times \{(c, \boldsymbol{\lambda})\} \subset \Omega$$

*and there exist $v_2$ and $v_3$ smooth maps defined on $\Omega_0^{r_0}$ such that for any $\boldsymbol{\mu} = (\mu, c, \boldsymbol{\lambda}) \in \Omega_0^{r_0}$,*

$$0 < v_2(\boldsymbol{\mu}) < v_0(c, \boldsymbol{\lambda}) < v_3(\boldsymbol{\mu}), \qquad \mu = \mathcal{W}(v_2(\boldsymbol{\mu}); c, \boldsymbol{\lambda}) = \mathcal{W}(v_3(\boldsymbol{\mu}); c, \boldsymbol{\lambda}),$$
$$\partial_v \mathcal{W}(v_2(\boldsymbol{\mu}); c, \boldsymbol{\lambda}) \neq 0, \qquad \partial_v \mathcal{W}(v_3(\boldsymbol{\mu}); c, \boldsymbol{\lambda}) \neq 0,$$
$$\forall v \in (v_2(\boldsymbol{\mu}), v_3(\boldsymbol{\mu})), \quad \mu \neq \mathcal{W}(v; c, \boldsymbol{\lambda}).$$

<u>Soliton limit</u> *Fix $(\underline{\mu}_s, \underline{c}_s, \boldsymbol{\underline{\lambda}}_s) \in \overline{\Omega}$ such that there exists $\underline{v}_s$ and $\underline{v}^s$ such that[10]*

$$0 < \underline{v}_s < \underline{v}^s, \qquad \underline{\mu}_s = \mathcal{W}(\underline{v}_s; \underline{c}_s, \boldsymbol{\underline{\lambda}}_s) = \mathcal{W}(\underline{v}^s; \underline{c}_s, \boldsymbol{\underline{\lambda}}_s),$$
$$\partial_v \mathcal{W}(\underline{v}_s; \underline{c}_s, \boldsymbol{\underline{\lambda}}_s) = 0, \qquad \partial_v^2 \mathcal{W}(\underline{v}_s; \underline{c}_s, \boldsymbol{\underline{\lambda}}_s) < 0,$$
$$\partial_v \mathcal{W}(\underline{v}^s; \underline{c}_s, \boldsymbol{\underline{\lambda}}_s) \neq 0, \qquad \text{and} \qquad \forall v \in (\underline{v}_s, \underline{v}^s), \quad \underline{\mu}_s \neq \mathcal{W}(v; \underline{c}_s, \boldsymbol{\underline{\lambda}}_s).$$

---

[9] Up to choosing the correct branch of parametrization and extending $\Omega$ if necessary. Implicitly $\underline{v}$ is chosen consistently.

[10] The choice that $\partial_v \mathcal{W}(\underline{v}_s; \underline{c}_s, \boldsymbol{\underline{\lambda}}_s) = 0$ and $\partial_v \mathcal{W}(\underline{v}^s; \underline{c}_s, \boldsymbol{\underline{\lambda}}_s) \neq 0$ instead of $\partial_v \mathcal{W}(\underline{v}_s; \underline{c}_s, \boldsymbol{\underline{\lambda}}_s) \neq 0$ and $\partial_v \mathcal{W}(\underline{v}^s; \underline{c}_s, \boldsymbol{\underline{\lambda}}_s) = 0$ is arbitrary and purely made for the sake of clarity and definiteness. There is no loss of generality since one may go from one case to the other by rewriting the equations for $\underline{v}$ in terms of $-\underline{v}$.

*Then there exists $\Lambda$ a connected open neighborhood of $(\underline{c}_s, \underline{\boldsymbol{\lambda}}_s)$ and smooth functions $v_s : \Lambda \to (0, \infty)$, $v^s : \Lambda \to (0, \infty)$ and $\mu_s : \Lambda \to \mathbb{R}$ such that $(v_s, v^s, \mu_0)(\underline{c}_s, \underline{\boldsymbol{\lambda}}_s) = (\underline{v}_s, \underline{v}^s, \underline{\mu}_s)$ and for any $(c, \boldsymbol{\lambda}) \in \Lambda$*

$$0 < v_s(c, \boldsymbol{\lambda}) < v^s(c, \boldsymbol{\lambda}), \qquad \mu_s(c, \boldsymbol{\lambda}) = \mathcal{W}(v_s(c, \boldsymbol{\lambda}); c, \boldsymbol{\lambda}) = \mathcal{W}(v^s(c, \boldsymbol{\lambda}); c, \boldsymbol{\lambda}),$$
$$\partial_v \mathcal{W}(v_s(c, \boldsymbol{\lambda}); c, \boldsymbol{\lambda}) = 0, \qquad \partial_v^2 \mathcal{W}(v_s(c, \boldsymbol{\lambda}); c, \boldsymbol{\lambda}) < 0,$$
$$\partial_v \mathcal{W}(v^s(c, \boldsymbol{\lambda}); c, \boldsymbol{\lambda}) \neq 0, \qquad \text{and} \qquad \forall v \in (v_s(c, \boldsymbol{\lambda}), v^s(c, \boldsymbol{\lambda})), \quad \mu_s(c, \boldsymbol{\lambda}) \neq \mathcal{W}(v; c, \boldsymbol{\lambda}).$$

*Moreover one may ensure that for some $r_0 > 0$*

$$\Omega_s^{r_0} := \bigcup_{(c, \boldsymbol{\lambda}) \in \Lambda} (\mu_s(c, \boldsymbol{\lambda}) - r_0, \mu_s(c, \boldsymbol{\lambda})) \times \{(c, \boldsymbol{\lambda})\} \subset \Omega$$

*and there exist $v_1$, $v_2$ and $v_3$ three smooth maps defined on $\Omega_s^{r_0}$ such that for any $\boldsymbol{\mu} = (\mu, c, \boldsymbol{\lambda}) \in \Omega_s^{r_0}$,*

$$0 < v_1(\boldsymbol{\mu}) < v_s(c, \boldsymbol{\lambda}) < v_2(\boldsymbol{\mu}) < v_3(\boldsymbol{\mu}) < v^s(c, \boldsymbol{\lambda}),$$
$$\mu = \mathcal{W}(v_1(\boldsymbol{\mu}); c, \boldsymbol{\lambda}) = \mathcal{W}(v_2(\boldsymbol{\mu}); c, \boldsymbol{\lambda}) = \mathcal{W}(v_3(\boldsymbol{\mu}); c, \boldsymbol{\lambda}),$$
$$\partial_v \mathcal{W}(v_1(\boldsymbol{\mu}); c, \boldsymbol{\lambda}) \neq 0, \qquad \partial_v \mathcal{W}(v_2(\boldsymbol{\mu}); c, \boldsymbol{\lambda}) \neq 0, \qquad \partial_v \mathcal{W}(v_3(\boldsymbol{\mu}); c, \boldsymbol{\lambda}) \neq 0,$$
$$\forall v \in (v_1(\boldsymbol{\mu}), v_2(\boldsymbol{\mu})) \cup (v_2(\boldsymbol{\mu}), v_3(\boldsymbol{\mu})), \quad \mu \neq \mathcal{W}(v; c, \boldsymbol{\lambda}).$$

For all $(c_*, \boldsymbol{\lambda}_*) \in \Lambda$, we consider

$$\text{either} \qquad \boldsymbol{\mu}_{0*} := (c_*, \boldsymbol{\lambda}_*, \mu_0(c_*, \boldsymbol{\lambda}_*)), \qquad \text{or} \qquad \boldsymbol{\mu}_{s*} := (c_*, \boldsymbol{\lambda}_*, \mu_s(c_*, \boldsymbol{\lambda}_*)),$$

which both belong to $\overline{\Omega}$, and the corresponding harmonic limit ($\boldsymbol{\mu} \xrightarrow{\Omega_0^{r_0}} \boldsymbol{\mu}_{0*}$) and soliton limit ($\boldsymbol{\mu} \xrightarrow{\Omega_s^{r_0}} \boldsymbol{\mu}_{s*}$). Actually it is more convenient and sufficient to *fix* $(c_*, \boldsymbol{\lambda}_*) \in \Lambda$ and consider either $\mu \to \mu_0(c_*, \boldsymbol{\lambda}_*)^+$ or $\mu \to \mu_s(c_*, \boldsymbol{\lambda}_*)^-$, provided one ensures local uniformity with respect to $(c_*, \boldsymbol{\lambda}_*) \in \Lambda$. By acting in this way, in [BGMR20] we derived asymptotic expansions in terms of two small parameters going to zero:

$$\delta := (v_3 - v_2)/2$$

in the harmonic limit and

$$\varrho := \frac{v_2 - v_1}{v_3 - v_2}$$

in the soliton limit. These expansions are recalled below after a few preliminaries.

First, for the sake of concision, in the case $N = 2$ we introduce notation $q(v; c, \lambda) := \mathcal{Q}(v, g(v; c, \lambda))$, still with $g(v; c, \lambda) = -((c/b) v + \lambda)/\tau(v)$. Note that in the sequel $g$ and $q$ are evaluated at $\lambda = \lambda_2$, the second component of $\boldsymbol{\lambda}$. For convenience we adopt a similar convention in the case $N = 1$ with merely $q(v) := \mathcal{Q}(v)$. In the statement that follows, we omit to write the dependence — if any — of these functions on the parameters $(c, \boldsymbol{\lambda})$ in order to shorten formulas a little bit and stress symmetry between cases $N = 1$ and $N = 2$. We also make use of the symmetric matrix $\mathbb{S}$ defined in (34).

Now we introduce a set of vectors that are crucially involved in the above-mentioned asymptotic expansions, and provide associated key cancellations proved in [BGMR20, Lemma 1].

**Proposition 4** ([BGMR20]). *For both indices $i = 0$ and $i = s$ we introduce the following vectors of $\mathbb{R}^{N+2}$: for $N = 2$*

$$\mathbf{V}_i := \begin{pmatrix} 1 \\ q(v_i) \\ v_i \\ g(v_i) \end{pmatrix}, \qquad \mathbf{W}_i := \begin{pmatrix} 0 \\ \partial_v q(v_i) \\ 1 \\ \partial_v g(v_i) \end{pmatrix}, \qquad \mathbf{Z}_i := \begin{pmatrix} 0 \\ \partial_v^2 q(v_i) \\ 0 \\ \partial_v^2 g(v_i) \end{pmatrix}, \qquad (35)$$

$$\mathbf{T}_i := \frac{1}{\sqrt{\tau(v_i)}} \begin{pmatrix} 0 \\ \frac{v_i}{b} \\ 0 \\ 1 \end{pmatrix}, \qquad \mathbf{E} := \begin{pmatrix} 1 \\ 0 \\ 0 \\ 0 \end{pmatrix} = \mathbb{S}^{-1} \mathbf{F}, \qquad \mathbf{F} := \begin{pmatrix} 0 \\ -1 \\ 0 \\ 0 \end{pmatrix};$$

*and for $N = 1$*

$$\mathbf{V}_i := \begin{pmatrix} 1 \\ q(v_i) \\ v_i \end{pmatrix}, \qquad \mathbf{W}_i := \begin{pmatrix} 0 \\ \partial_v q(v_i) \\ 1 \end{pmatrix}, \qquad \mathbf{Z}_i := \begin{pmatrix} 0 \\ \partial_v^2 q(v_i) \\ 0 \end{pmatrix}, \qquad (36)$$

$$\mathbf{T}_i := \begin{pmatrix} 0 \\ 0 \\ 0 \end{pmatrix}, \qquad \mathbf{E} := \begin{pmatrix} 1 \\ 0 \\ 0 \end{pmatrix} = \mathbb{S}^{-1}\mathbf{F}, \qquad \mathbf{F} := \begin{pmatrix} 0 \\ -1 \\ 0 \end{pmatrix}.$$

*These vectors are such that*

$$\begin{cases} \mathbf{V}_i \cdot \mathbb{S}^{-1}\mathbf{V}_i = 0, \quad \mathbf{V}_i \cdot \mathbb{S}^{-1}\mathbf{W}_i = 0, \quad \mathbf{V}_i \cdot \mathbb{S}^{-1}\mathbf{T}_i = 0, \\ \mathbf{V}_i \cdot \mathbb{S}^{-1}\mathbf{Z}_i = -\mathbf{W}_i \cdot \mathbb{S}^{-1}\mathbf{W}_i, \quad \mathbf{T}_i \cdot \mathbb{S}^{-1}\mathbf{T}_i = 0, \quad \mathbf{T}_i \cdot \mathbb{S}^{-1}\mathbf{Z}_i = 0, \\ \mathbf{E} \cdot \mathbf{V}_i = 1, \quad \mathbf{E} \cdot \mathbf{W}_i = 0, \quad \mathbf{E} \cdot \mathbf{Z}_i = 0, \quad \mathbf{E} \cdot \mathbf{T}_i = 0. \end{cases} \qquad (37)$$

At last, we introduce the Boussinesq moment of instability involved in solitary wave limits. We stress that it is both convenient and classical to parameterize solitary wave profiles $\underline{\mathbf{U}}^s$ not by $(c, \boldsymbol{\lambda})$ but by $(c, \mathbf{U}_s)$ with $\mathbf{U}_s$ the corresponding endstate. The associated $\boldsymbol{\lambda}$ is then recovered through

$$\boldsymbol{\lambda} = \boldsymbol{\lambda}_s(c, \mathbf{U}_s) := -\nabla_{\mathbf{U}}(\mathcal{H} + c\mathcal{Q})(\mathbf{U}_s, 0)$$

and $\mu_s$ is simply obtained as

$$\mu_s = -(\mathcal{H} + c\mathcal{Q})(\mathbf{U}_s, 0) + \nabla_{\mathbf{U}}(\mathcal{H} + c\mathcal{Q})(\mathbf{U}_s, 0) \cdot \mathbf{U}_s.$$

The Boussinesq moment of instability is then defined as

$$\begin{aligned} \mathcal{M}(c, \mathbf{U}_s) &= \int_{-\infty}^{+\infty} \left( \mathcal{H}[\underline{\mathbf{U}}^s] + c\mathcal{Q}(\underline{\mathbf{U}}^s) + \boldsymbol{\lambda}_s \cdot \underline{\mathbf{U}}^s + \mu_s \right) \mathrm{d}\xi \\ &= \int_{-\infty}^{+\infty} \left( (\mathcal{H} + c\mathcal{Q})[\underline{\mathbf{U}}^s] - (\mathcal{H} + c\mathcal{Q})(\mathbf{U}_s, 0) - \nabla_{\mathbf{U}}(\mathcal{H} + c\mathcal{Q})(\mathbf{U}_s, 0) \cdot (\underline{\mathbf{U}}^s - \mathbf{U}_s) \right) \mathrm{d}\xi. \end{aligned}$$

Note that, since $\delta(\mathcal{H} + c\mathcal{Q})[\underline{\mathbf{U}}^s] + \boldsymbol{\lambda}_s = 0$, we do have

$$\partial_c \mathcal{M}(c, \mathbf{U}_s) = \int_{-\infty}^{+\infty} \left( \mathcal{Q}(\underline{\mathbf{U}}^s) - \mathcal{Q}(\mathbf{U}_s) - \nabla_{\mathbf{U}}\mathcal{Q}(\mathbf{U}_s) \cdot (\underline{\mathbf{U}}^s - \mathbf{U}_s) \right) \mathrm{d}\xi$$

The following statement gathers elements from [BGMR20, Theorems 4 and 5] and their proofs.

**Theorem 2** ([BGMR20])**.** *Under Assumptions 1-2 we have the following asymptotics for the action derivatives.*
<u>Harmonic limit</u> *There exist real numbers $\mathfrak{a}_0$, $\mathfrak{b}_0$ and a positive number $\mathfrak{c}_0$ — depending smoothly on the parameters $(c, \boldsymbol{\lambda})$ — such that when $\delta$ goes to zero*

$$\frac{4\mathfrak{c}_0}{\Xi_0} \nabla_{\mu,c,\boldsymbol{\lambda}}\Theta = 4\mathfrak{c}_0 \mathbf{V}_0 + (\mathfrak{a}_0 \mathbf{V}_0 + \mathfrak{b}_0 \mathbf{W}_0 + \mathfrak{c}_0 \mathbf{Z}_0)\delta^2 + \mathcal{O}(\delta^4) \qquad (38)$$

$$\frac{1}{\Xi_0} \nabla^2_{\mu,c,\boldsymbol{\lambda}}\Theta = \mathfrak{a}_0 \mathbf{V}_0 \otimes \mathbf{V}_0 + \mathfrak{b}_0 (\mathbf{V}_0 \otimes \mathbf{W}_0 + \mathbf{W}_0 \otimes \mathbf{V}_0) - \mathbf{T}_0 \otimes \mathbf{T}_0 \qquad (39)$$
$$+ 2\mathfrak{c}_0 \mathbf{W}_0 \otimes \mathbf{W}_0 + \mathfrak{c}_0 (\mathbf{V}_0 \otimes \mathbf{Z}_0 + \mathbf{Z}_0 \otimes \mathbf{V}_0) + \mathcal{O}(\delta^2)$$

where $\Xi_0$ denotes the harmonic period at $v_0$, that is, $\Xi_0 = \sqrt{\kappa(v_0)/\partial_v^2\mathcal{W}(v_0;c,\boldsymbol{\lambda}))}$.

<u>Soliton limit</u> *There exist real numbers $\mathfrak{a}_s$, $\mathfrak{b}_s$, positive numbers $\mathfrak{c}_s$, $\mathfrak{h}_s$, a vector $\mathbf{X}_s$ and a symmetric matrix $\mathbb{O}_s$ — depending smoothly on the parameters $(c,\boldsymbol{\lambda})$ — such that*

$$\frac{\pi}{\Xi_s}\nabla_{\mu,c,\boldsymbol{\lambda}}\Theta = -\mathbf{V}_s \ln \varrho - \mathbf{X}_s + \frac{\varrho}{2}\mathbf{V}_s - \frac{1}{2\mathfrak{h}_s}\left(\mathfrak{a}_s\mathbf{V}_s + \mathfrak{b}_s\mathbf{W}_s + \mathfrak{c}_s\mathbf{Z}_s\right)\varrho^2\ln\varrho + \mathcal{O}(\varrho^2) \tag{40}$$

$$\frac{\pi}{\Xi_s}\nabla^2_{\mu,c,\boldsymbol{\lambda}}\Theta = \mathfrak{h}_s\frac{1+\varrho}{\varrho^2}\mathbf{V}_s \otimes \mathbf{V}_s + \left(\mathfrak{a}_s\mathbf{V}_s \otimes \mathbf{V}_s + \mathfrak{b}_s\left(\mathbf{V}_s \otimes \mathbf{W}_s + \mathbf{W}_s \otimes \mathbf{V}_s\right)\right)\ln\varrho \tag{41}$$

$$+ \left(\mathbf{T}_s \otimes \mathbf{T}_s + 2\mathfrak{c}_s\mathbf{W}_s \otimes \mathbf{W}_s + \mathfrak{c}_s\left(\mathbf{Z}_s \otimes \mathbf{V}_s + \mathbf{V}_s \otimes \mathbf{Z}_s\right)\right)\ln\varrho$$

$$+ \mathbb{O}_s + \mathcal{O}(\varrho \ln \varrho)$$

*when $\varrho$ goes to zero, where $\Xi_s$ denotes the harmonic period at $v_s$ of waves associated with the opposite 'capillarity' coefficient, that is, $\Xi_s := \sqrt{-\kappa(v_s)/\partial_v^2\mathcal{W}(v_s;c,\boldsymbol{\lambda}))}$. In addition, we have*[11]

$$\frac{\Xi_s}{\pi}(\mathbb{S}^{-1}\mathbf{V}_s) \cdot \mathbf{X}_s = \partial_c\mathcal{M}(c;\mathbf{U}_s), \qquad \frac{\Xi_s}{\pi}(\mathbb{S}^{-1}\mathbf{V}_s) \cdot \mathbb{O}_s\mathbb{S}^{-1}\mathbf{V}_s = \partial_c^2\mathcal{M}(c;\mathbf{U}_s), \tag{42}$$

*where $\mathbf{U}_s = v_s$ in the case $N=1$ and $\mathbf{U}_s = (v_s, g(v_s))$ in the case $N=2$.*

In the latter theorem and elsewhere in the present paper, for any two vectors $\mathbf{V}$ and $\mathbf{W}$ in $\mathbb{R}^d$, thought of as column vectors, $\mathbf{V} \otimes \mathbf{W}$ stands for the rank-one, square matrix of size $d$

$$\mathbf{V} \otimes \mathbf{W} = \mathbf{V}\mathbf{W}^\mathsf{T}$$

whatever $d$.

Observing that the matrices involved in the expansions of $\nabla^2_{\mu,c,\boldsymbol{\lambda}}\Theta$ in both the harmonic and the soliton limit have similar structures, we find useful to have at hand the following set of algebraic properties, which are either simple reformulations of relations in Proposition 4 or explicit computations from the definition of $\mathbb{A}$.

**Corollary 3.** • **Case $N=1$** *With*

$$\mathbb{P}_i := \mathbb{S}^{-1}\begin{pmatrix} \mathbf{F}_i & \mathbf{V}_i & \mathbf{W}_i \end{pmatrix}, \tag{43}$$

*we have*

$$\mathbb{D}_i := \mathbb{P}_i^\mathsf{T}\mathbb{S}\mathbb{P}_i = \mathbb{B}^{-1} = \left(\begin{array}{cc|c} 0 & 1 & 0 \\ 1 & 0 & 0 \\ \hline 0 & 0 & b^{-1} \end{array}\right),$$

$$\mathbb{P}_i^\mathsf{T}\mathbb{A} = \left(\begin{array}{c|c|c} -1/k & 0 & 0 \\ \hline \mathcal{Q}(v_i - \langle \underline{v} \rangle)/k & -k & (v_i - \langle \underline{v} \rangle)/b \\ \hline (v_i - \langle \underline{v} \rangle)/(bk) & 0 & 1/b \end{array}\right),$$

$$(\mathbb{P}_i^\mathsf{T}\mathbb{A})^{-1} = \left(\begin{array}{c|c|c} -k & 0 & 0 \\ \hline \mathcal{Q}(v_i - \langle \underline{v} \rangle)/k & -1/k & (v_i - \langle \underline{v} \rangle)/k \\ \hline v_i - \langle \underline{v} \rangle & 0 & b \end{array}\right),$$

---

[11]This comes from the proof of Theorem 5 in [BGMR20], the statement of which lacked the prefactor $\frac{\Xi_s}{\pi}$ in the relation between $\partial_c^2\mathcal{M}$ and the $\varrho^0$-term in the expansion of $\nabla^2_{\mu,c,\boldsymbol{\lambda}}\Theta$. We have corrected this omission in (42).

and, for any real numbers $(\mathfrak{a}, \mathfrak{b}, \mathfrak{c}, \mathfrak{m})$

$$\mathbb{P}_i{}^\mathsf{T}(\mathfrak{a}\mathbf{V}_i \otimes \mathbf{V}_i + \mathfrak{b}(\mathbf{V}_i \otimes \mathbf{W}_i + \mathbf{W}_i \otimes \mathbf{V}_i) + \mathfrak{m}\mathbf{W}_i \otimes \mathbf{W}_i + \mathfrak{c}(\mathbf{V}_i \otimes \mathbf{Z}_i + \mathbf{Z}_i \otimes \mathbf{V}_i))\mathbb{P}_i$$
$$= \left( \begin{array}{cc|c} \mathfrak{a} & -\mathfrak{c}\, b^{-1} & \mathfrak{b}\, b^{-1} \\ -\mathfrak{c}\, b^{-1} & 0 & 0 \\ \hline \mathfrak{b}\, b^{-1} & 0 & \mathfrak{m}\, b^{-2} \end{array} \right).$$

- **Case** $N = 2$ With

$$\mathbb{P}_i := \mathbb{S}^{-1} \begin{pmatrix} \mathbf{F}_i & \mathbf{V}_i & \mathbf{T}_i & \mathbf{W}_i \end{pmatrix}, \tag{44}$$

$$\begin{cases} \sigma_i &:= \mathbf{T}_i \cdot \mathbb{S}^{-1}\mathbf{W}_i = \dfrac{1}{b\sqrt{\tau(v)}}, \\ w_i &:= \mathbf{W}_i \cdot \mathbb{S}^{-1}\mathbf{W}_i = \dfrac{2g_v(v_i)}{b}, \\ \zeta_i &:= \mathbf{Z}_i \cdot \mathbb{S}^{-1}\mathbf{W}_i = \dfrac{g_{vv}(v_i)}{b}, \end{cases} \tag{45}$$

and

$$\mathbf{A}_i := \begin{pmatrix} 0 & 1/\sqrt{\tau(v_i)} \\ 1 & g_v(v_i) \end{pmatrix} = \begin{pmatrix} 0 & b\,\sigma_i \\ 1 & \frac{b}{2}\, w_i \end{pmatrix}, \tag{46}$$

we have

$$\mathbb{D}_i := \mathbb{P}_i{}^\mathsf{T}\mathbb{S}\mathbb{P}_i = \begin{pmatrix} 0 & 1 & 0 & 0 \\ 1 & 0 & 0 & 0 \\ \hline 0 & 0 & 0 & \sigma_i \\ 0 & 0 & \sigma_i & w_i \end{pmatrix},$$

$$\mathbb{P}_i{}^\mathsf{T}\mathbb{A} = \left( \begin{array}{c|c|c} -1/k & 0 & 0 \\ \hline \mathcal{Q}(\mathbf{U}_i - \mathbf{M})/k & -k & (\mathbf{U}_i - \mathbf{M})^\mathsf{T}\mathbf{B}^{-1} \\ \hline \mathbf{A}_i\,\mathbf{B}^{-1}(\mathbf{U}_i - \mathbf{M})/k & 0 & \mathbf{A}_i\,\mathbf{B}^{-1} \end{array} \right),$$

$$(\mathbb{P}_i{}^\mathsf{T}\mathbb{A})^{-1} = \left( \begin{array}{c|c|c} -k & 0 & 0 \\ \hline \mathcal{Q}(\mathbf{U}_i - \mathbf{M})/k & -1/k & (\mathbf{U}_i - \mathbf{M})^\mathsf{T}\mathbf{A}_i^{-1}/k \\ \hline \mathbf{U}_i - \mathbf{M} & 0 & \mathbf{B}\,\mathbf{A}_i^{-1} \end{array} \right),$$

and, for any real numbers $(\mathfrak{a}, \mathfrak{b}, \mathfrak{c}, \mathfrak{m}, \mathfrak{n})$

$$\mathbb{P}_i{}^\mathsf{T}(\mathfrak{a}\mathbf{V}_i \otimes \mathbf{V}_i + \mathfrak{b}(\mathbf{V}_i \otimes \mathbf{W}_i + \mathbf{W}_i \otimes \mathbf{V}_i) + \mathfrak{m}\mathbf{W}_i \otimes \mathbf{W}_i + \mathfrak{c}(\mathbf{V}_i \otimes \mathbf{Z}_i + \mathbf{Z}_i \otimes \mathbf{V}_i) + \mathfrak{n}\mathbf{T}_i \otimes \mathbf{T}_i)\mathbb{P}_i$$
$$= \left( \begin{array}{cc|cc} \mathfrak{a} & -\mathfrak{c}\, w_i & \mathfrak{b}\, \sigma_i & \mathfrak{b}\, w_i + \mathfrak{c}\, \zeta_i \\ -\mathfrak{c}\, w_i & 0 & 0 & 0 \\ \hline \mathfrak{b}\, \sigma_i & 0 & \mathfrak{m}\, \sigma_i^2 & \mathfrak{m}\, \sigma_i\, w_i \\ \mathfrak{b}\, w_i + \mathfrak{c}\, \zeta_i & 0 & \mathfrak{m}\, \sigma_i\, w_i & \mathfrak{m}\, w_i^2 + \mathfrak{n}\, \sigma_i^2 \end{array} \right).$$

For later reference, let us point out here that in any case

$$\mathbb{D}_i := \mathbb{P}_i{}^\mathsf{T}\mathbb{S}\mathbb{P}_i. \tag{47}$$

To unify cases $N = 1$ and $N = 2$, it is also useful to extend to $N = 1$ definitions in (45) and to set $\mathbf{A}_i = \begin{pmatrix} 1 \end{pmatrix}$ when $N = 1$.

## 3.2 Expansions of modulated variables

By (18) we have $\Xi = \partial_\mu \Theta$ so that we readily obtain expansions for the period $\Xi$ by projecting (38) and (40) onto their first component — which amounts to taking the inner product with $\mathbf{E}$. This gives

$$\Xi = \Xi_0 \left(1 + \frac{\mathfrak{a}_0}{4\mathfrak{c}_0} \delta^2 + \mathcal{O}(\delta^4)\right), \qquad \delta \to 0,$$

$$\Xi = \frac{\Xi_s}{\pi} \left(-\ln \varrho - \mathbf{E} \cdot \mathbf{X}_s + \frac{\varrho}{2} - \frac{\mathfrak{a}_s}{2\mathfrak{h}_s} \varrho^2 \ln \varrho + \mathcal{O}(\varrho^2)\right), \qquad \varrho \to 0,$$

from which we can of course infer expansions for the local wavenumber $k = 1/\Xi$

$$k = k_0 \left(1 - \frac{\mathfrak{a}_0}{4\mathfrak{c}_0} \delta^2 + \mathcal{O}(\delta^4)\right), \qquad \delta \to 0,$$

$$k = \frac{\pi}{\Xi_s} \left(-\frac{1}{\ln \varrho} + \frac{\mathbf{E} \cdot \mathbf{X}_s}{(\ln \varrho)^2} - \frac{\varrho}{2(\ln \varrho)^2} - \frac{(\mathbf{E} \cdot \mathbf{X}_s)^2}{(\ln \varrho)^3} + \frac{(\mathbf{E} \cdot \mathbf{X}_s) \varrho}{(\ln \varrho)^3} + \mathcal{O}\left(\frac{\varrho^2}{\ln \varrho}\right)\right), \qquad \varrho \to 0,$$

where $k_0 = 1/\Xi_0$.

Thanks to (18) again, the projections of (38) and (40) onto their intermediate and last components together with the expansions of $k$ yield expansions for the mean values $\langle \mathcal{Q}(\underline{\mathbf{U}}) \rangle$ and $\langle \underline{\mathbf{U}} \rangle$. To carry this out, it is convenient to introduce the $N \times (N+2)$ matrix

$$\mathbf{I} := \left(\begin{array}{cc|c} 0 & 0 & \mathbf{I}_N \end{array}\right)$$

of the projection onto last components, and to observe that taking the projection on the second component of vectors in $\mathbb{R}^{N+2}$ amounts to taking the inner product with $-\mathbf{F}$. We also recall that $\mathbf{U}_0 := \mathbf{I} \mathbf{V}_0$ and $\mathbf{U}_s := \mathbf{I} \mathbf{V}_s$.

Regarding the expansions of the mean value $\mathbf{M} = \langle \underline{\mathbf{U}} \rangle$ we get from (38) that

$$\begin{aligned} \mathbf{M} &= \frac{\Xi_0}{\Xi} \left(\mathbf{U}_0 + \tfrac{1}{4\mathfrak{c}_0} (\mathfrak{a}_0 \mathbf{U}_0 + \mathfrak{b}_0 \mathbf{I}\mathbf{W}_0 + \mathfrak{c}_0 \mathbf{I}\mathbf{Z}_0) \delta^2 + \mathcal{O}(\delta^4)\right) \\ &= \left(1 - \tfrac{\mathfrak{a}_0}{4\mathfrak{c}_0} \delta^2 + \mathcal{O}(\delta^4)\right) \left(\mathbf{U}_0 + \tfrac{1}{4\mathfrak{c}_0} (\mathfrak{a}_0 \mathbf{U}_0 + \mathfrak{b}_0 \mathbf{I}\mathbf{W}_0 + \mathfrak{c}_0 \mathbf{I}\mathbf{Z}_0) \delta^2 + \mathcal{O}(\delta^4)\right) \\ &= \mathbf{U}_0 + \mathbf{Y}_0 \delta^2 + \mathcal{O}(\delta^4), \end{aligned}$$

when $\delta$ goes to zero, with

$$\mathbf{Y}_0 := \frac{1}{4\mathfrak{c}_0} (\mathfrak{b}_0 \mathbf{I}\mathbf{W}_0 + \mathfrak{c}_0 \mathbf{I}\mathbf{Z}_0), \tag{48}$$

and from (40) that

$$\begin{aligned} \mathbf{M} &= \frac{\Xi_s}{\pi \Xi} \left(-\mathbf{U}_s \ln \varrho - \mathbf{I}\mathbf{X}_s + \frac{\varrho}{2} \mathbf{U}_s - \frac{1}{2\mathfrak{h}_s} (\mathfrak{a}_s \mathbf{U}_s + \mathfrak{b}_s \mathbf{I}\mathbf{W}_s + \mathfrak{c}_s \mathbf{I}\mathbf{Z}_s) \varrho^2 \ln \varrho + \mathcal{O}(\varrho^2)\right) \\ &= \left(-\frac{1}{\ln \varrho} + \frac{\mathbf{E} \cdot \mathbf{X}_s}{(\ln \varrho)^2} - \frac{\varrho}{2(\ln \varrho)^2} - \frac{(\mathbf{E} \cdot \mathbf{X}_s)^2}{(\ln \varrho)^3} + \frac{(\mathbf{E} \cdot \mathbf{X}_s) \varrho}{(\ln \varrho)^3} + \mathcal{O}\left(\frac{\varrho^2}{\ln \varrho}\right)\right) \times \\ &\qquad \left(-\mathbf{U}_s \ln \varrho - \mathbf{I}\mathbf{X}_s + \frac{\varrho}{2} \mathbf{U}_s - \frac{1}{2\mathfrak{h}_s} (\mathfrak{a}_s \mathbf{U}_s + \mathfrak{b}_s \mathbf{I}\mathbf{W}_s + \mathfrak{c}_s \mathbf{I}\mathbf{Z}_s) \varrho^2 \ln \varrho + \mathcal{O}(\varrho^2)\right) \\ &= \mathbf{U}_s + \frac{\mathbf{Y}_s}{\ln \varrho} - \frac{\mathbf{E} \cdot \mathbf{X}_s}{(\ln \varrho)^2} \mathbf{Y}_s + \frac{\varrho}{2(\ln \varrho)^2} \mathbf{Y}_s - \frac{(\mathbf{E} \cdot \mathbf{X}_s) \varrho}{(\ln \varrho)^3} \mathbf{I}\mathbf{X}_s + \mathcal{O}(\varrho^2), \end{aligned}$$

when $\varrho$ goes to zero, with

$$\mathbf{Y}_s := \mathbf{I}\mathbf{X}_s - (\mathbf{E} \cdot \mathbf{X}_s) \mathbf{U}_s. \tag{49}$$

Now that we have necessary pieces of notation, we gather in the following the behaviors found here above for $(k, \langle \underline{\mathbf{U}} \rangle)$ with the expansions proved below for $\alpha$.

**Corollary 4.** *Under Assumptions 1-2 and with notation from Theorem 2, (48) and (49) we have*
Harmonic limit
$$\begin{pmatrix} k \\ \alpha \\ \mathbf{M} \end{pmatrix} = \begin{pmatrix} k_0 \\ 0 \\ \mathbf{U}_0 \end{pmatrix} + \frac{\delta^2}{4\mathfrak{c}_0} \begin{pmatrix} -k_0 \mathfrak{a}_0 \\ \frac{1}{k_0} \mathbf{W}_0 \cdot \mathbb{S}^{-1} \mathbf{W}_0 \\ \mathbf{Y}_0 \end{pmatrix} + \mathcal{O}(\delta^4) \tag{50}$$
*when $\delta$ goes to zero.*
Soliton limit
$$\begin{pmatrix} k \\ \alpha \\ \mathbf{M} \end{pmatrix} = \begin{pmatrix} 0 \\ \partial_c \mathcal{M}(c; \mathbf{U}_s) \\ \mathbf{U}_s \end{pmatrix} + \frac{1}{\ln \varrho} \begin{pmatrix} -\frac{\pi}{\Xi_s} \\ \frac{\Xi_s}{\pi} \mathcal{Q}(\mathbf{Y}_s) \\ \mathbf{Y}_s \end{pmatrix} + \mathcal{O}\left(\frac{1}{(\ln \varrho)^2}\right) \tag{51}$$
*when $\varrho$ goes to zero.*

**Remark 11.** One can always look at the soliton limit as being the limit when $k$ goes to zero and Corollary 4 in particular contains in this regime
$$\alpha = \partial_c \mathcal{M}(c; \mathbf{U}_s) + \mathcal{O}(k), \qquad \langle \underline{\mathbf{U}} \rangle = \mathbf{U}_s + \mathcal{O}(k).$$
Likewise, assuming that
$$w_0 := \mathbf{W}_0 \cdot \mathbb{S}^{-1} \mathbf{W}_0 \tag{52}$$
is nonzero, we can equivalently look at the harmonic limit as the limit $\alpha$ goes to zero and then
$$k = k_0 + \mathcal{O}(\alpha), \qquad \langle \underline{\mathbf{U}} \rangle = \mathbf{U}_0 + \mathcal{O}(\alpha).$$

**Remark 12.** As already observed in [BGMR20], we can check in practical cases that $w_0 = \mathbf{W}_0 \cdot \mathbb{S}^{-1} \mathbf{W}_0 \neq 0$. Indeed, $\mathbf{W}_0 \cdot \mathbb{S}^{-1} \mathbf{W}_0 = 1/b$ in the case $N = 1$ and in the case $N = 2$, $\mathbf{W}_0 \cdot \mathbb{S}^{-1} \mathbf{W}_0 = 2\partial_v g(v_0)/b$ is nonzero both when $\tau$ is constant — which is the case for the Euler–Korteweg system in mass Lagrangian coordinates — and $c \neq 0$ and when $\tau$ is linear in $v$ — which is the case for the Euler–Korteweg system in Eulerian coordinates — and $\lambda_2 \neq 0$. We stress that the latter conditions are exactly the same conditions encountered in Proposition 2 where the sign of $\alpha$ was investigated. In particular as pointed out in Remark 5 both conditions are conjugated by the passage between mass Lagrangian and Eulerian formulations.

*Proof.* The only thing left is to expand our variable $\alpha$. In order to do so, by using (18) we can conveniently write it as
$$\alpha = -\mathbf{F} \cdot \nabla_{\mu,c,\boldsymbol{\lambda}} \Theta - \frac{\mathbf{I} \nabla_{\mu,c,\boldsymbol{\lambda}} \Theta}{2 \mathbf{E} \cdot \nabla_{\mu,c,\boldsymbol{\lambda}} \Theta} \cdot \mathbf{B}^{-1} \mathbf{I} \nabla_{\mu,c,\boldsymbol{\lambda}} \Theta.$$
From (38) and Proposition 37 we obtain
$$\begin{aligned}
\frac{\alpha}{\Xi_0} &= \mathcal{Q}(\mathbf{U}_0) + \frac{1}{4\mathfrak{c}_0} \left( \mathfrak{a}_0 \mathcal{Q}(\mathbf{U}_0) - \mathfrak{b}_0 \mathbf{F} \cdot \mathbf{W}_0 - \mathfrak{c}_0 \mathbf{F} \cdot \mathbf{Z}_0 \right) \delta^2 + \mathcal{O}(\delta^4) \\
&\quad - \frac{1}{1 + \frac{\mathfrak{a}_0}{4\mathfrak{c}_0} \delta^2 + \mathcal{O}(\delta^4)} \mathcal{Q}\left(\mathbf{U}_0 + \frac{1}{4\mathfrak{c}_0} (\mathfrak{a}_0 \mathbf{U}_0 + \mathfrak{b}_0 \mathbf{I} \mathbf{W}_0 + \mathfrak{c}_0 \mathbf{I} \mathbf{Z}_0) \delta^2 + \mathcal{O}(\delta^4)\right) \\
&= \mathcal{Q}(\mathbf{U}_0) + \frac{1}{4\mathfrak{c}_0} \left( \mathfrak{a}_0 \mathcal{Q}(\mathbf{U}_0) - \mathfrak{b}_0 \mathbf{F} \cdot \mathbf{W}_0 - \mathfrak{c}_0 \mathbf{F} \cdot \mathbf{Z}_0 \right) \delta^2 + \mathcal{O}(\delta^4) \\
&\quad - \frac{1}{1 + \frac{\mathfrak{a}_0}{4\mathfrak{c}_0} \delta^2 + \mathcal{O}(\delta^4)} \left( \mathcal{Q}(\mathbf{U}_0) + \mathbf{U}_0 \cdot \mathbf{B}^{-1} \frac{1}{4\mathfrak{c}_0} (\mathfrak{a}_0 \mathbf{U}_0 + \mathfrak{b}_0 \mathbf{I} \mathbf{W}_0 + \mathfrak{c}_0 \mathbf{I} \mathbf{Z}_0) \delta^2 + \mathcal{O}(\delta^4) \right)
\end{aligned}$$
in which there are some simplifications because by (37)
$$\begin{aligned}
\mathbf{F} \cdot \mathbf{W}_0 + \mathbf{U}_0 \mathbf{B}^{-1} \mathbf{I} \mathbf{W}_0 &= \mathbf{V}_0 \cdot \mathbb{S}^{-1} \mathbf{W}_0 + \mathcal{Q}(\mathbf{U}_0) \mathbf{E} \cdot \mathbf{W}_0 = 0, \\
\mathbf{F} \cdot \mathbf{Z}_0 + \mathbf{U}_0 \mathbf{B}^{-1} \mathbf{I} \mathbf{Z}_0 &= \mathbf{V}_0 \cdot \mathbb{S}^{-1} \mathbf{Z}_0 + \mathcal{Q}(\mathbf{U}_0) \mathbf{E} \cdot \mathbf{Z}_0 = -\mathbf{W}_0 \cdot \mathbb{S}^{-1} \mathbf{W}_0.
\end{aligned}$$

So we eventually find that
$$\alpha = \frac{\mathbf{W}_0 \cdot \mathbb{S}^{-1}\mathbf{W}_0}{4\mathfrak{c}_0\, k_0} \delta^2 + \mathcal{O}(\delta^4).$$

Likewise, from (40) and Proposition 37 we get
$$\frac{\pi}{\Xi_s}\alpha = -\mathcal{Q}(\mathbf{U}_s)\ln\varrho + \mathbf{F}\cdot\mathbf{X}_s + \mathcal{O}(\varrho)$$
$$+ \frac{\ln\varrho}{2} \frac{\mathbf{U}_s + \frac{\mathbf{IX}_s}{\ln\varrho} + \mathcal{O}\bigl(\frac{\varrho}{\ln\varrho}\bigr)}{1 + \frac{\mathbf{E}\cdot\mathbf{X}_s}{\ln\varrho} + \mathcal{O}\bigl(\frac{\varrho}{\ln\varrho}\bigr)} \cdot \mathbf{B}^{-1}\left(\mathbf{U}_s + \frac{\mathbf{IX}_s}{\ln\varrho} + \mathcal{O}\bigl(\frac{\varrho}{\ln\varrho}\bigr)\right),$$

which eventually simplifies into
$$\frac{\pi}{\Xi_s}\alpha = \mathbf{F}\cdot\mathbf{X}_s + \mathbf{U}_s\cdot\mathbf{B}^{-1}\mathbf{IX}_s - \mathcal{Q}(\mathbf{U}_s)(\mathbf{E}\cdot\mathbf{X}_s) + \frac{\mathcal{Q}(\mathbf{Y}_s)}{\ln\varrho} + \mathcal{O}\Bigl(\frac{1}{(\ln\varrho)^2}\Bigr),$$

or equivalently, since $\mathbf{F}\cdot\mathbf{X}_s + \mathbf{U}_s\cdot\mathbf{B}^{-1}\mathbf{IX}_s - \mathcal{Q}(\mathbf{U}_s)(\mathbf{E}\cdot\mathbf{X}_s) = (\mathbb{S}^{-1}\mathbf{V}_s)\cdot\mathbf{X}_s$,
$$\alpha = \partial_c \mathcal{M}(c;\mathbf{U}_s) + \frac{\Xi_s\mathcal{Q}(\mathbf{Y}_s)}{\pi\ln\varrho} + \mathcal{O}\Bigl(\frac{1}{(\ln\varrho)^2}\Bigr)$$

thanks to (42). $\square$

### 3.3 Extending the parametrization

We can even go further and show that $(k,\alpha,\mathbf{M})$ are 'good' variables up to the limits $k=0$ and $\alpha=0$.

**Theorem 3.** *Under Assumptions 1-2 and with notation from Theorem 2 we have*
<u>Harmonic limit</u> *The continuous extension of*
$$(\delta^2, c, \boldsymbol{\lambda}) \mapsto (k,\alpha,\mathbf{M})$$
*or equivalently of*
$$(\mu - \mu_0(c,\boldsymbol{\lambda}), c, \boldsymbol{\lambda}) \mapsto (k,\alpha,\mathbf{M})$$
*to $\{0\}\times\Lambda$ defines a $\mathcal{C}^1$ map in a connected open neighborhood (in $\mathbb{R}_+\times\Lambda$) of $\{0\}\times\Lambda$, which, provided that[12] $w_0$ does not vanish on $\Lambda$, is also a $\mathcal{C}^1$-diffeomorphism.*
<u>Soliton limit</u> *The continuous extension of the map*
$$(-\tfrac{1}{\ln\varrho}, c, \boldsymbol{\lambda}) \mapsto (k,\alpha,\mathbf{M})$$
*or equivalently of the map*
$$(-\tfrac{1}{\ln(\mu_s(c,\boldsymbol{\lambda})-\mu)}, c, \boldsymbol{\lambda}) \mapsto (k,\alpha,\mathbf{M})$$
*to $\{0\}\times\Lambda$ defines a $\mathcal{C}^1$ map in a connected open neighborhood (in $\mathbb{R}_+\times\Lambda$) of $\{0\}\times\Lambda$, which, provided that, for any $(c,\boldsymbol{\lambda})\in\Lambda$, $\partial_c^2\mathcal{M}(c;\mathbf{U}_s(c,\boldsymbol{\lambda}))\neq 0$, is also a $\mathcal{C}^1$-diffeomorphism.*

*Proof.* Expansions (50) and (51) show that the maps under consideration possess continuous extensions. To prove that these extensions are $\mathcal{C}^1$, we only need to prove that their Jacobian maps also extend continuously to $\Lambda\times\{0\}$. After that, by the Inverse Function Theorem, the proof will be achieved provided we also derive from extra assumptions that at any point of $\Lambda\times\{0\}$ the limit of the Jacobian map is nonsingular. In both limits our starting point is (31), that yields
$$\nabla_{\mu,c,\boldsymbol{\lambda}}\begin{pmatrix} k \\ \alpha \\ \mathbf{M} \end{pmatrix} = k\,(\nabla^2_{\mu,c,\boldsymbol{\lambda}}\Theta)\,(\mathbb{A}^{\mathsf{T}})^{-1}, \tag{53}$$

---
[12] See Remark 12.

with $\mathbb{A}$ as in Proposition 3.

In the harmonic limit we set $\epsilon := \delta^2$ and observe that the chain rule yields follows from

$$= k \left( \nabla_{\mu,c,\boldsymbol{\lambda}}\epsilon \;\middle|\; \begin{array}{c|c} 0 & 0 \\ 1 & 0 \\ 0 & \mathbf{I}_N \end{array} \right)^{-1} (\nabla^2_{\mu,c,\boldsymbol{\lambda}}\Theta)\,(\mathbb{A}^\mathsf{T})^{-1}.$$

Since $k$ has a nonzero limit $k_0$, it is trivial to check that both $k$ and $\mathbb{A}$ admit invertible limits when $\delta \to 0$. To deal with the factor involving $\nabla_{\mu,c,\boldsymbol{\lambda}}\epsilon$ we extract from [BGMR20, Proposition 4], expressed in our current notation,

$$\nabla_{\mu,c,\boldsymbol{\lambda}}\epsilon = 4\mathfrak{c}_0\mathbf{V}_0 + \mathcal{O}(\delta^2)$$

that implies readily that

$$\left( \nabla_{\mu,c,\boldsymbol{\lambda}}\epsilon \;\middle|\; \begin{array}{c|c} 0 & 0 \\ 1 & 0 \\ 0 & \mathbf{I}_N \end{array} \right)$$

possesses an invertible limit. At last, the fact that $\nabla^2_{\mu,c,\boldsymbol{\lambda}}\Theta$ possesses a limit when $\delta \to 0$ is a direct consequence of (39). The invertibility of the corresponding limit when $w_0$ (defined in (52)) is nonzero follows from straightforward computations based on the limit of $\mathbb{P}_0^\mathsf{T}(\nabla^2_{\mu,c,\boldsymbol{\lambda}}\Theta)\,\mathbb{P}_0$ obtained from Corollary 3.

More delicate is the soliton limit, in which all matrices involved in (53) blow up. To begin with we set $\epsilon := -1/\ln\varrho$ and extract from [BGMR20, Proposition 5]

$$\nabla_{\mu,c,\boldsymbol{\lambda}}\epsilon = -\frac{\mathfrak{h}_s}{(\varrho\ln\varrho)^2}((1+\tfrac{3}{2}\varrho)\mathbf{V}_s + \mathcal{O}(\varrho^2)),$$

To make the most of computations already carried out in Corollary 3, we use the factorization

$$\nabla_{\epsilon,c,\boldsymbol{\lambda}}\begin{pmatrix} k \\ \alpha \\ \mathbf{M} \end{pmatrix} = k\left( \mathbb{P}_s^\mathsf{T}\left( \nabla_{\mu,c,\boldsymbol{\lambda}}\epsilon \;\middle|\; \begin{array}{c|c} 0 & 0 \\ 1 & 0 \\ 0 & \mathbf{I}_N \end{array} \right) \right)^{-1} \mathbb{P}_s^\mathsf{T}(\nabla^2_{\mu,c,\boldsymbol{\lambda}}\Theta)\,\mathbb{P}_s\,((\mathbb{P}_s^\mathsf{T}\mathbb{A})^{-1})^\mathsf{T}. \tag{54}$$

stemming from the chain rule. To do so, first we observe that

$$\mathbb{P}_s^\mathsf{T}\left( \nabla_{\mu,c,\boldsymbol{\lambda}}\epsilon \;\middle|\; \begin{array}{c|c} 0 & 0 \\ 1 & 0 \\ 0 & \mathbf{I}_N \end{array} \right) = \left( \begin{array}{c|c} \frac{1}{(\varrho\ln\varrho)^2} & 0 \\ \hline 0 & \mathbf{I}_{N+1} \end{array} \right) \left( \begin{array}{c|c} -\mathfrak{h}_s + \mathcal{O}(\varrho) & 0 \\ \hline \mathcal{O}(\frac{1}{(\ln\varrho)^2}) & \mathbb{K}_s \end{array} \right) \tag{55}$$

with

$$\mathbb{K}_s = \left( \,0 \;\middle|\; \mathbf{I}_{N+1}\, \right)\mathbb{P}_s^\mathsf{T}\left( \begin{array}{c} 0 \\ \hline \mathbf{I}_{N+1} \end{array} \right)$$

easily seen to be invertible so that the last matrix in (55) possesses an invertible limit when $\varrho \to 0$. Now we stress that

$$\mathbb{P}_s^\mathsf{T}\mathbb{A} = \left( \begin{array}{c|c|c} 1/k & 0 & 0 \\ \hline 0 & k & 0 \\ \hline 0 & 0 & \mathbf{I}_N \end{array} \right) \left( \begin{array}{c|c|c} -1 & 0 & 0 \\ \hline \mathcal{Q}(\mathbf{U}_s - \mathbf{M})/k^2 & -1 & (\mathbf{U}_s - \mathbf{M})^\mathsf{T}\mathbf{B}^{-1}/k \\ \hline \mathbf{A}_s\mathbf{B}^{-1}(\mathbf{U}_s - \mathbf{M})/k & 0 & \mathbf{A}_s\mathbf{B}^{-1} \end{array} \right) \tag{56}$$

and that it follows from (51) that the last matrix in (56) possesses an invertible limit when $\varrho \to 0$.

Combining (54)-(55)-(56) with (51) reduces the issue to the inspection of the matrix

$$\mathbb{L} := \frac{1}{\ln(\varrho)} \begin{pmatrix} (\varrho \ln \varrho)^2 & 0 \\ \hline 0 & \mathbf{I}_{N+1} \end{pmatrix} \mathbb{P}_s{}^{\mathsf{T}} (\nabla^2_{\mu,c,\boldsymbol{\lambda}} \Theta) \mathbb{P}_s \begin{pmatrix} \frac{1}{\ln(\varrho)} & 0 & 0 \\ \hline 0 & \ln(\varrho) & 0 \\ \hline 0 & 0 & \mathbf{I}_N \end{pmatrix}. \qquad (57)$$

It follows from Theorem 2 and Corollary 3 that, when $N = 2$

$$\frac{\pi}{\Xi_s} \mathbb{L} = \begin{pmatrix} \mathfrak{h}_s & 0 & 0 & 0 \\ 0 & (\mathbb{S}^{-1}\mathbf{V}_s) \cdot \mathbb{O}_s \mathbb{S}^{-1}\mathbf{V}_s & 0 & 0 \\ 0 & (\mathbb{S}^{-1}\mathbf{T}_s) \cdot \mathbb{O}_s \mathbb{S}^{-1}\mathbf{V}_s & 2\mathfrak{c}_s \sigma_s^2 & 2\mathfrak{c}_s \sigma_s w_s \\ 0 & (\mathbb{S}^{-1}\mathbf{W}_s) \cdot \mathbb{O}_s \mathbb{S}^{-1}\mathbf{V}_s & 2\mathfrak{c}_s \sigma_s w_s & 2\mathfrak{c}_s w_s^2 + \sigma_s^2 \end{pmatrix} + \mathcal{O}\left(\frac{1}{(\ln \varrho)}\right)$$

while when $N = 1$

$$\frac{\pi}{\Xi_s} \mathbb{L} = \begin{pmatrix} \mathfrak{h}_s & 0 & 0 \\ 0 & (\mathbb{S}^{-1}\mathbf{V}_s) \cdot \mathbb{O}_s \mathbb{S}^{-1}\mathbf{V}_s & 0 \\ 0 & (\mathbb{S}^{-1}\mathbf{W}_s) \cdot \mathbb{O}_s \mathbb{S}^{-1}\mathbf{V}_s & 2\mathfrak{c}_s b^{-2} \end{pmatrix} + \mathcal{O}\left(\frac{1}{(\ln \varrho)}\right).$$

This implies that $\mathbb{L}$ possesses a limit when $\varrho \to 0$ and that this limit is invertible provided that $(\mathbb{S}^{-1}\mathbf{V}_s) \cdot \mathbb{O}_s \mathbb{S}^{-1}\mathbf{V}_s \neq 0$, and finishes the proof. $\square$

# 4 Asymptotics of the modulation system

## 4.1 Extending the averaged Hamiltonian

Our goal is now to show that the averaged Hamiltonian H extends as a $\mathcal{C}^2$ function of $(k, \alpha, \mathbf{M})$ both to the zero-amplitude regime $\alpha = 0$ and to the zero-wavenumber regime $k = 0$.

Under natural assumptions required by Theorem 1, it is quite elementary, by using the definition of H and relations (21), to check that H does extend as a $\mathcal{C}^1$ map both to $\alpha = 0$ and to $k = 0$. This is already sufficient to take the relevant limits of the conservative form (22) of the modulated system. Yet to ensure that hyperbolic properties of the limiting system do transfer to the original ones in relevant regimes one needs to be able to take limits in the quasilinear form (28) hence to prove the $\mathcal{C}^2$ extension property we discuss now.

To state the following theorem in a precise way, let us denote, in the harmonic limit, as $\Lambda_0$ the image of $\Lambda$ by $(c, \boldsymbol{\lambda}) \mapsto (k_0, 0, \mathbf{U}_0)$ and, in the soliton limit, as $\Lambda_s$ the image of $\Lambda$ by $(c, \boldsymbol{\lambda}) \mapsto (0, \partial_c \mathcal{M}(c; \mathbf{U}_s), \mathbf{U}_s)$.

**Theorem 4.** *Under Assumptions 1-2 and with notation from Theorem 2 we have*
<u>Harmonic limit</u> *Provided that $w_0$ does not vanish on $\Lambda$,*
*the averaged Hamiltonian H extends as a $\mathcal{C}^2$ function of $(k, \alpha, \mathbf{M})$ to a connected open neighborhood (in[13] $\mathbb{R} \times (\text{sign}(w_0)\,\mathbb{R}_+) \times \mathbb{R})$ of $\Lambda_0$.*
<u>Soliton limit</u> *Provided that, for any $(c, \boldsymbol{\lambda}) \in \Lambda$, $\partial_c^2 \mathcal{M}(c; \mathbf{U}_s(c, \boldsymbol{\lambda})) \neq 0$,*
*the averaged Hamiltonian H extends as a $\mathcal{C}^2$ function of $(k, \alpha, \mathbf{M})$ to a connected open neighborhood (in $\mathbb{R}_+ \times \mathbb{R} \times \mathbb{R})$ of $\Lambda_s$.*

*Proof.* The proof is similar to the one of Theorem 3. In particular the issue is readily reduced to checking that the assumptions of Theorem 4 ensure that in the relevant regimes $\nabla^2_{\mu,c,\boldsymbol{\lambda}} \Theta$ is invertible and $\nabla^2_{k,\alpha,\mathbf{M}} \text{H}$ possesses a limit, the study of the latter relying on (29).

In the harmonic limit, we have already checked all the required claims along the proof of Theorem 3 since there we have checked that $\nabla^2_{\mu,c,\boldsymbol{\lambda}} \Theta$ possesses an invertible limit.

---
[13]Note that since $\Lambda$ is connected, $w_0$ has a definite sign on $\Lambda_0$.

The soliton limit requires slightly more work. We already know from the proof of Theorem 3 that both

$$\left(\begin{array}{c|c} (\varrho \ln \varrho)^2 & 0 \\ \hline 0 & \mathbf{I}_{N+1} \end{array}\right) k \, \mathbb{P}_s{}^\mathsf{T} (\nabla^2_{\mu,c,\boldsymbol{\lambda}} \Theta) (\mathbb{A}^\mathsf{T})^{-1} \quad \text{and} \quad \left(\begin{array}{c|c|c} 1/\ln(\varrho) & 0 & 0 \\ \hline 0 & \ln(\varrho) & 0 \\ \hline 0 & 0 & \mathbf{I}_N \end{array}\right) \mathbb{P}_s{}^\mathsf{T} \mathbb{A}$$

possess invertible limits. Thus the result stems from

$$\nabla^2_{k,\alpha,\mathbf{M}} \mathrm{H} + c\, \mathbb{B}^{-1}$$

$$= - \left( \left(\begin{array}{c|c} (\varrho \ln \varrho)^2 & 0 \\ \hline 0 & \mathbf{I}_{N+1} \end{array}\right) k\, \mathbb{P}_s{}^\mathsf{T} \nabla^2_{\mu,c,\boldsymbol{\lambda}} \Theta\, (\mathbb{A}^\mathsf{T})^{-1} \right)^{-1}$$

$$\times \left(\begin{array}{c|c|c} \varrho^2(\ln \varrho)^3 & 0 & 0 \\ \hline 0 & 1/\ln(\varrho) & 0 \\ \hline 0 & 0 & \mathbf{I}_N \end{array}\right) \times \left(\begin{array}{c|c|c} 1/\ln(\varrho) & 0 & 0 \\ \hline 0 & \ln(\varrho) & 0 \\ \hline 0 & 0 & \mathbf{I}_N \end{array}\right) \mathbb{P}_s{}^\mathsf{T} \mathbb{A}$$

derived from (29). □

## 4.2 Basic features of the limiting modulated systems

A detailed inspection of the proof of Theorem 4 provides explicit formulas for limiting values of $\nabla^2_{k,\alpha,\mathbf{M}} \mathrm{H}$ thus of the Whitham matrix $\mathbb{W} = -\mathbb{B}\, \nabla^2_{k,\alpha,\mathbf{M}} \mathrm{H}$ in terms of coefficients from Theorem 2. Yet first we restrain from giving these and focus instead on what can be derived from more elementary arguments, using only the conclusion from Theorem 4, that is, $\mathcal{C}^2$ regularity of H. Along the discussion we shall still denote as H the extension of H to either $\Lambda_0$ or $\Lambda_s$.

To do so, we first point out the elementary

$$\partial_k \mathrm{H} = \Theta - \alpha c = 0 \quad \text{on } \Lambda_0, \qquad \partial_\alpha \mathrm{H} = -kc = 0 \quad \text{on } \Lambda_s,$$

which by differentiating tangentially yield

$$\partial_k^2 \mathrm{H} = 0 \quad \text{and} \quad \partial_k \nabla_\mathbf{M} \mathrm{H} = 0 \quad \text{on } \Lambda_0, \qquad \partial_\alpha^2 \mathrm{H} = 0 \quad \text{and} \quad \partial_\alpha \nabla_\mathbf{M} \mathrm{H} = 0 \quad \text{on } \Lambda_s.$$

In particular,

$$\mathbb{B}\, \nabla^2_{k,\alpha,\mathbf{M}} \mathrm{H} = \left(\begin{array}{cc|c} \partial^2_{k\alpha} \mathrm{H} & \partial^2_\alpha \mathrm{H} & \partial_\alpha \nabla_\mathbf{M} \mathrm{H}^\mathsf{T} \\ 0 & \partial^2_{k\alpha} \mathrm{H} & 0 \\ \hline 0 & \mathbf{B} \partial_\alpha \nabla_\mathbf{M} \mathrm{H} & \mathbf{B} \nabla^2_\mathbf{M} \mathrm{H} \end{array}\right) \quad \text{on } \Lambda_0,$$

$$\mathbb{B}\, \nabla^2_{k,\alpha,\mathbf{M}} \mathrm{H} = \left(\begin{array}{cc|c} \partial^2_{k\alpha} \mathrm{H} & 0 & 0 \\ \partial^2_k \mathrm{H} & \partial^2_{k\alpha} \mathrm{H} & \partial_k \nabla_\mathbf{M} \mathrm{H}^\mathsf{T} \\ \hline \mathbf{B} \partial_k \nabla_\mathbf{M} \mathrm{H} & 0 & \mathbf{B} \nabla^2_\mathbf{M} \mathrm{H} \end{array}\right) \quad \text{on } \Lambda_s,$$

and a direct computation of a characteristic polynomial shows that, on either $\Lambda_0$ or $\Lambda_s$, the spectrum of $\mathbb{W}$ is the union, with multiplicity, of the spectrum of $-\mathbf{B} \nabla^2_\mathbf{M} \mathrm{H}$ and twice $-\partial^2_{k\alpha} \mathrm{H}$.

The fact that some second order derivatives of H are easier to compute is no accident. Since it is easy to extend H as a $\mathcal{C}^1$ map it is also straightforward to extend its second-order derivatives that contain at most one normal derivative. The hard parts of Theorem 4 are the extensions of $\partial_\alpha^2 \mathrm{H}$ to $\Lambda_0$ and, even more, of $\partial_k^2 \mathrm{H}$ to $\Lambda_s$. To illustrate this further let us stress that for $i \in \{0, s\}$

$$\mathrm{H}(k, \alpha, \mathbf{U}_i) = \mathcal{H}(\mathbf{U}_i, 0) \quad \text{on } \Lambda_i, \qquad \text{thus} \qquad \nabla_\mathbf{M}^2 \mathrm{H}(k, \alpha, \mathbf{U}_i) = \nabla_\mathbf{U}^2 \mathcal{H}(\mathbf{U}_i, 0) \quad \text{on } \Lambda_i$$

so that, on $\Lambda_i$, $-\mathbb{B}\nabla_\mathbf{M}^2 \mathrm{H}(k, \alpha, \mathbf{U}_i)$ is the characteristic matrix at $\mathbf{U}_i$ of the dispersionless system (7). Likewise, for any $(k_0, \mathbf{U}_0)$ such that $(k_0, 0, \mathbf{U}_0) \in \Lambda_0$

$$\begin{aligned}
\partial_{k\alpha}^2 \mathrm{H}(k_0, 0, \mathbf{U}_0) &= -\partial_k(k\, c_0)(k_0, \mathbf{U}_0) = -v_g(k_0, \mathbf{U}_0)\,, \\
\partial_\alpha \nabla_\mathbf{M} \mathrm{H}(k_0, 0, \mathbf{U}_0) &= -k_0 \nabla_\mathbf{U} c_0(k_0, \mathbf{U}_0)\,,
\end{aligned}$$

where $v_g(k_0, \mathbf{U}_0)$ is the linear group velocity and $c_0(k_0, \mathbf{U}_0) = c(k, 0, \mathbf{U}_0)$ is the harmonic phase velocity of the harmonic wavetrain on $\mathbf{U}_0$ at wavenumber $k_0$. Similarly, at the soliton limit, we have for any $(c_s, \mathbf{U}_s) \in \Lambda$

$$\partial_k \mathrm{H}(0, \partial_c \mathcal{M}(c_s; \mathbf{U}_s), \mathbf{U}_s) = \mathcal{M}(c_s; \mathbf{U}_s) - c_s\, \partial_c \mathcal{M}(c_s; \mathbf{U}_s)$$

so that for any $(c_s, \mathbf{U}_s) \in \Lambda$

$$\begin{aligned}
\partial_{k\alpha}^2 \mathrm{H}(0, \partial_c \mathcal{M}(c_s; \mathbf{U}_s), \mathbf{U}_s) &= -c_s\,, \\
\partial_k \nabla_\mathbf{M} \mathrm{H}(0, \partial_c \mathcal{M}(c_s; \mathbf{U}_s), \mathbf{U}_s) &= \nabla_\mathbf{U} \mathcal{M}(c_s; \mathbf{U}_s)\,.
\end{aligned}$$

Going back to the cancellations in $\mathbb{B}\, \nabla_{k,\alpha,\mathbf{M}}^2 \mathrm{H}$, we make the following elementary algebraic observation, whose proof follows from a short computation — left to the reader.

**Lemma 1.** *For any real numbers $v$ and $a$, any vectors $\mathbf{l}_0$, $\mathbf{r}_0$ in $\mathbb{R}^N$ and any $N \times N$ matrix $\mathbb{M}$ such that $v$ is not an eigenvalue of $\mathbb{M}$ we have*

$$\left(\begin{array}{cc|c} 1 & 0 & \mathbf{l}^\mathsf{T} \\ 0 & 1 & 0 \\ \hline 0 & \mathbf{r} & \mathbf{I}_N \end{array}\right)^{-1} \left(\begin{array}{cc|c} v & a & \mathbf{l}_0^\mathsf{T} \\ 0 & v & 0 \\ \hline 0 & \mathbf{r}_0 & \mathbb{M} \end{array}\right) \left(\begin{array}{cc|c} 1 & 0 & \mathbf{l}^\mathsf{T} \\ 0 & 1 & 0 \\ \hline 0 & \mathbf{r} & \mathbf{I}_N \end{array}\right) = \left(\begin{array}{cc|c} v & a - \mathbf{l}_0^\mathsf{T}(\mathbb{M} - v\,\mathbf{I}_N)^{-1}\mathbf{r}_0 & 0 \\ 0 & v & 0 \\ \hline 0 & 0 & \mathbb{M} \end{array}\right)$$

$$\left(\begin{array}{cc|c} 1 & 0 & 0 \\ 0 & 1 & \mathbf{l}^\mathsf{T} \\ \hline \mathbf{r} & 0 & \mathbf{I}_N \end{array}\right)^{-1} \left(\begin{array}{cc|c} v & 0 & 0 \\ a & v & \mathbf{l}_0^\mathsf{T} \\ \hline \mathbf{r}_0 & 0 & \mathbb{M} \end{array}\right) \left(\begin{array}{cc|c} 1 & 0 & 0 \\ 0 & 1 & \mathbf{l}^\mathsf{T} \\ \hline \mathbf{r} & 0 & \mathbf{I}_N \end{array}\right) = \left(\begin{array}{cc|c} v & 0 & 0 \\ a - \mathbf{l}_0^\mathsf{T}(\mathbb{M} - v\,\mathbf{I}_N)^{-1}\mathbf{r}_0 & v & 0 \\ \hline 0 & 0 & \mathbb{M} \end{array}\right)$$

*with*

$$\mathbf{l} := (\mathbb{M}^\mathsf{T} - v\,\mathbf{I}_N)^{-1}\mathbf{l}_0\,, \qquad \mathbf{r} := -(\mathbb{M} - v\,\mathbf{I}_N)^{-1}\mathbf{r}_0\,.$$

Note that the two matrices considered in Lemma 1 are obtained one from the other merely by exchanging the first and second coordinates. We have introduced these two cases just to emphasize that this algebraic lemma applies to both kinds of limit.

For convenience, let us summarize part of the foregoing findings in the following statement.

**Corollary 5.** *Under the assumptions of Theorem 4, still denoting by H its extension to either $\Lambda_0$ or $\Lambda_s$, we have*
<u>Harmonic limit</u> *At any point $(k_0, 0, \mathbf{U}_0)$ of $\Lambda_0$, the spectrum of the characteristic matrix of the modulation system $-\mathbb{B}\,\nabla_{k,\alpha,\mathbf{M}}^2 \mathrm{H}(k_0, 0, \mathbf{U}_0)$ is given, with algebraic multiplicity, by the spectrum of the dispersionless characteristic matrix $-\mathbb{B}\,\nabla_\mathbf{U}^2 \mathcal{H}(\mathbf{U}_0, 0)$ and twice the linear group velocity $v_g(k_0, \mathbf{U}_0) = -\partial_{k\alpha}^2 \mathrm{H}(k_0, 0, \mathbf{U}_0)$, so that in particular the modulation system is weakly hyperbolic if and only if the dispersionless system is so.*

*Moreover, $v_g(k_0, \mathbf{U}_0)$ is a semisimple characteristic of (22) if and only if $\partial_\alpha^2 H(k_0, 0, \mathbf{U}_0)$ coincides with[14]*

$$(\partial_\alpha \nabla_\mathbf{M} H)^\mathsf{T} (\nabla_\mathbf{M}^2 H - (\partial_{k\alpha}^2 H) \mathbf{B}^{-1})^{-1} \partial_\alpha \nabla_\mathbf{M} H$$
$$= k_0^2 (\nabla_\mathbf{U} c_0)^\mathsf{T} (\nabla_\mathbf{U}^2 \mathcal{H}(\mathbf{U}_0, 0) + v_g \mathbf{B}^{-1})^{-1} \nabla_\mathbf{U} c_0$$

*so that the modulation system is hyperbolic if and only if the foregoing condition is satisfied and the dispersionless system is hyperbolic.*

<u>Soliton limit</u> *For any $(c_s, \mathbf{U}_s)$ of $\Lambda$, the spectrum of the characteristic matrix of the modulation system $-\mathbf{B}\,\nabla_{k,\alpha,\mathbf{M}}^2 H(0, \partial_c \mathcal{M}(c_s; \mathbf{U}_s), \mathbf{U}_s)$ is given, with algebraic multiplicity, by the spectrum of the dispersionless characteristic matrix $-\mathbf{B}\,\nabla_\mathbf{U}^2 \mathcal{H}(\mathbf{U}_s, 0)$ and twice the soliton velocity $c_s = -\partial_{k\alpha}^2 H(0, \partial_c \mathcal{M}(c_s; \mathbf{U}_s), \mathbf{U}_s)$, so that in particular the modulation system is always weakly hyperbolic.*

*Moreover, $c_s$ is a semisimple characteristic of (22) if and only if $\partial_k^2 H(0, \partial_c \mathcal{M}(c_s; \mathbf{U}_s), \mathbf{U}_s)$ coincides with[15]*

$$(\partial_k \nabla_\mathbf{M} H)^\mathsf{T} (\nabla_\mathbf{M}^2 H - (\partial_{k\alpha}^2 H) \mathbf{B}^{-1})^{-1} \partial_k \nabla_\mathbf{M} H$$
$$= (\nabla_\mathbf{U} \mathcal{M})^\mathsf{T} (\nabla_\mathbf{U}^2 \mathcal{H}(\mathbf{U}_s, 0) + c_s \mathbf{B}^{-1})^{-1} \nabla_\mathbf{U} \mathcal{M}$$

*so that the modulation system is hyperbolic if and only if the foregoing condition is satisfied.*

*Proof.* At the soliton limit, the only thing left is to check that from the assumptions of Theorem 4 stem that $c_s$ is not an eigenvalue of $-\mathbf{B}\,\nabla_\mathbf{U}^2 \mathcal{H}(\mathbf{U}_s, 0)$ and that eigenvalues of the latter matrix are real and distinct. Yet a relatively direct computation (for which the reader is referred to Appendix A and [BGMR20, Appendix A]) shows that

$$\det(\mathbf{B}\,\nabla_\mathbf{U}^2 \mathcal{H}(\mathbf{U}_s, 0) + c_s\,\mathbf{I}_N) = \begin{cases} b\,\partial_v^2 \mathcal{W}(v_s; c_s, \boldsymbol{\lambda}_s) & \text{if} \quad N = 1 \\ b^2\,\tau(v_s)\,\partial_v^2 \mathcal{W}(v_s; c_s, \boldsymbol{\lambda}_s) & \text{if} \quad N = 2 \end{cases}$$

so that the conditions stem from Assumption 2 that contains $\partial_v^2 \mathcal{W}(v_s; c_s, \boldsymbol{\lambda}_s) < 0$.

At the harmonic limit, we only need to check that from the assumptions of Theorem 4 stems that $v_g$ is not an eigenvalue of $-\mathbf{B}\,\nabla_\mathbf{U}^2 \mathcal{H}(\mathbf{U}_0, 0)$. We first stress that the relation pointed out above also holds for $\mathbf{U}_0$ (instead of $\mathbf{U}_s$) so that

$$\det(\mathbf{B}\,\nabla_\mathbf{U}^2 \mathcal{H}(\mathbf{U}_0, 0) + c_0(k_0, \mathbf{U}_0)\,\mathbf{I}_N) = \begin{cases} b\,(2\pi)^2 k_0^2\,\kappa(v_0) & \text{if} \quad N = 1 \\ b^2\,\tau(v_0)\,(2\pi)^2 k_0^2\,\kappa(v_0) & \text{if} \quad N = 2 \end{cases}.$$

Multiplying first the latter by $k_0$ then differentiating it with respect to $k_0$ yield by the $N$-linearity of the determinant

$$\det(\mathbf{B}\,\nabla_\mathbf{U}^2 \mathcal{H}(\mathbf{U}_0, 0) + v_g(k_0, \mathbf{U}_0)\,\mathbf{I}_N)$$
$$= \begin{cases} 3\,b\,(2\pi)^2 k_0^2\,\kappa(v_0) & \text{if} \quad N = 1 \\ 3\,b^2\,\tau(v_0)\,(2\pi)^2 k_0^2\,\kappa(v_0) + (k_0 \partial_k c_0(k_0, \mathbf{U}_0))^2 & \text{if} \quad N = 2 \end{cases}.$$

This proves the claim. □

**Remark 13.** Concerning the case $N = 2$, at the harmonic limit, note that $\mathbf{B}\,\nabla_\mathbf{U}^2 \mathcal{H}(\mathbf{U}_0, 0)$ is never diagonal and that

$$(\mathrm{tr}(\mathbf{B}\,\nabla_\mathbf{U}^2 \mathcal{H}(\mathbf{U}_0, 0)))^2 - 4\,\det(\mathbf{B}\,\nabla_\mathbf{U}^2 \mathcal{H}(\mathbf{U}_0, 0)) = 4\,b^2 \tau(v_0)\left(f''(v_0) + \frac{1}{2}\tau''(v_0)\,u_0^2\right)$$

so that the dispersionless system is weakly hyperbolic (resp. hyperbolic) at $\mathbf{U}_0$ if and only if $f''(v_0) + \frac{1}{2}\tau''(v_0)\,u_0^2 \geq 0$ (resp. $f''(v_0) + \frac{1}{2}\tau''(v_0)\,u_0^2 > 0$). In particular when $\tau$ is affine, as is the case for Euler-Korteweg systems, this condition reduces to the requirement that $f$ be convex, which is the usual hyperbolicity condition for the Euler systems in terms of pressure monotonicity.

---

[14] The left-hand side being evaluated at $(k_0, 0, \mathbf{U}_0)$ and the right-hand side at $(k_0, \mathbf{U}_0)$.
[15] The left-hand side being evaluated at $(0, \partial_c \mathcal{M}(c_s; \mathbf{U}_s), \mathbf{U}_s)$ and the right-hand side at $(c_s, \mathbf{U}_s)$.

Despite the symmetry of Corollary 5 with respect to permutation of the variables $k$ and $\alpha$, the harmonic and soliton limits differ significantly in terms of the hyperbolic nature of the limiting system. Indeed as it follows from the analysis expounded in next subsection, the condition of Corollary 5 ensuring the semisimplicity of the characteristic value $-\partial^2_{k\alpha} H$ is always satisfied at the soliton limit whereas in general it fails at the harmonic limit. In particular as we show in Appendix A the latter condition does fail for the classical Korteweg–de Vries equation.

We stress however that in both cases limiting systems contain two subsystems of size $3 \times 3$ whose hyperbolicity is directly determined by the hyperbolicity of the dispersionless system at the limiting constant. On one hand, as expected, in the harmonic limit (resp. the soliton limit) the limiting modulation system leaves invariant the hyperplane $\alpha = 0$ (resp. $k = 0$) and its restriction fits the claim. On the other hand, in the harmonic limit (resp. the soliton limit) one may leave out the $\partial_t k$ equation (resp. $\partial_t \alpha$ equation) and obtain a closed system in $(\alpha, \mathbf{M})$ (resp. $(k, \mathbf{M})$) that also fits the claim. We warn the reader that a careless[16] choice of variables, effectively loosing one dimension at the limit, could suggest as a limiting system a system apparently of size $4 \times 4$ but actually reducible to the latter $3 \times 3$ system and hence fail to capture losses of hyperbolicity at the limit.

Note however that the direct consequences of the above discrepancy on the original modulation systems (and not their limiting extensions) are almost immaterial. Indeed whereas the failure of weak hyperbolicity (as potentially caused here by the failure of weak hyperbolicity of the dispersionless system) is stable under perturbation, neither hyperbolicity nor failure of hyperbolicity are stable phenomena in the presence of a multiple root. The determination of the nature of the original modulation systems will require an even finer analysis than the one carried out in next subsection.

The reader is referred to Appendix B for further comments on the foregoing phenomena, with fully workout details illustrated on a $2 \times 2$ model system.

## 4.3 Explicit formulas for the limiting modulation systems

Now, to push our analysis a bit further, we extract from the proof of Theorem 4 explicit formulas for the limiting values of $\nabla^2_{k,\alpha,\mathbf{M}} H$, in particular for $\partial^2_\alpha H$ in the harmonic limit and for $\partial^2_k H$ in the soliton limit.

Let us begin with the harmonic limit. For concision's sake we first introduce

$$\Sigma_0 = \begin{cases} \begin{pmatrix} 2\mathfrak{c}_0\, b^{-2} \end{pmatrix} & \text{if } N = 1 \\ \begin{pmatrix} 2\mathfrak{c}_0\, \sigma_0^2 & 2\mathfrak{c}_0\, \sigma_0\, w_0 \\ 2\mathfrak{c}_0\, \sigma_0\, w_0 & 2\mathfrak{c}_0\, w_0^2 - \sigma_0^2 \end{pmatrix} & \text{if } N = 2 \end{cases}$$

and

$$\mathbf{x}_0 = \begin{cases} \begin{pmatrix} \mathfrak{b}_0\, b^{-1} \end{pmatrix} & \text{if } N = 1 \\ \begin{pmatrix} \mathfrak{b}_0\, \sigma_0 \\ \mathfrak{b}_0\, w_0 + \mathfrak{c}_0\, \zeta_0 \end{pmatrix} & \text{if } N = 2 \end{cases}$$

so that it follows from Theorem 2 and Corollary 3 that

$$k\, \mathbb{P}_0^\mathsf{T} (\nabla^2_{\mu,c,\boldsymbol{\lambda}} \Theta)\, (\mathbb{A}^\mathsf{T})^{-1} \mathbb{B}^{-1} \;=\; \left( \begin{array}{cc|c} \mathfrak{c}_0 w_0\, \Xi_0 & -\dfrac{\mathfrak{a}_0}{\Xi_0} & \mathbf{x}_0^\mathsf{T} (\mathbf{A}_0^\mathsf{T})^{-1} \\ 0 & \dfrac{\mathfrak{c}_0 w_0}{\Xi_0} & 0 \\ \hline 0 & -\dfrac{\mathbf{x}_0}{\Xi_0} & \Sigma_0 (\mathbf{A}_0^\mathsf{T})^{-1} \end{array} \right) + \mathcal{O}(\delta^2)$$

---

[16]In contrast, Theorem 3 ensures that the set of variables $(k, \alpha, \mathbf{M})$ does not suffer from such a flaw.

thus its inverse is

$$\left(\begin{array}{cc|c} \frac{1}{\mathfrak{c}_0 w_0 \Xi_0} & \frac{1}{\mathfrak{c}_0^2 w_0^2 \Xi_0} \left(\mathfrak{a}_0 - \mathbf{x}_0{}^\mathsf{T}(\Sigma_0)^{-1}\mathbf{x}_0\right) & -\frac{1}{\mathfrak{c}_0 w_0 \Xi_0}\mathbf{x}_0{}^\mathsf{T}(\Sigma_0)^{-1} \\ 0 & \frac{\Xi_0}{\mathfrak{c}_0 w_0} & 0 \\ \hline 0 & \frac{1}{\mathfrak{c}_0 w_0}\mathbf{A}_0{}^\mathsf{T}(\Sigma_0)^{-1}\mathbf{x}_0 & \mathbf{A}_0{}^\mathsf{T}(\Sigma_0)^{-1} \end{array}\right) + \mathcal{O}(\delta^2).$$

From this stems that $\frac{1}{k}\mathbb{B}\mathbb{A}^\mathsf{T}(\nabla^2_{\mu,c,\boldsymbol{\lambda}}\Theta)^{-1}\mathbb{A}$ equals

$$\left(\begin{array}{cc|c} -\frac{1}{\mathfrak{c}_0 w_0} & -\frac{1}{\mathfrak{c}_0^2 w_0^2 \Xi_0^2}\left(\mathfrak{a}_0 - \mathbf{x}_0{}^\mathsf{T}(\Sigma_0)^{-1}\mathbf{x}_0\right) & -\frac{1}{\mathfrak{c}_0 w_0 \Xi_0}\mathbf{x}_0{}^\mathsf{T}(\Sigma_0)^{-1}\mathbf{A}_0\mathbf{B}^{-1} \\ 0 & -\frac{1}{\mathfrak{c}_0 w_0} & 0 \\ \hline 0 & -\frac{1}{\mathfrak{c}_0 w_0 \Xi_0}\mathbf{A}_0{}^\mathsf{T}(\Sigma_0)^{-1}\mathbf{x}_0 & \mathbf{A}_0{}^\mathsf{T}(\Sigma_0)^{-1}\mathbf{A}_0\mathbf{B}^{-1} \end{array}\right) + \mathcal{O}(\delta^2).$$

In particular, on $\Lambda_0$,

$$\nabla^2_{k,\alpha,\mathbf{M}}\mathrm{H} =$$

$$\left(\begin{array}{cc|c} 0 & -c + \frac{1}{\mathfrak{c}_0 w_0} & 0 \\ -c + \frac{1}{\mathfrak{c}_0 w_0} & \frac{1}{\mathfrak{c}_0^2 w_0^2 \Xi_0^2}\left(\mathfrak{a}_0 - \mathbf{x}_0{}^\mathsf{T}(\Sigma_0)^{-1}\mathbf{x}_0\right) & \frac{1}{\mathfrak{c}_0 w_0 \Xi_0}\mathbf{x}_0{}^\mathsf{T}(\Sigma_0)^{-1}\mathbf{A}_0\mathbf{B}^{-1} \\ \hline 0 & \frac{1}{\mathfrak{c}_0 w_0 \Xi_0}\mathbf{B}^{-1}\mathbf{A}_0{}^\mathsf{T}(\Sigma_0)^{-1}\mathbf{x}_0 & -c\mathbf{B}^{-1} - \mathbf{B}^{-1}\mathbf{A}_0{}^\mathsf{T}(\Sigma_0)^{-1}\mathbf{A}_0\mathbf{B}^{-1} \end{array}\right),$$

$$\mathbb{W} = \left(\begin{array}{cc|c} c - \frac{1}{\mathfrak{c}_0 w_0} & -\frac{1}{\mathfrak{c}_0^2 w_0^2 \Xi_0^2}\left(\mathfrak{a}_0 - \mathbf{x}_0{}^\mathsf{T}(\Sigma_0)^{-1}\mathbf{x}_0\right) & -\frac{1}{\mathfrak{c}_0 w_0 \Xi_0}\mathbf{x}_0{}^\mathsf{T}(\Sigma_0)^{-1}\mathbf{A}_0\mathbf{B}^{-1} \\ 0 & c - \frac{1}{\mathfrak{c}_0 w_0} & 0 \\ \hline 0 & -\frac{1}{\mathfrak{c}_0 w_0 \Xi_0}\mathbf{A}_0{}^\mathsf{T}(\Sigma_0)^{-1}\mathbf{x}_0 & c\mathbf{I}_N + \mathbf{A}_0{}^\mathsf{T}(\Sigma_0)^{-1}\mathbf{A}_0\mathbf{B}^{-1} \end{array}\right).$$

Translating, by identification, the foregoing computations into the notation of Corollary 5 yields the following result.

**Theorem 5.** *Under the assumptions of Theorem 4, let us still denote by* $\mathrm{H}$ *and* $\mathbb{W}$ *their extensions to*

$\Lambda_0$. Then at any point $(k_0, 0, \mathbf{U}_0)$ of $\Lambda_0$, we have[17]

$$\partial_\alpha^2 H = k_0^4 (\partial_k c_0)^2 \, \mathfrak{a}_0 + k_0^2 \, \nabla_{\mathbf{U}} c_0^{\mathsf{T}} \, (\nabla_{\mathbf{U}}^2 \mathcal{H}(\mathbf{U}_0, 0) + c_0 \mathbf{B}^{-1})^{-1} \nabla_{\mathbf{U}} c_0$$

$$\nabla^2_{k,\alpha,\mathbf{M}} H = \left( \begin{array}{cc|c} 0 & -c_0 - k_0 \partial_k c_0 & 0 \\ -c_0 - k_0 \partial_k c_0 & \partial_\alpha^2 H & -k_0 \, \nabla_{\mathbf{U}} c_0^{\mathsf{T}} \\ \hline 0 & -k_0 \, \nabla_{\mathbf{U}} c_0 & \nabla_{\mathbf{U}}^2 \mathcal{H}(\mathbf{U}_0, 0) \end{array} \right)$$

$$\mathbb{W} = \left( \begin{array}{cc|c} c_0 + k_0 \partial_k c_0 & -\partial_\alpha^2 H & k_0 \, \nabla_{\mathbf{U}} c_0^{\mathsf{T}} \\ 0 & c_0 + k_0 \partial_k c_0 & 0 \\ \hline 0 & k_0 \, \mathbf{B} \nabla_{\mathbf{U}} c_0 & -\mathbf{B} \nabla_{\mathbf{U}}^2 \mathcal{H}(\mathbf{U}_0, 0) \end{array} \right)$$

where $\mathfrak{a}_0$ is as in Theorem 2

$$\mathfrak{a}_0 = k_0 \, \partial_\mu^2 \Theta(\mu_0, c_0, \boldsymbol{\lambda}_0) \,.$$

In particular, at any $(k_0, 0, \mathbf{U}_0)$ of $\Lambda_0$, denoting $v_g = c_0 + k_0 \partial_k c_0(k_0, \mathbf{U}_0)$, we have

$$\widetilde{\mathbb{P}}_0^{-1} \, \mathbb{W} \, \widetilde{\mathbb{P}}_0 = \left( \begin{array}{cc|c} v_g & \widetilde{\mathfrak{a}}_0 & 0 \\ 0 & v_g & 0 \\ \hline 0 & 0 & -\mathbf{B} \, \nabla_{\mathbf{U}}^2 \mathcal{H}(\mathbf{U}_0, 0) \end{array} \right)$$

with

$$\widetilde{\mathfrak{a}}_0 := -k_0^4 (\partial_k c_0)^2 \, \mathfrak{a}_0 - k_0^2 \, \nabla_{\mathbf{U}} c_0^{\mathsf{T}} \, (\nabla_{\mathbf{U}}^2 \mathcal{H}(\mathbf{U}_0, 0) + c_0 \mathbf{B}^{-1})^{-1} \nabla_{\mathbf{U}} c_0$$
$$+ k_0^2 \, \nabla_{\mathbf{U}} c_0^{\mathsf{T}} \, (\nabla_{\mathbf{U}}^2 \mathcal{H}(\mathbf{U}_0, 0) + v_g \mathbf{B}^{-1})^{-1} \nabla_{\mathbf{U}} c_0$$

$$\widetilde{\mathbb{P}}_0 := \left( \begin{array}{cc|c} 1 & 0 & -k_0 \, \nabla_{\mathbf{U}} c_0^{\mathsf{T}} (\nabla_{\mathbf{U}}^2 \mathcal{H}(\mathbf{U}_0, 0) + v_g \mathbf{B}^{-1})^{-1} \mathbf{B}^{-1} \\ 0 & 1 & 0 \\ \hline 0 & k_0 \, (\nabla_{\mathbf{U}}^2 \mathcal{H}(\mathbf{U}_0, 0) + v_g \mathbf{B}^{-1})^{-1} \nabla_{\mathbf{U}} c_0 & \mathbf{I}_N \end{array} \right) \,.$$

As already announced, in general $\widetilde{\mathfrak{a}}_0$ is not zero and $\mathbb{W}$ possesses a Jordan block associated with $v_g$. In particular, in Appendix A we check that no vanishing occurs for the classical KdV equation.

Now we turn to the soliton limit. For concision's sake we first introduce

$$\Sigma_s = \begin{cases} \begin{pmatrix} 2\mathfrak{c}_s \, b^{-2} \end{pmatrix} & \text{if } N = 1 \\ \begin{pmatrix} 2\mathfrak{c}_s \sigma_s^2 & 2\mathfrak{c}_s \sigma_s w_s \\ 2\mathfrak{c}_s \sigma_s w_s & 2\mathfrak{c}_s w_s^2 + \sigma_s^2 \end{pmatrix} & \text{if } N = 2 \end{cases}$$

and

$$\mathbf{y}_s = \begin{cases} \begin{pmatrix} (\mathbb{S}^{-1} \mathbf{W}_s) \cdot \mathbb{O}_s \, \mathbb{S}^{-1} \mathbf{V}_s \end{pmatrix} & \text{if } N = 1 \\ \begin{pmatrix} (\mathbb{S}^{-1} \mathbf{T}_s) \cdot \mathbb{O}_s \, \mathbb{S}^{-1} \mathbf{V}_s \\ (\mathbb{S}^{-1} \mathbf{W}_s) \cdot \mathbb{O}_s \, \mathbb{S}^{-1} \mathbf{V}_s \end{pmatrix} & \text{if } N = 2 \end{cases}$$

---

[17] The left-hand side being evaluated at $(k_0, 0, \mathbf{U}_0)$ and the right-hand side at $(k_0, \mathbf{U}_0)$.

so that it follows from Theorem 2 and Corollary 3 that

$$\left(\begin{array}{c|c}(\varrho\ln\varrho)^2 & 0 \\ \hline 0 & \mathbf{I}_{N+1}\end{array}\right) k\,\mathbb{P}_s{}^\mathsf{T}(\nabla^2_{\mu,c,\boldsymbol{\lambda}}\Theta)\,(\mathbb{A}^\mathsf{T})^{-1}\mathbb{B}^{-1}$$

$$=-\left(\begin{array}{cc|c} -\frac{\Xi_s}{\pi}\mathfrak{h}_s\mathcal{Q}(\mathbf{Y}_s) & \frac{\pi}{\Xi_s}\mathfrak{h}_s & -\mathfrak{h}_s\mathbf{Y}_s{}^\mathsf{T}\mathbf{B}^{-1} \\ \partial^2_c\mathcal{M} & 0 & 0 \\ \hline \frac{\Xi_s}{\pi}\left(\mathbf{y}_s + \Sigma_s(\mathbf{A}_s{}^\mathsf{T})^{-1}\mathbf{Y}_s\right) & 0 & \Sigma_s(\mathbf{A}_s{}^\mathsf{T})^{-1}\end{array}\right)+\mathcal{O}\left(\frac{1}{(\ln\varrho)}\right)$$

and its inverse is

$$\left(\begin{array}{cc|c} 0 & -\frac{1}{\partial^2_c\mathcal{M}} & 0 \\ -\frac{\Xi_s}{\pi\mathfrak{h}_s} & \frac{\Xi_s^2}{\pi^2\partial^2_c\mathcal{M}}\left(\mathcal{Q}(\mathbf{Y}_s)+\mathbf{Y}_s{}^\mathsf{T}\mathbf{B}^{-1}\mathbf{A}_s{}^\mathsf{T}(\Sigma_s)^{-1}\mathbf{y}_s\right) & -\frac{\Xi_s}{\pi}\mathbf{Y}_s{}^\mathsf{T}\mathbf{B}^{-1}\mathbf{A}_s{}^\mathsf{T}(\Sigma_s)^{-1} \\ \hline 0 & \frac{\Xi_s}{\pi\partial^2_c\mathcal{M}}\left(\mathbf{Y}_s+\mathbf{A}_s{}^\mathsf{T}(\Sigma_s)^{-1}\mathbf{y}_s\right) & -\mathbf{A}_s{}^\mathsf{T}(\Sigma_s)^{-1}\end{array}\right)$$

$$+\mathcal{O}\left(\frac{1}{(\ln\varrho)}\right).$$

From this stems that $\frac{1}{k}\mathbb{B}\mathbb{A}^\mathsf{T}(\nabla^2_{\mu,c,\boldsymbol{\lambda}}\Theta)^{-1}\mathbb{A}$ equals

$$-\left(\begin{array}{cc|c} 0 & 0 & 0 \\ \frac{\Xi_s^2}{\pi^2}\mathbf{Y}_s{}^\mathsf{T}\mathbf{B}^{-1}\mathbf{A}_s{}^\mathsf{T}(\Sigma_s)^{-1}\mathbf{A}_s\mathbf{B}^{-1}\mathbf{Y}_s & 0 & \frac{\Xi_s}{\pi}\mathbf{Y}_s{}^\mathsf{T}\mathbf{B}^{-1}\mathbf{A}_s{}^\mathsf{T}(\Sigma_s)^{-1}\mathbf{A}_s\mathbf{B}^{-1} \\ \hline \frac{\Xi_s}{\pi}\mathbf{A}_s{}^\mathsf{T}(\Sigma_s)^{-1}\mathbf{A}_s\mathbf{B}^{-1}\mathbf{Y}_s & 0 & \mathbf{A}_s{}^\mathsf{T}(\Sigma_s)^{-1}\mathbf{A}_s\mathbf{B}^{-1}\end{array}\right)+\mathcal{O}\left(\frac{1}{(\ln\varrho)}\right).$$

In particular, on $\Lambda_s$,

$$\nabla^2_{k,\alpha,\mathbf{M}}\mathrm{H}=\left(\begin{array}{cc|c} \frac{\Xi_s^2}{\pi^2}\mathbf{Y}_s{}^\mathsf{T}\mathbf{B}^{-1}\mathbf{A}_s{}^\mathsf{T}(\Sigma_s)^{-1}\mathbf{A}_s\mathbf{B}^{-1}\mathbf{Y}_s & -c & \frac{\Xi_s}{\pi}\mathbf{Y}_s{}^\mathsf{T}\mathbf{B}^{-1}\mathbf{A}_s{}^\mathsf{T}(\Sigma_s)^{-1}\mathbf{A}_s\mathbf{B}^{-1} \\ -c & 0 & 0 \\ \hline \frac{\Xi_s}{\pi}\mathbf{B}^{-1}\mathbf{A}_s{}^\mathsf{T}(\Sigma_s)^{-1}\mathbf{A}_s\mathbf{B}^{-1}\mathbf{Y}_s & 0 & -c\mathbf{B}^{-1}+\mathbf{B}^{-1}\mathbf{A}_s{}^\mathsf{T}(\Sigma_s)^{-1}\mathbf{A}_s\mathbf{B}^{-1}\end{array}\right)$$

$$\mathbb{W}=\left(\begin{array}{cc|c} c & 0 & 0 \\ -\frac{\Xi_s^2}{\pi^2}\mathbf{Y}_s{}^\mathsf{T}\mathbf{B}^{-1}\mathbf{A}_s{}^\mathsf{T}(\Sigma_s)^{-1}\mathbf{A}_s\mathbf{B}^{-1}\mathbf{Y}_s & c & -\frac{\Xi_s}{\pi}\mathbf{Y}_s{}^\mathsf{T}\mathbf{B}^{-1}\mathbf{A}_s{}^\mathsf{T}(\Sigma_s)^{-1}\mathbf{A}_s\mathbf{B}^{-1} \\ \hline -\frac{\Xi_s}{\pi}\mathbf{A}_s{}^\mathsf{T}(\Sigma_s)^{-1}\mathbf{A}_s\mathbf{B}^{-1}\mathbf{Y}_s & 0 & c\mathbf{I}_N-\mathbf{A}_s{}^\mathsf{T}(\Sigma_s)^{-1}\mathbf{A}_s\mathbf{B}^{-1}\end{array}\right).$$

Translating again, by identification, into the notation of Corollary 5 yields the following result.

**Theorem 6.** *Under the assumptions of Theorem 4, let us still denote by* $\mathrm{H}$ *and* $\mathbb{W}$ *their extensions to* $\Lambda_s$. *Then for any* $(c_s,\mathbf{U}_s)$ *of* $\Lambda$, *we have*[18]

$$\nabla^2_{k,\alpha,\mathbf{M}}\mathrm{H}=\left(\begin{array}{cc|c} \nabla_\mathbf{U}\mathcal{M}^\mathsf{T}(\nabla^2_\mathbf{U}\mathcal{H}(\mathbf{U}_0,0)+c_s\mathbf{B}^{-1})^{-1}\nabla_\mathbf{U}\mathcal{M} & -c_s & \nabla_\mathbf{U}\mathcal{M}^\mathsf{T} \\ -c_s & 0 & 0 \\ \hline \nabla_\mathbf{U}\mathcal{M} & 0 & \nabla^2_\mathbf{U}\mathcal{H}(\mathbf{U}_0,0)\end{array}\right)$$

$$\mathbb{W}=\left(\begin{array}{cc|c} c_s & 0 & 0 \\ -\nabla_\mathbf{U}\mathcal{M}^\mathsf{T}(\nabla^2_\mathbf{U}\mathcal{H}(\mathbf{U}_0,0)+c_s\mathbf{B}^{-1})^{-1}\nabla_\mathbf{U}\mathcal{M} & c_s & -\nabla_\mathbf{U}\mathcal{M}^\mathsf{T} \\ \hline -\mathbf{B}\nabla_\mathbf{U}\mathcal{M} & 0 & -\mathbf{B}\nabla^2_\mathbf{U}\mathcal{H}(\mathbf{U}_0,0)\end{array}\right)$$

---

[18]The left-hand side being evaluated at $(0,\partial_c\mathcal{M}(c_s;\mathbf{U}_s),\mathbf{U}_s)$ and the right-hand side at $(c_s,\mathbf{U}_s)$.

*thus*

$$\widetilde{\mathbb{P}}_s^{-1}\,\mathbb{W}\,\widetilde{\mathbb{P}}_s \;=\; \left(\begin{array}{c|c} c_s\mathbf{I}_2 & 0 \\ \hline 0 & -\mathbf{B}\,\nabla_{\mathbf{U}}^2\mathcal{H}(\mathbf{U}_0,0) \end{array}\right)$$

*with*

$$\widetilde{\mathbb{P}}_s := \left(\begin{array}{cc|c} 1 & 0 & 0 \\ 0 & 1 & \nabla_{\mathbf{U}}\mathcal{M}^{\mathsf{T}}(\nabla_{\mathbf{U}}^2\mathcal{H}(\mathbf{U}_0,0)+c_s\mathbf{B}^{-1})^{-1}\mathbf{B}^{-1} \\ \hline -(\nabla_{\mathbf{U}}^2\mathcal{H}(\mathbf{U}_0,0)+c_s\mathbf{B}^{-1})^{-1}\nabla_{\mathbf{U}}\mathcal{M} & 0 & \mathbf{I}_N \end{array}\right).$$

**Remark 14.** In both limits, the $2\times 2$ principal block of $\nabla^2\mathrm{H}$ has negative determinant hence signature $(1,1)$, and therefore $\nabla^2\mathrm{H}$ is neither positive definite nor negative definite in either regime.

## 5 Asymptotics of the modulation eigenfields

Since limiting characteristic matrices exhibit double roots, we need to perform a higher-order asymptotic analysis so as to determine the hyperbolic nature of modulation systems not at the limit of interest but near the distinguished limit. We undertake this task now.

### 5.1 Small amplitude regime

In the harmonic regime, the $N$ eigenvalues arising from those of $-\mathbf{B}\nabla_{\mathbf{U}}^2\mathcal{H}(\mathbf{U}_0,0)$ may be analyzed by standard spectral perturbation analysis. We only need to blow up the two eigenvalues near $v_g$ and we shall do it by inverting and scaling $\mathbb{W}-v_g\mathbf{I}_{N+2}$ so as to reduce the problem to the spectral perturbation of simple eigenvalues.

The scaling process will reveal the prominent role played by some of the higher-order correctors not made explicit in Theorem 2. With this in mind, note that the proof of Theorem 2, in [BGMR20], also gives that under the assumptions of Theorem 2, $\nabla^3_{\mu,c,\boldsymbol{\lambda}}\Theta$ possesses a limit with convergence rate $\mathcal{O}(\delta)$ when $\delta\to 0$. This implies that $\Theta$ possesses as a function of $(\mu,c,\boldsymbol{\lambda})$ a $\mathcal{C}^3$ extension to the limit $\delta=0$ with convergence rate $\mathcal{O}(\delta)$. In turn this implies, under the assumptions of Theorem 3, that $\mathrm{H}$ as a function of $(k,\alpha,\mathbf{M})$ possesses a $\mathcal{C}^3$ extension to $\Lambda_0$, with convergence rate $\mathcal{O}(\delta)$. Then proceeding as in Subsection 4.2, we deduce that

$$\partial_k^2\mathrm{H}(k_0,\alpha,\mathbf{U}_0) = \partial^3_{kk\alpha}\mathrm{H}(k_0,0,\mathbf{U}_0)\,\alpha + \mathcal{O}(\delta^3),$$
$$\partial_k\nabla_{\mathbf{M}}\mathrm{H}(k_0,\alpha,\mathbf{U}_0) = \partial^2_{k\alpha}\nabla_{\mathbf{M}}\mathrm{H}(k_0,0,\mathbf{U}_0)\,\alpha + \mathcal{O}(\delta^3),$$

with

$$\begin{aligned}\partial^3_{kk\alpha}\mathrm{H}(k_0,0,\mathbf{U}_0) &= -\partial_k^2(k\,c_0)(k_0,\mathbf{U}_0) &&= -2\partial_k c_0(k_0,\mathbf{U}_0) - k_0\,\partial_k^2 c_0(k_0,\mathbf{U}_0),\\ \partial^2_{k\alpha}\nabla_{\mathbf{M}}\mathrm{H}(k_0,0,\mathbf{U}_0) &= -\partial_k\nabla_{\mathbf{U}}(k\,c_0)(k_0,\mathbf{U}_0) &&= -\nabla_{\mathbf{U}}c_0(k_0,\mathbf{U}_0) - k_0\,\partial_k\nabla_{\mathbf{U}}c_0(k_0,\mathbf{U}_0).\end{aligned}$$

As a consequence, with notation from Theorem 5, we have

$$\widetilde{\mathbb{P}}_0^{-1}\,\mathbb{W}(k_0,\alpha,\mathbf{U}_0)\,\widetilde{\mathbb{P}}_0 \;=\;$$

$$\left(\begin{array}{cc|c} v_g + \mathcal{O}(\delta^2) & \widetilde{\mathfrak{a}}_0 + \mathcal{O}(\delta^2) & \mathcal{O}(\delta^2) \\ -\partial^3_{kk\alpha}\mathrm{H}(k_0,0,\mathbf{U}_0)\,\alpha + \mathcal{O}(\delta^3) & v_g + \mathcal{O}(\delta^2) & -\partial^2_{k\alpha}\nabla_{\mathbf{M}}\mathrm{H}(k_0,0,\mathbf{U}_0)^{\mathsf{T}}\,\alpha + \mathcal{O}(\delta^3) \\ \hline \mathcal{O}(\delta^2) & \mathcal{O}(\delta^2) & -\mathbf{B}\,\nabla_{\mathbf{U}}^2\mathcal{H}(\mathbf{U}_0,0) + \mathcal{O}(\delta^2) \end{array}\right)$$

so that when $\widetilde{\mathfrak{a}}_0 \, \partial^2_{k\alpha} \nabla_{\mathbf{M}} H(k_0, 0, \mathbf{U}_0) \neq 0$

$$\frac{1}{\sqrt{\alpha}} \begin{pmatrix} 1 & 0 & 0 \\ 0 & \frac{1}{\sqrt{\alpha}} & 0 \\ \hline 0 & 0 & \mathbf{I}_N \end{pmatrix} \widetilde{\mathbb{P}}_0^{-1} \left[ \mathbb{W}(k_0, \alpha, \mathbf{U}_0) - v_g \, \mathbf{I}_{N+2} \right] \widetilde{\mathbb{P}}_0 \begin{pmatrix} 1 & 0 & 0 \\ 0 & \sqrt{\alpha} & 0 \\ \hline 0 & 0 & \mathbf{I}_N \end{pmatrix}$$

$$= \begin{pmatrix} \mathbf{I}_2 & 0 \\ \hline 0 & \frac{1}{\sqrt{\alpha}} \mathbf{I}_N \end{pmatrix} \left[ \begin{pmatrix} 0 & \widetilde{\mathfrak{a}}_0 & 0 \\ -\partial^3_{kk\alpha} H(k_0, 0, \mathbf{U}_0) & 0 & -\partial^2_{k\alpha} \nabla_{\mathbf{M}} H(k_0, 0, \mathbf{U}_0)^{\mathsf{T}} \\ \hline 0 & 0 & -\mathbf{B} \, \nabla^2_{\mathbf{U}} \mathcal{H}(\mathbf{U}_0, 0) - v_g \mathbf{I}_N \end{pmatrix} + \mathcal{O}(\delta) \right]$$

is invertible (provided that $\delta$ is sufficiently small) and its inverse is

$$\begin{pmatrix} 0 & -\frac{1}{\partial^3_{kk\alpha} H(k_0, 0, \mathbf{U}_0)} & 0 \\ \frac{1}{\widetilde{\mathfrak{a}}_0} & 0 & 0 \\ \hline 0 & 0 & 0 \end{pmatrix} + \mathcal{O}(\delta) \, .$$

At last we may apply elementary spectral perturbation theory to the latter matrix to study its two simple eigenvalues near $\pm 1/\sqrt{-\widetilde{\mathfrak{a}}_0 \partial^3_{kk\alpha} H(k_0, 0, \mathbf{U}_0)}$ (where here $\sqrt{\cdot}$ denotes any determination of the square root function). This leads to the following result.

**Theorem 7.** *Under the assumptions of Theorem 4, let us still denote by $H$ its extension to $\Lambda_0$ and consider $(k_0, 0, \mathbf{U}_0) \in \Lambda_0$, with associated linear group velocity*

$$v_g(k_0, \mathbf{U}_0) = -\partial^2_{k\alpha} H(k_0, 0, \mathbf{U}_0) = c_0(k_0, \mathbf{U}_0) + k_0 \partial_k c_0(k_0, \mathbf{U}_0) \, .$$

*Then in the small amplitude regime, the spectrum of the Whitham matrix $\mathbb{W}(k_0, \alpha, \mathbf{U}_0)$ contains*

1. *two eigenvalues near $v_g$, that expand as*

$$v_g \pm \sqrt{\alpha \Delta_{MI}} + \mathcal{O}(\alpha)$$

   *(where here $\sqrt{\cdot}$ denotes some determination of the square root function), with corresponding eigenvectors*

$$\begin{pmatrix} 1 + \mathcal{O}(\sqrt{\alpha}) \\ \mp \frac{\partial^3_{kk\alpha} H(k_0, 0, \mathbf{U}_0)}{\sqrt{\Delta_{MI}}} \sqrt{\alpha} + \mathcal{O}(\alpha) \\ \mathcal{O}(\sqrt{\alpha}) \end{pmatrix} ,$$

   *provided that the modulational-instability index $\Delta_{MI}(k_0, \mathbf{U}_0)$, given by*[19]

$$\Delta_{MI} := \Big( -k_0^5 (\partial_k c_0)^2 \, \partial^2_\mu \Theta(\mu_0, c_0, \boldsymbol{\lambda}_0) - k_0^2 \, \nabla_{\mathbf{U}} c_0^{\mathsf{T}} \, (\nabla^2_{\mathbf{U}} \mathcal{H}(\mathbf{U}_0, 0) + c_0 \mathbf{B}^{-1})^{-1} \nabla_{\mathbf{U}} c_0$$
$$+ k_0^2 \, \nabla_{\mathbf{U}} c_0^{\mathsf{T}} \, (\nabla^2_{\mathbf{U}} \mathcal{H}(\mathbf{U}_0, 0) + v_g \mathbf{B}^{-1})^{-1} \nabla_{\mathbf{U}} c_0 \Big)$$
$$\times \big( 2\partial_k c_0 + k_0 \, \partial^2_k c_0 \big)$$
$$= \Big( \partial^2_\alpha H - (\partial_\alpha \nabla_{\mathbf{M}} H)^{\mathsf{T}} (\nabla^2_{\mathbf{M}} H - (\partial^2_{k\alpha} H) \, \mathbf{B}^{-1})^{-1} \partial_\alpha \nabla_{\mathbf{M}} H \Big) \times \partial^3_{kk\alpha} H \, ,$$

   *is not zero ;*

2. *and $N$ eigenvalues near the eigenvalues of the dispersionless characteristic matrix $-\mathbf{B} \, \nabla^2_{\mathbf{U}} \mathcal{H}(\mathbf{U}_0, 0)$, that expand as*

$$z_j + \mathcal{O}(\alpha) \, , \qquad j \in \{1, N\}$$

---

[19] With evaluation either at $(k_0, \mathbf{U}_0)$ or at $(k_0, 0, \mathbf{U}_0)$, depending on terms.

*with associated eigenvectors*

$$\left(\begin{array}{c} -k_0 \nabla_{\mathbf{U}} c_0{}^{\mathsf{T}} (\nabla_{\mathbf{U}}^2 \mathcal{H}(\mathbf{U}_0, 0) + v_g \mathbf{B}^{-1})^{-1} \mathbf{B}^{-1} \mathbf{r}_j + \mathcal{O}(\alpha) \\ \mathcal{O}(\alpha) \\ \hline \mathbf{r}_j + \mathcal{O}(\alpha) \end{array}\right), \quad j \in \{1, N\}$$

*where $z_j$, $j \in \{1, N\}$, are the eigenvalues of $-\mathbf{B} \nabla_{\mathbf{U}}^2 \mathcal{H}(\mathbf{U}_0, 0)$, with corresponding eigenvectors, $\mathbf{r}_j$, $j \in \{1, N\}$, provided that these $N$ eigenvalues are distinct.*

*Moreover all the bounds are locally uniform with respect to $(k_0, \mathbf{U}_0)$.*

Note that the existence of an expansion into powers of $\sqrt{\alpha}$ of the eigenvalues of an $\mathcal{O}(\alpha)$ perturbation of a matrix possessing a double root from which they emerge is consistent with the general — worst-case — algebraicity theory for the spectrum of matrices.

**Remark 15.** Instead of using $\Delta_{MI}$, a simplified criterion on $\partial_\alpha^2 \mathrm{H} \times \partial_{kk\alpha}^3 \mathrm{H} = \partial_\alpha \omega \times \partial_k^2 \omega$ is sometimes incorrectly invoked. This is based on the deceptive guess that relevant conclusions may be derived from the consideration of the (artificially uncoupled) $2 \times 2$ block of the Whitham matrix concerning the wavenumber and the amplitude (see for instance [Whi99, p.490]). The reader is also referred to Remark A.i for a more concrete discussion of the latter.

**Remark 16.** We recall that it was proved in [BGNR14] that the failure of weak hyperbolicity of the modulation system does imply a slow side-band[20] instability of the background periodic wave, hence the use of the term modulational instability here. It follows from our analysis that such an instability occurs near the harmonic limit when the dispersionless system fails to be weakly hyperbolic or when $\mathsf{sign}(w_0) \Delta_{MI}$ is negative - recall from Corollary 4 and Remark 12 that the sign of $w_0$ dictates the one of $\alpha$ in the harmonic limit. For this reason, given its practical importance, we make the latter sign more explicit in Appendix A.

## 5.2 Small wavenumber regime

As in the harmonic regime, the $N$ eigenvalues arising from those of $-\mathbf{B} \nabla_{\mathbf{U}}^2 \mathcal{H}(\mathbf{U}_s, 0)$ may be analyzed by standard spectral perturbation analysis. We only need to blow up the two eigenvalues near $c$ and we shall do it by inverting and scaling $\mathbb{W} - c \mathbf{I}_{N+2}$.

To do so we first observe that, since $\mathbb{S} = \mathbb{A}\mathbb{B}\mathbb{A}^{\mathsf{T}}$ and $\mathbb{D}_s = \mathbb{P}_s{}^{\mathsf{T}} \mathbb{S} \mathbb{P}_s$, we have

$$\mathbb{W} - c \mathbf{I}_{N+2} = (\mathbb{P}_s{}^{\mathsf{T}} \mathbb{A})^{-1} \left( k \, \mathbb{P}_s{}^{\mathsf{T}} (\nabla_{\mu,c,\boldsymbol{\lambda}}^2 \Theta) \mathbb{P}_s \mathbb{D}_s^{-1} \right)^{-1} \mathbb{P}_s{}^{\mathsf{T}} \mathbb{A}$$

so that it is equivalent to study two blowing-up eigenvalues of $k \, \mathbb{P}_s{}^{\mathsf{T}} (\nabla_{\mu,c,\boldsymbol{\lambda}}^2 \Theta) \mathbb{P}_s \mathbb{D}_s^{-1}$.

Now for concision's sake we introduce

$$\mathbf{D}_s = \begin{cases} \begin{pmatrix} b^{-1} \end{pmatrix} & \text{if } N = 1 \\ \begin{pmatrix} 0 & \sigma_s \\ \sigma_s & w_s \end{pmatrix} & \text{if } N = 2 \end{cases}$$

and

$$\mathbf{x}_s = \begin{cases} \begin{pmatrix} \mathfrak{b}_s \, b^{-1} \end{pmatrix} & \text{if } N = 1 \\ \begin{pmatrix} \mathfrak{b}_s \, \sigma_s \\ \mathfrak{b}_s \, w_s + \mathfrak{c}_s \, \zeta_s \end{pmatrix} & \text{if } N = 2 \end{cases}.$$

---
[20]That is, with small spectral parameter and small Floquet exponent.

It follows from Theorem 2 and Corollary 3 that

$$\frac{\pi}{\Xi_s}\mathbb{P}_s{}^{\mathsf{T}}(\nabla^2_{\mu,c,\boldsymbol{\lambda}}\Theta)\mathbb{P}_s\mathbb{D}_s^{-1}$$
$$=\left(\begin{array}{cc|c} -\mathfrak{c}_s w_s\, b^{-1}\ln(\varrho) & \mathfrak{h}_s\frac{1+\varrho}{\varrho^2}+\mathfrak{a}_s\ln(\varrho) & \mathbf{x}_s{}^{\mathsf{T}}\mathbf{D}_s^{-1}\ln(\varrho) \\ 0 & -\mathfrak{c}_s w_s\, b^{-1}\ln(\varrho) & 0 \\ \hline 0 & \mathbf{x}_s\ln(\varrho) & \Sigma_s\,\mathbf{D}_s^{-1}\ln(\varrho)\end{array}\right)+\frac{\pi}{\Xi_s}\mathbb{P}_s{}^{\mathsf{T}}\mathbb{O}_s\mathbb{P}_s\mathbb{D}_s^{-1}+\mathcal{O}\left(\varrho\left(\ln\varrho\right)\right)$$

with

$$\mathbb{P}_s{}^{\mathsf{T}}\mathbb{O}_s\mathbb{P}_s\mathbb{D}_s^{-1}=\left(\begin{array}{cc|c} \frac{\Xi_s}{\pi}\partial_c^2\mathcal{M}(c_s,\mathbf{U}_s) & * & \mathbf{y}_s{}^{\mathsf{T}}\mathbf{D}_s^{-1} \\ * & * & * \\ \hline * & * & *\end{array}\right)$$

so that

$$\frac{\varrho}{\sqrt{1+\varrho}}\left(\begin{array}{cc|c} 1 & 0 & 0 \\ 0 & \frac{\sqrt{1+\varrho}}{\varrho} & 0 \\ \hline 0 & 0 & \mathbf{I}_N\end{array}\right)\frac{\pi}{\Xi_s}\mathbb{P}_s{}^{\mathsf{T}}(\nabla^2_{\mu,c,\boldsymbol{\lambda}}\Theta)\mathbb{P}_s\mathbb{D}_s^{-1}\left(\begin{array}{cc|c} 1 & 0 & 0 \\ 0 & \frac{\varrho}{\sqrt{1+\varrho}} & 0 \\ \hline 0 & 0 & \mathbf{I}_N\end{array}\right)$$
$$=\left(\begin{array}{cc|c} 0 & \mathfrak{h}_s & 0 \\ \frac{\Xi_s}{\pi}\partial_c^2\mathcal{M}(c_s,\mathbf{U}_s) & 0 & \mathbf{y}_s{}^{\mathsf{T}}\mathbf{D}_s^{-1} \\ \hline 0 & 0 & 0\end{array}\right)+\mathcal{O}\left(\varrho\left(\ln\varrho\right)\right).$$

To the latter matrix we may apply elementary spectral perturbation analysis to study the two simple eigenvalues arising from $\pm\sqrt{\mathfrak{h}_s\,\Xi_s\,\partial_c^2\mathcal{M}/\pi}$ (where here $\sqrt{\cdot}$ denotes any determination of the square root function).

**Theorem 8.** *Under the assumptions of Theorem 4, consider $(c_s,\mathbf{U}_s)\in\Lambda$ such that $\partial_c^2\mathcal{M}(c_s,\mathbf{U}_s)\neq 0$. Then in the large period regime, the spectrum of the Whitham matrix $\mathbb{W}(k,\partial_c^2\mathcal{M}(c_s,\mathbf{U}_s),\mathbf{U}_s)$ is given by*

1. *two eigenvalues expanding as*[21]

$$c_s\pm\frac{\varrho\sqrt{\pi}}{\sqrt{\mathfrak{h}_s\,\Xi_s\,\partial_c^2\mathcal{M}(c_s,\mathbf{U}_s),\mathbf{U}_s)}}+\mathcal{O}(\varrho^2\,\ln(\varrho))$$

*(where here $\sqrt{\cdot}$ denotes some determination of the square root function), with corresponding eigenvectors*[22]

$$\mathbb{P}_s{}^{\mathsf{T}}\mathbb{A}\left(\begin{array}{c} 1+\mathcal{O}(\varrho\ln(\varrho)) \\ \pm\varrho\sqrt{\frac{\Xi_s\,\partial_c^2\mathcal{M}(c_s,\mathbf{U}_s),\mathbf{U}_s)}{\mathfrak{h}_s\,\pi}}+\mathcal{O}(\varrho^2\,\ln(\varrho)) \\ \mathcal{O}(\varrho\ln(\varrho))\end{array}\right)$$

2. *and $N$ eigenvalues expanding as*

$$z_j+\mathcal{O}(k),\qquad j\in\{1,N\}$$

---

[21]We recall that $k\sim-\pi/(\Xi_s\ln(\varrho))$ in the solitary wave limit.
[22]We recall that

$$\mathbb{P}_s{}^{\mathsf{T}}\mathbb{A}=\left(\begin{array}{cc|c} -1/k & 0 & 0 \\ \mathcal{Q}(\mathbf{U}_s-\mathbf{M})/k & -k & (\mathbf{U}_s-\mathbf{M})^{\mathsf{T}}\mathbf{B}^{-1} \\ \hline \mathbf{A}_s\,\mathbf{B}^{-1}(\mathbf{U}_s-\mathbf{M})/k & 0 & \mathbf{A}_s\,\mathbf{B}^{-1}\end{array}\right)=\left(\begin{array}{cc|c} -1/k & 0 & 0 \\ \mathcal{O}(k) & -k & \mathcal{O}(k) \\ \hline \mathcal{O}(1) & 0 & \mathcal{O}(1)\end{array}\right).$$

*with associate eigenvectors*

$$\begin{pmatrix} \nabla_{\mathbf{U}}\mathcal{M}^{\mathsf{T}}(\nabla_{\mathbf{U}}^2\mathcal{H}(\mathbf{U}_0,0) + c_s\mathbf{B}^{-1})^{-1}\mathbf{B}^{-1}\mathbf{r}_j + \mathcal{O}(k) \\ \mathcal{O}(k) \\ \hline \mathbf{r}_j + \mathcal{O}(k) \end{pmatrix}, \qquad j \in \{1, N\}$$

*where $z_j$, $j \in \{1, N\}$, are the distinct and real eigenvalues of $-\mathbf{B}\nabla_{\mathbf{U}}^2\mathcal{H}(\mathbf{U}_s, 0)$, with corresponding eigenvectors, $\mathbf{r}_j$, $j \in \{1, N\}$.*

*Moreover all the bounds are locally uniform with respect to $(c_s, \mathbf{U}_s)$.*

**Remark 17.** Though diagonilizability of the limiting modulation systems has little direct impact on the hyperbolicity of modulation systems near the limit, in the reverse direction the expansions derived in Theorems 7 and 8 shed some light on the asymmetry between the harmonic and the soliton limits in terms of diagonalizability of the asymptotic systems. Indeed, in the latter limit, the convergence of the eigenvalues towards the double root occurs exponentially faster — as $\varrho \ln \varrho$ — than the convergence of eigenvectors, which converge as $1/\ln \varrho$ and this may be proved to imply *per se* persistence of diagonalizability at the limit. In contrast, in the former limit the perturbations of eigenvectors and eigenvalues are of the same order — namely $\delta$ — leaving room for a limiting Jordan block.

# Appendix

## A  Explicit formula for the modulational-instability index

The goal of this section is to make explicit both $\widetilde{\mathfrak{a}}_0$ and $\Delta_{MI}$ that are involved in the hyperbolicity of the Whitham system near or at the harmonic limit.

This requires the extraction from [BGMR20] of an explicit value for the coefficient $\mathfrak{a}_0$ in Theorem 2 (denoted $\alpha_0$ in [BGMR20]). First we recall from [BGMR20] that

$$\mathfrak{a}_0 := -\frac{1}{3}\frac{\partial_v^3\mathcal{W}(v_0)}{(\partial_v^2\mathcal{W}(v_0))^2}\frac{\partial_v\mathscr{Y}^0 + 2\partial_z\mathscr{Y}^0}{\mathscr{Y}^0} + \frac{2}{\partial_v^2\mathcal{W}(v_0)}\frac{\frac{1}{4}\partial_v^2\mathscr{Y}^0 + \partial_z^2\mathscr{Y}^0 - \partial_{wz}^2\mathscr{Y}^0}{\mathscr{Y}^0}$$

where

$$\mathscr{Y}(v,w,z) := \sqrt{\frac{2\kappa(v)}{\mathscr{R}(v,w,z)}},$$

$$\mathscr{R}(v,w,z) := \int_0^1\int_0^1 t\partial_v^2\mathcal{W}(w + t(z-w) + ts(v-z))\,\mathrm{d}s\mathrm{d}t.$$

Here we omit to specify the dependence of $\mathcal{W}$, $\mathscr{R}$ and $\mathscr{Y}$ on parameters $(c_0, \boldsymbol{\lambda}_0)$ since they are held fixed along the computation, and the exponent $^0$ denotes that functions of $(v, w, z)$ are evaluated at $(v_0, v_0, v_0)$.

First we recall from [BGMR20, Appendix B] that $\mathscr{R}$ is a symmetric function and we observe that

$$\mathscr{R}^0 = \frac{1}{2}\partial_v^2\mathcal{W}(v_0), \qquad\qquad \partial_v\mathscr{R}^0 = \frac{1}{6}\partial_v^3\mathcal{W}(v_0),$$

$$\partial_v^2\mathscr{R}^0 = \frac{1}{12}\partial_v^4\mathcal{W}(v_0), \qquad\qquad \partial_{wz}^2\mathscr{R}^0 = \frac{1}{24}\partial_v^4\mathcal{W}(v_0).$$

Moreover direct computations yield

$$\frac{\partial_v\mathscr{Y}^0}{\mathscr{Y}^0} = \frac{1}{2}\frac{\kappa'(v_0)}{\kappa(v_0)} - \frac{1}{2}\frac{\partial_v\mathscr{R}^0}{\mathscr{R}^0}, \qquad\qquad \frac{\partial_z\mathscr{Y}^0}{\mathscr{Y}^0} = -\frac{1}{2}\frac{\partial_v\mathscr{R}^0}{\mathscr{R}^0},$$

$$\frac{\partial_v\mathscr{Y}^0 + 2\partial_z\mathscr{Y}^0}{\mathscr{Y}^0} = \frac{1}{2}\frac{\kappa'(v_0)}{\kappa(v_0)} - \frac{3}{2}\frac{\partial_v\mathscr{R}^0}{\mathscr{R}^0} = \frac{1}{2}\frac{\kappa'(v_0)}{\kappa(v_0)} - \frac{1}{2}\frac{\partial_v^3\mathcal{W}(v_0)}{\partial_v^2\mathcal{W}(v_0)},$$

and

$$\frac{\partial_v^2 \mathscr{Y}^0}{\mathscr{Y}^0} = \left(\frac{\partial_v \mathscr{Y}^0}{\mathscr{Y}^0}\right)^2 + \frac{1}{2}\left(\frac{\kappa''(v_0)}{\kappa(v_0)} - \left(\frac{\kappa'(v_0)}{\kappa(v_0)}\right)^2\right) - \frac{1}{2}\left(\frac{\partial_v^2 \mathscr{R}^0}{\mathscr{R}^0} - \left(\frac{\partial_v \mathscr{R}^0}{\mathscr{R}^0}\right)^2\right)$$

$$= \frac{1}{2}\left(\frac{\kappa''(v_0)}{\kappa(v_0)} - \frac{1}{2}\left(\frac{\kappa'(v_0)}{\kappa(v_0)}\right)^2\right) - \frac{1}{2}\frac{\kappa'(v_0)}{\kappa(v_0)}\frac{\partial_v \mathscr{R}^0}{\mathscr{R}^0} - \frac{1}{2}\left(\frac{\partial_v^2 \mathscr{R}^0}{\mathscr{R}^0} - \frac{3}{2}\left(\frac{\partial_v \mathscr{R}^0}{\mathscr{R}^0}\right)^2\right),$$

$$\frac{\partial_z^2 \mathscr{Y}^0}{\mathscr{Y}^0} = -\frac{1}{2}\left(\frac{\partial_v^2 \mathscr{R}^0}{\mathscr{R}^0} - \frac{3}{2}\left(\frac{\partial_v \mathscr{R}^0}{\mathscr{R}^0}\right)^2\right),$$

$$\frac{\partial_{wz}^2 \mathscr{Y}^0}{\mathscr{Y}^0} = -\frac{1}{2}\left(\frac{\partial_{wz}^2 \mathscr{R}^0}{\mathscr{R}^0} - \frac{3}{2}\left(\frac{\partial_v \mathscr{R}^0}{\mathscr{R}^0}\right)^2\right) = -\frac{1}{2}\left(\frac{1}{2}\frac{\partial_v^2 \mathscr{R}^0}{\mathscr{R}^0} - \frac{3}{2}\left(\frac{\partial_v \mathscr{R}^0}{\mathscr{R}^0}\right)^2\right),$$

so that

$$\frac{\frac{1}{4}\partial_v^2 \mathscr{Y}^0 + \partial_z^2 \mathscr{Y}^0 - \partial_{wz}^2 \mathscr{Y}^0}{\mathscr{Y}^0} =$$

$$\frac{1}{8}\left(\frac{\kappa''(v_0)}{\kappa(v_0)} - \frac{1}{2}\left(\frac{\kappa'(v_0)}{\kappa(v_0)}\right)^2\right) - \frac{1}{8}\frac{\kappa'(v_0)}{\kappa(v_0)}\frac{\partial_v \mathscr{R}^0}{\mathscr{R}^0} - \frac{3}{8}\left(\frac{\partial_v^2 \mathscr{R}^0}{\mathscr{R}^0} - \frac{1}{2}\left(\frac{\partial_v \mathscr{R}^0}{\mathscr{R}^0}\right)^2\right)$$

$$= \frac{1}{8}\left(\frac{\kappa''(v_0)}{\kappa(v_0)} - \frac{1}{2}\left(\frac{\kappa'(v_0)}{\kappa(v_0)}\right)^2\right) - \frac{1}{24}\frac{\kappa'(v_0)}{\kappa(v_0)}\frac{\partial_v^3 \mathcal{W}(v_0)}{\partial_v^2 \mathcal{W}(v_0)} - \frac{1}{16}\left(\frac{\partial_v^4 \mathcal{W}(v_0)}{\partial_v^2 \mathcal{W}(v_0)} - \frac{1}{3}\left(\frac{\partial_v^3 \mathcal{W}(v_0)}{\partial_v^2 \mathcal{W}(v_0)}\right)^2\right)$$

thus

$$\partial_v^2 \mathcal{W}(v_0)\,\mathfrak{a}_0 =$$

$$\frac{1}{4}\left(\frac{\kappa''(v_0)}{\kappa(v_0)} - \frac{1}{2}\left(\frac{\kappa'(v_0)}{\kappa(v_0)}\right)^2\right) - \frac{1}{4}\frac{\kappa'(v_0)}{\kappa(v_0)}\frac{\partial_v^3 \mathcal{W}(v_0)}{\partial_v^2 \mathcal{W}(v_0)} - \frac{1}{8}\frac{\partial_v^4 \mathcal{W}(v_0)}{\partial_v^2 \mathcal{W}(v_0)} + \frac{5}{24}\left(\frac{\partial_v^3 \mathcal{W}(v_0)}{\partial_v^2 \mathcal{W}(v_0)}\right)^2.$$

To go further with computations we find it convenient to separate the scalar and system case.

### A.1  The scalar case

In the scalar case, note that the computations in the proof of Corollary 5 provide

$$\partial_v^2 \mathcal{W}(v_0) = -f''(v_0) - \frac{c_0(k_0, v_0)}{b} = (2\pi)^2\, k_0^2\, \kappa(v_0)$$

$$b\,\partial_v^2 \mathcal{H}(v_0, 0) + c_0(k_0, v_0) = -b\,\partial_v^2 \mathcal{W}(v_0)$$

$$b\,\partial_v^2 \mathcal{H}(v_0, 0) + v_g(k_0, v_0) = -3\,b\,\partial_v^2 \mathcal{W}(v_0)$$

and observe that when $\ell \geq 3$, $\partial_v^\ell \mathcal{W}(v_0) = -f^{(\ell)}(v_0)$. From this one readily derives

$$k_0\,\partial_k c_0(k_0, v_0) = -2\,b\,\partial_v^2 \mathcal{W}(v_0),$$

$$k_0\left(-2\partial_k c_0(k_0, v_0) - k_0 \partial_k^2 c_0(k_0, v_0)\right) = 6\,b\,\partial_v^2 \mathcal{W}(v_0),$$

$$\partial_{v_0} c_0(k_0, v_0) = b\,\partial_v^3 \mathcal{W}(v_0) - b\,\frac{\kappa'(v_0)}{\kappa(v_0)}\,\partial_v^2 \mathcal{W}(v_0),$$

so that

$$\frac{\widetilde{\mathfrak{a}}_0}{b^2 k_0^2\,\partial_v^2 \mathcal{W}(v_0)} = -4\,\partial_v^2 \mathcal{W}(v_0)\,\mathfrak{a}_0 + \frac{2}{3}\left(\frac{\partial_v^3 \mathcal{W}(v_0)}{\partial_v^2 \mathcal{W}(v_0)} - \frac{\kappa'(v_0)}{\kappa(v_0)}\right)^2$$

$$= -\left(\frac{\kappa''(v_0)}{\kappa(v_0)} - \frac{1}{2}\left(\frac{\kappa'(v_0)}{\kappa(v_0)}\right)^2\right) + \frac{\kappa'(v_0)}{\kappa(v_0)}\frac{\partial_v^3 \mathcal{W}(v_0)}{\partial_v^2 \mathcal{W}(v_0)} + \frac{1}{2}\frac{\partial_v^4 \mathcal{W}(v_0)}{\partial_v^2 \mathcal{W}(v_0)} - \frac{5}{6}\left(\frac{\partial_v^3 \mathcal{W}(v_0)}{\partial_v^2 \mathcal{W}(v_0)}\right)^2$$

$$= -\frac{\kappa''(v_0)}{\kappa(v_0)} + \frac{5}{6}\left(\frac{\kappa'(v_0)}{\kappa(v_0)}\right)^2 + \frac{1}{3}\frac{\kappa'(v_0)}{\kappa(v_0)}\frac{f'''(v_0)}{(2\pi)^2 k_0^2 \kappa(v_0)} - \frac{1}{6}\left(\frac{f'''(v_0)}{(2\pi)^2 k_0^2 \kappa(v_0)}\right)^2 + \frac{1}{2}\frac{f''''(v_0)}{(2\pi)^2 k_0^2 \kappa(v_0)}.$$

From the foregoing computations we also derive that

$$\Delta_{MI} = \widetilde{\mathfrak{a}}_0 \left(2\partial_k c_0 + k_0\, \partial_k^2 c_0\right) = 6\, b^3\, k_0\, (\partial_v^2 \mathcal{W}(v_0))^2 \tag{58}$$

$$\times \left[\frac{\kappa''(v_0)}{\kappa(v_0)} - \frac{5}{6}\left(\frac{\kappa'(v_0)}{\kappa(v_0)}\right)^2 - \frac{1}{3}\frac{\kappa'(v_0)}{\kappa(v_0)}\frac{f'''(v_0)}{(2\pi)^2 k_0^2 \kappa(v_0)} + \frac{1}{6}\left(\frac{f'''(v_0)}{(2\pi)^2 k_0^2 \kappa(v_0)}\right)^2 - \frac{1}{2}\frac{f''''(v_0)}{(2\pi)^2 k_0^2 \kappa(v_0)}\right].$$

Recall from Proposition 2 that in the scalar case the sign of $\alpha$ is given by the sign of $b$ so that we are interested in the sign of $b\,\Delta_{MI}$. We stress moreover that this sign may be determined by considering a second-order polynomial in the unknown $k_0^2$, that varies in $(0,\infty)$, with coefficients depending on $v_0$.

We leave this general discussion to the reader and focus now on the most classical case when $\kappa$ is constant. To begin, note that for the 'genuine' Korteweg-de Vries equation, $f$ is cubic and $\kappa$ is constant so that both $-\widetilde{\mathfrak{a}}_0$ and $b\,\Delta_{MI}$ are positive. Likewise, when $\kappa$ is constant, we have

- when either ($f'''(v_0) \neq 0$ and $f''''(v_0) = 0$) or $f''''(v_0) < 0$, $-\widetilde{\mathfrak{a}}_0$ and $b\,\Delta_{MI}$ are positive;

- when $f'''(v_0) = 0$ and $f''''(v_0) > 0$, $-\widetilde{\mathfrak{a}}_0$ and $b\,\Delta_{MI}$ are negative;

- when $f'''(v_0) \neq 0$ and $f''''(v_0) > 0$, the common sign of $-\widetilde{\mathfrak{a}}_0$ and $b\,\Delta_{MI}$ depends on the harmonic wavenumber $k_0$, modulational instability occurring for wavenumbers $k_0$ larger than the critical wavenumber

$$k_c(v_0) := \frac{1}{\sqrt{3}}\frac{|f'''(v_0)|}{2\pi \sqrt{\kappa(v_0)\, f''''(v_0)}}.$$

It is worth pointing out that the general case when $\kappa$ is arbitrary is richer and that there are situations when two critical wavenumbers appear in the analysis.

**Remark A.i.** For comparison, let us compute the incorrect index by extrapolating [Whi99, p.490]

$$b\, \partial_\alpha \omega \times \partial_k^2 \omega = b\, \partial_\alpha^2 \mathrm{H} \times \partial_{kk\alpha}^3 \mathrm{H}$$

$$= 6\, b^4\, k_0\, (\partial_v^2 \mathcal{W}(v_0))^2 \left(-\frac{\widetilde{\mathfrak{a}}_0}{b^2 k_0^2\, \partial_v^2 \mathcal{W}(v_0)} - \frac{1}{3}\left(\frac{f'''(v_0)}{\partial_v^2 \mathcal{W}(v_0)} + \frac{\kappa'(v_0)}{\kappa(v_0)}\right)^2\right)$$

$$= 6\, b^4\, k_0\, (\partial_v^2 \mathcal{W}(v_0))^2$$

$$\times \left[\frac{\kappa''(v_0)}{\kappa(v_0)} - \frac{7}{6}\left(\frac{\kappa'(v_0)}{\kappa(v_0)}\right)^2 - \frac{\kappa'(v_0)}{\kappa(v_0)}\frac{f'''(v_0)}{(2\pi)^2 k_0^2 \kappa(v_0)} - \frac{1}{6}\left(\frac{f'''(v_0)}{(2\pi)^2 k_0^2 \kappa(v_0)}\right)^2 - \frac{1}{2}\frac{f''''(v_0)}{(2\pi)^2 k_0^2 \kappa(v_0)}\right].$$

Note that even for the 'genuine' Korteweg-de Vries equation this predicts deceptively instability. More generally, when $\kappa$ is constant, the associated wrong criterion predicts instability when one of the following conditions is satisfied

- ($f'''(v_0) \neq 0$ and $f''''(v_0) = 0$) or $f''''(v_0) > 0$;

- $f''''(v_0) < 0$ and the wavenumber $k_0$ is smaller than

$$\frac{1}{\sqrt{3}}\frac{|f'''(v_0)|}{2\pi \sqrt{\kappa(v_0)\, |f''''(v_0)|}}.$$

This is in strong contrast with the correct conclusions drawn above.

## A.2 The system case

As a preliminary to computations in the system case, we recall that

$$b\,\tau(v)\,g(v;c,\lambda) = -c\,v - b\,\lambda\,,$$
$$\mathcal{W}(v;c,\boldsymbol{\lambda}) = -f(v) - \frac{1}{2}\,\tau(v)\,(g(v;c,\lambda_2))^2 - \frac{c}{b}\,v\,g(v;c,\lambda_2) - \boldsymbol{\lambda}\cdot(v,g(v;c,\lambda_2))\,,$$

so that

$$b\,\tau(v)\,\partial_v g(v;c,\lambda) = -c - b\,\tau'(v)\,g(v;c,\lambda)\,,$$
$$b\,\tau(v)\,\partial_v^2 g(v;c,\lambda) = -b\,\tau''(v)\,g(v;c,\lambda) - 2\,b\,\tau'(v)\,\partial_v g(v;c,\lambda)\,,$$
$$b\,\tau(v)\,\partial_v^3 g(v;c,\lambda) = -b\,\tau'''(v)\,g(v;c,\lambda) - 3\,b\,\tau''(v)\,\partial_v g(v;c,\lambda) - 3\,b\,\tau'(v)\,\partial_v^2 g(v;c,\lambda)\,,$$

and

$$\partial_v \mathcal{W}(v;c,\boldsymbol{\lambda}) = -f'(v) - \frac{1}{2}\,\tau'(v)\,(g(v;c,\lambda_2))^2 - \frac{c}{b}\,g(v;c,\lambda_2) - \lambda_1\,,$$
$$\partial_v^2 \mathcal{W}(v;c,\boldsymbol{\lambda}) = -\partial_v^2 \mathcal{H}((v,g(v;c,\lambda_2)),0) + \tau(v)\,(\partial_v g(v;c,\lambda_2))^2\,,$$
$$\partial_v^3 \mathcal{W}(v;c,\boldsymbol{\lambda}) = -\partial_v^3 \mathcal{H}((v,g(v;c,\lambda_2)),0) - 3\,\tau''(v)\,g(v;c,\lambda_2)\,\partial_v g(v;c,\lambda_2)$$
$$- 3\tau'(v)\,(\partial_v g(v;c,\lambda_2))^2\,,$$
$$\partial_v^4 \mathcal{W}(v;c,\boldsymbol{\lambda}) = -\partial_v^4 \mathcal{H}((v,g(v;c,\lambda_2)),0) - 4\,\tau'''(v)\,g(v;c,\lambda_2)\,\partial_v g(v;c,\lambda_2)$$
$$- 6\,\tau''(v)\,(\partial_v g(v;c,\lambda_2))^2 - 3\,\tau''(v)\,g(v;c,\lambda_2)\,\partial_v^2 g(v;c,\lambda_2)$$
$$- 6\tau'(v)\,\partial_v g(v;c,\lambda_2)\,\partial_v^2 g(v;c,\lambda_2)\,.$$

In particular, it follows recursively that for $\ell \geq 2$, $\partial_v^\ell g(v;c,\lambda_2)$ and $\partial_v^\ell \mathcal{W}(v;c,\boldsymbol{\lambda})$ may be written as functions of $v$, $g(v;c,\lambda_2)$ and $\partial_v g(v;c,\lambda)$, independently of $\boldsymbol{\lambda}$ and $c$. Note that it follows from Remark 13 that $\partial_v^2 \mathcal{H}(\mathbf{U}_0,0) < 0$ yields modulational instability through non hyperbolicity of the dispersionless system at $\mathbf{U}_0$ so that from now on when discussing signs of $\Delta_{MI}$ we may assume without mention that $\partial_v^2 \mathcal{H}(\mathbf{U}_0,0) \geq 0$ hence $\partial_v^2 \mathcal{W} \leq \tau(\partial_v g)^2$.

We also observe, essentially as in the proof of Corollary 5, that

$$\mathrm{tr}(\mathbf{B}\,\nabla_{\mathbf{U}}^2 \mathcal{H}((v,g(v;c,\lambda)),0) + c\,\mathbf{I}_N) = -2b\,\tau(v)\,\partial_v g(v;c,\lambda_2)\,,$$
$$\det(\mathbf{B}\,\nabla_{\mathbf{U}}^2 \mathcal{H}((v,g(v;c,\lambda)),0) + c\,\mathbf{I}_N) = b^2\,\tau(v)\,\partial_v^2 \mathcal{W}(v;c,\boldsymbol{\lambda})\,.$$

Now we point out that $c_0(k_0,\mathbf{U}_0)$ is defined by

$$\det(\mathbf{B}\,\nabla_{\mathbf{U}}^2 \mathcal{H}(\mathbf{U}_0,0) + c_0(k_0,\mathbf{U}_0)\,\mathbf{I}_N) = b^2\,\tau(v_0)\,(2\pi)^2\,k_0^2\,\kappa(v_0)$$

more explicitly written as

$$(c_0(k_0,\mathbf{U}_0))^2 + c_0(k_0,\mathbf{U}_0)\,\mathrm{tr}(\mathbf{B}\,\nabla_{\mathbf{U}}^2 \mathcal{H}(\mathbf{U}_0,0)) + \det(\mathbf{B}\,\nabla_{\mathbf{U}}^2 \mathcal{H}(\mathbf{U}_0,0)) = b^2\,\tau(v_0)\,(2\pi)^2\,k_0^2\,\kappa(v_0)$$

and that this definition makes sense if and only if

$$k_0^2 \geq \frac{\det(\mathbf{B}\,\nabla_{\mathbf{U}}^2 \mathcal{H}(\mathbf{U}_0,0)) - \frac{1}{4}(\mathrm{tr}(\mathbf{B}\,\nabla_{\mathbf{U}}^2 \mathcal{H}(\mathbf{U}_0,0)))^2}{b^2\,\tau(v_0)\,(2\pi)^2\,\kappa(v_0)} = -\frac{f''(v_0) + \frac{1}{2}\tau''(v_0)\,u_0^2}{(2\pi)^2\,\kappa(v_0)}\,.$$

Yet the latter prescribes a minimal value for $k_0$ only if $f''(v_0) + \frac{1}{2}\tau''(v_0)\,u_0^2 < 0$, that is, only if the corresponding dispersionless system already fails to be hyperbolic. Moreover when the inequality on $k_0^2$ is strictly satisfied there are actually two possible values for $c_0(k_0,\mathbf{U}_0)$. This defines two branches for $c_0$ and henceforth we follow one such branch.

By differentiating the relation defining $c_0$ with respect to $k$ we derive

$$\partial_v g(v_0)\, k_0\, \partial_k c_0(k_0, \mathbf{U}_0) = -b\, \partial_v^2 \mathcal{W}(v_0),$$

$$\partial_v g(v_0)\, k_0^2\, \partial_k^2 c_0(k_0, \mathbf{U}_0) = -b\, \partial_v^2 \mathcal{W}(v_0) + \frac{k_0^2 (\partial_k c_0(k_0, \mathbf{U}_0))^2}{b\, \tau(v_0)},$$

$$(-2\partial_k c_0(k_0, \mathbf{U}_0) - k_0\, \partial_k^2 c_0(k_0, \mathbf{U}_0)) = \frac{b\, \partial_v^2 \mathcal{W}\, (-\partial_v^2 \mathcal{W} + 3\tau (\partial_v g)^2)}{k_0\, \tau\, (\partial_v g)^3},$$

(where again here and from now on we omit to mark dependencies on $c$ and $\boldsymbol{\lambda}$ on $g$ and $\mathcal{W}$).

At this stage, we could differentiate with respect to $\mathbf{U}$ to compute $\nabla_{\mathbf{U}} c_0(k_0, \mathbf{U}_0)$ and conclude as in the scalar case. Yet, instead we shall directly use the relatively explicit formula derived in Subsection 4.3. The only missing piece to carry out this task is to extract from [BGMR20] a formula for the coefficient $\mathfrak{b}_0$ from Theorem 2 (denoted $\beta_0$ in [BGMR20]). With notation introduced above,

$$\mathfrak{b}_0 := -\frac{1}{3} \frac{\partial_v^3 \mathcal{W}(v_0)}{(\partial_v^2 \mathcal{W}(v_0))^2} + \frac{1}{\partial_v^2 \mathcal{W}(v_0)} \frac{\partial_v \mathscr{Y}^0}{\mathscr{Y}^0}$$

$$= -\frac{1}{3} \frac{\partial_v^3 \mathcal{W}(v_0)}{(\partial_v^2 \mathcal{W}(v_0))^2} + \frac{1}{\partial_v^2 \mathcal{W}(v_0)} \left( \frac{1}{2} \frac{\kappa'(v_0)}{\kappa(v_0)} - \frac{1}{6} \frac{\partial_v^3 \mathcal{W}(v_0)}{\partial_v^2 \mathcal{W}(v_0)} \right)$$

$$= \frac{1}{2} \frac{1}{\partial_v^2 \mathcal{W}(v_0)} \left( \frac{\kappa'(v_0)}{\kappa(v_0)} - \frac{\partial_v^3 \mathcal{W}(v_0)}{\partial_v^2 \mathcal{W}(v_0)} \right).$$

Now with notation from Subsection 4.3

$$\Sigma_0^{-1} = -\frac{1}{2\mathfrak{c}_0\, \sigma_0^4} \begin{pmatrix} 2\mathfrak{c}_0\, w_0^2 - \sigma_0^2 & -2\mathfrak{c}_0\, \sigma_0\, w_0 \\ -2\mathfrak{c}_0\, \sigma_0\, w_0 & 2\mathfrak{c}_0\, \sigma_0^2 \end{pmatrix}$$

$$-\frac{k_0^2}{\mathfrak{c}_0^2\, w_0^2} \mathbf{x}_0^\mathsf{T} \Sigma_0^{-1} \mathbf{x}_0 = \frac{k_0^2}{2\mathfrak{c}_0^3\, w_0^2} \left( -\mathfrak{b}_0^2 + 2\mathfrak{c}_0^3 \frac{\zeta_0^2}{\sigma_0^2} \right)$$

$$= \frac{\frac{1}{4} b^2 k_0^2}{\partial_v^2 \mathcal{W} (\partial_v g)^2} \left( -\left( \frac{\kappa'}{\kappa} \partial_v^2 \mathcal{W} - \partial_v^3 \mathcal{W} \right)^2 + (\partial_v^2 \mathcal{W})\, \tau\, (\partial_v^2 g)^2 \right)$$

thus at the harmonic limit

$$\partial_\alpha^2 \mathrm{H} = \frac{k_0^2}{\mathfrak{c}_0^2\, w_0^2} \left( \mathfrak{a}_0 - \mathbf{x}_0^\mathsf{T} (\Sigma_0)^{-1} \mathbf{x}_0 \right)$$

$$= \frac{\frac{1}{4} b^2 k_0^2}{\partial_v^2 \mathcal{W} (\partial_v g)^2} \Bigg[ \left( \frac{\kappa''}{\kappa} - \frac{1}{2} \left( \frac{\kappa'}{\kappa} \right)^2 \right) (\partial_v^2 \mathcal{W})^2 - \frac{\kappa'}{\kappa} \partial_v^3 \mathcal{W}\, \partial_v^2 \mathcal{W}$$

$$- \frac{1}{2} \partial_v^4 \mathcal{W}\, \partial_v^2 \mathcal{W} + \frac{5}{6} (\partial_v^3 \mathcal{W})^2 - \left( \frac{\kappa'}{\kappa} \partial_v^2 \mathcal{W} - \partial_v^3 \mathcal{W} \right)^2 + (\partial_v^2 \mathcal{W})\, \tau\, (\partial_v^2 g)^2 \Bigg]$$

$$= \frac{\frac{1}{4} b^2 k_0^2}{\partial_v^2 \mathcal{W} (\partial_v g)^2} \Bigg[ -\frac{1}{2} \partial_v^4 \mathcal{W}\, \partial_v^2 \mathcal{W} - \frac{1}{6} (\partial_v^3 \mathcal{W})^2 + \frac{\kappa'}{\kappa} \partial_v^3 \mathcal{W}\, \partial_v^2 \mathcal{W}$$

$$+ \left( \frac{\kappa''}{\kappa} - \frac{3}{2} \left( \frac{\kappa'}{\kappa} \right)^2 \right) (\partial_v^2 \mathcal{W})^2 + (\partial_v^2 \mathcal{W})\, \tau\, (\partial_v^2 g)^2 \Bigg].$$

To proceed we now consider at the harmonic limit

$$-(\partial_\alpha \nabla_\mathbf{M} \mathrm{H})^\mathsf{T} (\nabla_\mathbf{M}^2 \mathrm{H} - (\partial_{k\alpha}^2 \mathrm{H})\, \mathbf{B}^{-1})^{-1} \partial_\alpha \nabla_\mathbf{M} \mathrm{H}$$

$$= \frac{k_0^2}{\mathfrak{c}_0^2\, w_0^2} \mathbf{x}_0^\mathsf{T} (\Sigma_0)^{-1} \mathbf{A}_0 \mathbf{B}^{-1} \left( \mathbf{B}^{-1} \mathbf{A}_0^\mathsf{T} (\Sigma_0)^{-1} \mathbf{A}_0 \mathbf{B}^{-1} + \frac{1}{\mathfrak{c}_0 w_0} \mathbf{B}^{-1} \right)^{-1} \mathbf{B}^{-1} \mathbf{A}_0^\mathsf{T} (\Sigma_0)^{-1} \mathbf{x}_0$$

and observe that on one hand

$$\mathbf{A}_0 \mathbf{B}^{-1} = \begin{pmatrix} \sigma_0 & 0 \\ \frac{w_0}{2} & \frac{1}{b} \end{pmatrix}$$

$$\mathbf{B}^{-1}\mathbf{A}_0{}^\mathsf{T}(\Sigma_0)^{-1}\mathbf{A}_0\mathbf{B}^{-1} = -\frac{1}{2\mathfrak{c}_0\,\sigma_0^2} \begin{pmatrix} \frac{\mathfrak{c}_0\,w_0^2}{2} - \sigma_0^2 & -\frac{\mathfrak{c}_0\,w_0}{b} \\ -\frac{\mathfrak{c}_0\,w_0}{b} & \frac{2\mathfrak{c}_0}{b^2} \end{pmatrix}$$

$$\mathbf{B}^{-1}\mathbf{A}_0{}^\mathsf{T}(\Sigma_0)^{-1}\mathbf{A}_0\mathbf{B}^{-1} + \frac{1}{\mathfrak{c}_0 w_0}\mathbf{B}^{-1} = \frac{1}{2\mathfrak{c}_0\,\sigma_0^2} \begin{pmatrix} -\frac{\mathfrak{c}_0\,w_0^2}{2} + \sigma_0^2 & \frac{\mathfrak{c}_0\,w_0}{b} + \frac{2\sigma_0^2}{bw_0} \\ \frac{\mathfrak{c}_0\,w_0}{b} + \frac{2\sigma_0^2}{bw_0} & -\frac{2\mathfrak{c}_0}{b^2} \end{pmatrix}$$

$$\left(\mathbf{B}^{-1}\mathbf{A}_0{}^\mathsf{T}(\Sigma_0)^{-1}\mathbf{A}_0\mathbf{B}^{-1} + \frac{1}{\mathfrak{c}_0 w_0}\mathbf{B}^{-1}\right)^{-1} = \frac{b^2 \mathfrak{c}_0 w_0^2}{2\sigma_0^2 + 3\mathfrak{c}_0\,w_0^2} \begin{pmatrix} \frac{2\mathfrak{c}_0}{b^2} & \frac{\mathfrak{c}_0\,w_0}{b} + \frac{2\sigma_0^2}{bw_0} \\ \frac{\mathfrak{c}_0\,w_0}{b} + \frac{2\sigma_0^2}{bw_0} & \frac{\mathfrak{c}_0\,w_0^2}{2} - \sigma_0^2 \end{pmatrix}$$

and that on the other hand

$$(\Sigma_0)^{-1}\mathbf{x}_0 = -\frac{1}{2\mathfrak{c}_0\sigma_0^3}\begin{pmatrix} -\mathfrak{b}_0\sigma_0^2 - 2\mathfrak{c}_0^2\zeta_0 w_0 \\ 2\mathfrak{c}_0^2\sigma_0\zeta_0 \end{pmatrix}$$

$$\mathbf{B}^{-1}\mathbf{A}_0{}^\mathsf{T}(\Sigma_0)^{-1}\mathbf{x}_0 = \frac{1}{2\mathfrak{c}_0\sigma_0^2}\begin{pmatrix} \mathfrak{b}_0\sigma_0^2 + \mathfrak{c}_0^2\zeta_0 w_0 \\ -\frac{2\mathfrak{c}_0^2\zeta_0}{b} \end{pmatrix}.$$

Thus

$$-(\partial_\alpha \nabla_\mathbf{M} \mathrm{H})^\mathsf{T}(\nabla_\mathbf{M}^2 \mathrm{H} - (\partial_{k\alpha}^2 \mathrm{H})\,\mathbf{B}^{-1})^{-1}\partial_\alpha \nabla_\mathbf{M} \mathrm{H}$$

$$= \frac{b^2 k_0^2}{4\mathfrak{c}_0^3\sigma_0^4(2\sigma_0^2 + 3\mathfrak{c}_0\,w_0^2)}\left[\frac{2\mathfrak{c}_0}{b^2}(\mathfrak{b}_0\sigma_0^2 + \mathfrak{c}_0^2\zeta_0 w_0)^2 - \frac{4\mathfrak{c}_0^2\zeta_0}{b}(\mathfrak{b}_0\sigma_0^2 + \mathfrak{c}_0^2\zeta_0 w_0)\left(\frac{\mathfrak{c}_0\,w_0}{b} + \frac{2\sigma_0^2}{bw_0}\right)\right.$$

$$\left. + \frac{4\mathfrak{c}_0^4\zeta_0^2}{b^2}\left(\frac{\mathfrak{c}_0\,w_0^2}{2} - \sigma_0^2\right)\right]$$

$$= \frac{k_0^2}{2\mathfrak{c}_0^2 w_0 \sigma_0^2 (2\sigma_0^2 + 3\mathfrak{c}_0\,w_0^2)}\left[\mathfrak{b}_0^2\sigma_0^2 w_0 - 4\mathfrak{b}_0\mathfrak{c}_0\sigma_0^2\zeta_0 - 6\mathfrak{c}_0^3\zeta_0^2 w_0\right]$$

$$= \frac{\frac{1}{4}b^2 k_0^2 \tau}{\partial_v g(\partial_v^2 \mathcal{W})(\partial_v^2 \mathcal{W} + 3\tau\,(\partial_v g)^2)}\left[\partial_v g\left(\frac{\kappa'}{\kappa}\partial_v^2 \mathcal{W} - \partial_v^3 \mathcal{W}\right)^2\right.$$

$$\left. - 2\,\partial_v^2 g\,\partial_v^2 \mathcal{W}\left(\frac{\kappa'}{\kappa}\partial_v^2 \mathcal{W} - \partial_v^3 \mathcal{W}\right) - 3\partial_v^2 \mathcal{W}\,\tau\,\partial_v g\,(\partial_v^2 g)^2\right]$$

$$= \frac{\frac{1}{4}b^2 k_0^2 \tau}{\partial_v g(\partial_v^2 \mathcal{W})(\partial_v^2 \mathcal{W} + 3\tau\,(\partial_v g)^2)}\left[\partial_v g\,(\partial_v^3 \mathcal{W})^2 + 2\partial_v^3 \mathcal{W}\,\partial_v^2 \mathcal{W}\left(-\frac{\kappa'}{\kappa}\partial_v g + \partial_v^2 g\right)\right.$$

$$\left. + (\partial_v^2 \mathcal{W})^2\left(\left(\frac{\kappa'}{\kappa}\right)^2\partial_v g - 2\frac{\kappa'}{\kappa}\partial_v^2 g\right) - 3\partial_v^2 \mathcal{W}\,\tau\,\partial_v g\,(\partial_v^2 g)^2\right].$$

Finally

$$-\widetilde{\mathfrak{a}}_0 = \partial_\alpha^2 \mathrm{H} - (\partial_\alpha \nabla_\mathbf{M} \mathrm{H})^\mathsf{T} (\nabla_\mathbf{M}^2 \mathrm{H} - (\partial_{k\alpha}^2 \mathrm{H}) \, \mathbf{B}^{-1})^{-1} \partial_\alpha \nabla_\mathbf{M} \mathrm{H}$$

$$= \frac{\frac{1}{4} b^2 k_0^2}{(\partial_v g)^2 (\partial_v^2 \mathcal{W})(\partial_v^2 \mathcal{W} + 3\tau \, (\partial_v g)^2)} \bigg[ \tau (\partial_v g)^2 \, (\partial_v^3 \mathcal{W})^2$$

$$+ 2 \partial_v^3 \mathcal{W} \, \partial_v^2 \mathcal{W} \tau (\partial_v g) \left( -\frac{\kappa'}{\kappa} \partial_v g + \partial_v^2 g \right)$$

$$+ (\partial_v^2 \mathcal{W})^2 \tau \partial_v g \left( \left( \frac{\kappa'}{\kappa} \right)^2 \partial_v g - 2 \frac{\kappa'}{\kappa} \partial_v^2 g \right) - 3 \partial_v^2 \mathcal{W} \, \tau^2 \, (\partial_v g)^2 \, (\partial_v^2 g)^2$$

$$+ (\partial_v^2 \mathcal{W} + 3\tau \, (\partial_v g)^2) \bigg( -\frac{1}{2} \partial_v^4 \mathcal{W} \, \partial_v^2 \mathcal{W} - \frac{1}{6} \, (\partial_v^3 \mathcal{W})^2 + \frac{\kappa'}{\kappa} \partial_v^3 \mathcal{W} \, \partial_v^2 \mathcal{W}$$

$$+ \left( \frac{\kappa''}{\kappa} - \frac{3}{2} \left( \frac{\kappa'}{\kappa} \right)^2 \right) (\partial_v^2 \mathcal{W})^2 + (\partial_v^2 \mathcal{W}) \, \tau \, (\partial_v^2 g)^2 \bigg) \bigg]$$

$$= \frac{\frac{1}{4} b^2 k_0^2}{(\partial_v g)^2 (\partial_v^2 \mathcal{W})(\partial_v^2 \mathcal{W} + 3\tau \, (\partial_v g)^2)} \bigg[ -\frac{1}{2} \partial_v^4 \mathcal{W} \, \partial_v^2 \mathcal{W} (\partial_v^2 \mathcal{W} + 3\tau \, (\partial_v g)^2)$$

$$- \frac{1}{6} \, (\partial_v^3 \mathcal{W})^2 \, (\partial_v^2 \mathcal{W} - 3\tau \, (\partial_v g)^2) + \partial_v^3 \mathcal{W} \, \partial_v^2 \mathcal{W} \left( \frac{\kappa'}{\kappa} (\partial_v^2 \mathcal{W} + \tau \, (\partial_v g)^2) + 2\tau (\partial_v g) \partial_v^2 g \right)$$

$$+ \left( \frac{\kappa''}{\kappa} - \frac{3}{2} \left( \frac{\kappa'}{\kappa} \right)^2 \right) (\partial_v^2 \mathcal{W})^3$$

$$+ (\partial_v^2 \mathcal{W})^2 \left( \tau \, (\partial_v g)^2 \left( 3 \frac{\kappa''}{\kappa} - \frac{7}{2} \left( \frac{\kappa'}{\kappa} \right)^2 \right) - 2 \frac{\kappa'}{\kappa} \tau \partial_v g \, \partial_v^2 g + \tau (\partial_v^2 g)^2 \right) \bigg]$$

and

$$\frac{4\tau (\partial_v g)^5}{b^3 k_0} \frac{\partial_v^2 \mathcal{W} + 3\tau \, (\partial_v g)^2}{-\partial_v^2 \mathcal{W} + 3\tau \, (\partial_v g)^2} \times \Delta_{MI}$$

$$= -\frac{1}{2} \partial_v^4 \mathcal{W} \, \partial_v^2 \mathcal{W} (\partial_v^2 \mathcal{W} + 3\tau \, (\partial_v g)^2)$$

$$- \frac{1}{6} \, (\partial_v^3 \mathcal{W})^2 \, (\partial_v^2 \mathcal{W} - 3\tau \, (\partial_v g)^2) + \partial_v^3 \mathcal{W} \, \partial_v^2 \mathcal{W} \left( \frac{\kappa'}{\kappa} (\partial_v^2 \mathcal{W} + \tau \, (\partial_v g)^2) + 2\tau (\partial_v g) \partial_v^2 g \right)$$

$$+ \left( \frac{\kappa''}{\kappa} - \frac{3}{2} \left( \frac{\kappa'}{\kappa} \right)^2 \right) (\partial_v^2 \mathcal{W})^3$$

$$+ (\partial_v^2 \mathcal{W})^2 \left( \tau \, (\partial_v g)^2 \left( 3 \frac{\kappa''}{\kappa} - \frac{7}{2} \left( \frac{\kappa'}{\kappa} \right)^2 \right) - 2 \frac{\kappa'}{\kappa} \tau \partial_v g \, \partial_v^2 g + \tau (\partial_v^2 g)^2 \right).$$

Recall that $\alpha$ is of the sign of $w_0 = 2\, \partial_v g/b$ so that this is the sign of $\partial_v g \, \Delta_{MI}/b$, hence of the quantity written above, that matters here. We observe that in order to write this criterion directly in terms of $(k_0, \mathbf{U}_0)$, one may use that

$$\partial_v^2 \mathcal{W}(v_0) = (2\pi)^2 \, k_0^2 \, \kappa(v_0),$$

$$\partial_v g(v_0) = \pm \frac{1}{\sqrt{\tau(v_0)}} \sqrt{\partial_v^2 \mathcal{H}(\mathbf{U}_0, 0) + (2\pi)^2 \, k_0^2 \, \kappa(v_0)}$$

with the sign choice corresponding to the choice of a branch for $c_0(k_0, \mathbf{U}_0)$ and that all other quantities have already been expressed in terms of $\mathbf{U}_0$ and $\partial_v g(v_0)$. Note however that since $\partial_v g(v_0)$, thus $c_0$, is not a polynomial function of $k_0$ the range of possibilities is significantly harder to analyze in terms of $(k_0, \mathbf{U}_0)$ than in the scalar case. It may be preferable instead to express the criterion in terms of $(\mathbf{U}_0, \partial_v g(v_0))$.

Alternatively, since the general computations are somewhat tedious, from now on we shall rather make the extra assumption, satisfied by the most standard cases that $\tau$ is affine. This ensures that the expression to study is indeed a polynomial in $k_0^2$ with coefficients depending on $\mathbf{U}_0$. In this direction, note that in this case

$$\tau \left(\partial_v g(v_0)\right)^2 = f''(v_0) + \partial_v^2 \mathcal{W}(v_0),$$

$$\partial_v^2 g(v_0) = -2 \frac{\tau'(v_0)}{\tau(v_0)} \partial_v g(v_0),$$

$$\partial_v^3 \mathcal{W}(v_0) = -f'''(v_0) - 3 \frac{\tau'(v_0)}{\tau(v_0)} \tau(v_0)(\partial_v g(v_0))^2,$$

$$\partial_v^4 \mathcal{W}(v_0) = -f''''(v_0) + 12 \left(\frac{\tau'(v_0)}{\tau(v_0)}\right)^2 \tau(v_0)(\partial_v g(v_0))^2.$$

Thus under the same assumption the range of admissible parameters is described by $\partial_v^2 \mathcal{W} \geq -f''$ and we have

$$\frac{4\tau(\partial_v g)^5}{b^3 k_0} \frac{4\partial_v^2 \mathcal{W} + 3f''}{2\partial_v^2 \mathcal{W} + 3f''} \times \Delta_{MI} \tag{59}$$

$$= \left(\frac{1}{2}f'''' - 6\left(\frac{\tau'}{\tau}\right)^2 (\partial_v^2 \mathcal{W} + f'')\right) \partial_v^2 \mathcal{W}(4\partial_v^2 \mathcal{W} + 3f'')$$

$$+ \frac{1}{6} \left(f''' + 3\frac{\tau'}{\tau} (\partial_v^2 \mathcal{W} + f'')\right)^2 (2\partial_v^2 \mathcal{W} + 3f'')$$

$$- \left(f''' + 3\frac{\tau'}{\tau} (\partial_v^2 \mathcal{W} + f'')\right) \partial_v^2 \mathcal{W} \left(\frac{\kappa'}{\kappa} (2\partial_v^2 \mathcal{W} + f'') - 4\frac{\tau'}{\tau}(\partial_v^2 \mathcal{W} + f'')\right)$$

$$+ \left(\frac{\kappa''}{\kappa} - \frac{3}{2}\left(\frac{\kappa'}{\kappa}\right)^2\right) (\partial_v^2 \mathcal{W})^3$$

$$+ (\partial_v^2 \mathcal{W})^2(\partial_v^2 \mathcal{W} + f'') \left(3\frac{\kappa''}{\kappa} - \frac{7}{2}\left(\frac{\kappa'}{\kappa}\right)^2 + 4\frac{\kappa'}{\kappa}\frac{\tau'}{\tau} + 4\left(\frac{\tau'}{\tau}\right)^2\right)$$

$$= (\partial_v^2 \mathcal{W})^3 \left(-5\left(\frac{\tau'}{\tau}\right)^2 - 2\frac{\kappa'}{\kappa}\frac{\tau'}{\tau} - 5\left(\frac{\kappa'}{\kappa}\right)^2 + 4\frac{\kappa''}{\kappa}\right)$$

$$+ (\partial_v^2 \mathcal{W})^2 \left(f'' \left(-\frac{7}{2}\left(\frac{\tau'}{\tau}\right)^2 - 5\frac{\kappa'}{\kappa}\frac{\tau'}{\tau} - \frac{7}{2}\left(\frac{\kappa'}{\kappa}\right)^2 + 3\frac{\kappa''}{\kappa}\right) + f''' \left(6\frac{\tau'}{\tau} - 2\frac{\kappa'}{\kappa}\right) + 2f''''\right)$$

$$+ \partial_v^2 \mathcal{W} \left((f'')^2 \left(6\left(\frac{\tau'}{\tau}\right)^2 - 3\frac{\kappa'}{\kappa}\frac{\tau'}{\tau}\right) + f'' f''' \left(9\frac{\tau'}{\tau} - \frac{\kappa'}{\kappa}\right) + \frac{1}{3}(f''')^2 + \frac{3}{2}f'' f''''\right)$$

$$+ \frac{1}{2}f'' \left(f''' + 3\frac{\tau'}{\tau}f''\right)^2.$$

Recall that $\partial_v^2 \mathcal{W}(v_0) = \kappa(v_0)(2\pi)^2 k_0^2$ so that the latter expression is indeed a third-order polynomial expression in $k_0^2$ with coefficients depending on $v_0$, $k_0^2$ being allowed to vary in $(\max(\{0, -f''(v_0)/(\kappa(v_0)(2\pi)^2)\}), \infty)$ and that this is negativity of the expression that yields modulational instability.

Note that if one specializes to the cases arising from the hydrodynamic formulation of a nonlinear Schrödinger equation (see [BG13] for instance)

$$i \partial_t \psi = -\partial_x^2 \psi + f'(|\psi|^2) \psi,$$

then $\tau = \mathrm{Id}$ and $\kappa$ is given by $\kappa(v) = 1/(4v)$ so that the foregoing expression is reduced to $v_0^{-2}$ times the second-order polynomial

$$\begin{aligned}&(\partial_v^2 \mathcal{W})^2 \left(4v_0^{-2} f'' + 8v_0^{-1} f''' + 2 f''''\right) \\ &+ \partial_v^2 \mathcal{W} \left(9 v_0^{-2} (f'')^2 + 10 v_0^{-1} f'' f''' + \frac{1}{3}(f''')^2 + \frac{3}{2} f'' f''''\right) \\ &+ \frac{1}{2} f'' \left(f''' + 3 v_0^{-1} f''\right)^2.\end{aligned}$$

We remind the reader that $f''(v_0) < 0$ is already known to yield modulational instability through non hyperbolicity of the dispersionless system. We observe furthermore that in the case under consideration when $f''(v_0) > 0$, $f'''(v_0) \geq 0$ and $f''''(v_0) \geq 0$ then any $k_0$ is admissible and no modulational instability occurs. In particular for the hydrodynamic formulations of cubic Schrödinger equations, that is, when $f'$ is an affine function, modulational instability is completely decided by the sign of $f''(v_0)$ independently of $k_0$, that is, it is driven by the focusing/defocusing nature of the equation.

Going back to the general case (when $\kappa$ is arbitrary and $\tau$ is affine), we stress, as in Remark 5, the consistency of the foregoing computations with the Eulerian/mass Lagrangian conjugation (see [BG13, BGNR14]). To be more explicit, we denote with subscripts $_E$ and $_L$ quantities corresponding to each formulation. First we observe that $b_E = -1$ and $\tau_E = \mathrm{Id}$, whereas $b_L = 1$ and $\tau_L \equiv 0$. Moreover

$$f_L(v) = v\, f_E\left(\frac{1}{v}\right), \qquad \kappa_L(v) = \frac{1}{v^5} \kappa_E\left(\frac{1}{v}\right)$$

and at the harmonic limit

$$(v_L)_0 = \frac{1}{(v_E)_0}, \qquad (k_L)_0 = \frac{(k_E)_0}{(v_E)_0}.$$

Our observation is that when going from mass Lagrangian to Eulerian formulations the third-order polynomial is simply multiplied by $((v_L)_0)^{-11}$.

## B  About splitting of double roots

In the present section we give elementary results and address some warnings and caveats concerning the perturbation of double characteristics. Our goal is to shed some light on the asymptotics of modulation eigenfields in the harmonic and soliton limits.

For our present purpose it it sufficient to consider a $2 \times 2$-toy model. Thus we introduce an $\epsilon$-family of systems

$$\begin{cases} \partial_t A + v \partial_x A + (\widetilde{\mathfrak{a}} + \epsilon \delta') \partial_x B &= 0, \\ \partial_t B + \epsilon \delta \partial_x A + v \partial_x B &= 0, \end{cases} \qquad (60)$$

for some real constants $(v, \widetilde{\mathfrak{a}}, \delta, \delta')$ and the limiting system

$$\begin{cases} \partial_t A + v \partial_x A + \widetilde{\mathfrak{a}} \partial_x B &= 0, \\ \partial_t B + v \partial_x B &= 0. \end{cases} \qquad (61)$$

System (61) is hyperbolic if and only if $\widetilde{\mathfrak{a}} = 0$. Note that in this case Systems (60) are hyperbolic when $\epsilon > 0$ if and only if either $\delta \delta' > 0$ or ($\delta = 0$ and $\delta' = 0$). In the former case, pairs of eigenvalues/eigenvectors may be chosen as

$$v \pm \epsilon \sqrt{\delta\, \delta'}, \qquad \begin{pmatrix} 1 \\ \pm\sqrt{\frac{\delta}{\delta'}} \end{pmatrix}.$$

We now turn to the case when $\widetilde{\mathfrak{a}} \neq 0$ so that System (61) fails to be hyperbolic. In this case Systems (60) are hyperbolic when $\epsilon > 0$ is sufficiently small if and only if $\delta\widetilde{\mathfrak{a}} > 0$. Under this condition, pairs of eigenvalues/eigenvectors may be chosen so as to expand as

$$v \pm \sqrt{\epsilon\delta\widetilde{\mathfrak{a}}} + \mathcal{O}(\epsilon^{\frac{3}{2}}), \qquad \begin{pmatrix} 1 \\ \pm\sqrt{\epsilon\frac{\delta}{\mathfrak{a}}} + \mathcal{O}(\epsilon^{\frac{3}{2}}) \end{pmatrix}$$

when $\epsilon \to 0$.

We would like to add a few comments on the case $\delta\widetilde{\mathfrak{a}} > 0$ and we enforce this condition from now on. We first observe that System (61) leaves invariant the line $B = 0$ and that its restriction to this line reduces to $\partial_t A + v\partial_x A = 0$ hence is hyperbolic. Our second observation is that one may obviously solve System (61) by first solving $\partial_t B + v\partial_x B = 0$ and then recover $A$ from $B$; this leads to the fact that, for any $s \in \mathbb{R}$, though System (61) is ill-posed in $H^s(\mathbb{R}) \times H^s(\mathbb{R})$, it is well-posed in $H^{s+1}(\mathbb{R}) \times H^s(\mathbb{R})$. Our last observation is that there are $\epsilon$-dependent change of variables that provides a form of System (60) converging to an hyperbolic diagonal system. Yet these change of variables are singular in one way or the other and our claim is that they lead to spurious conclusions.

We first illustrate this with a rather naive change of variables. For $\epsilon > 0$, let us set $(\widetilde{A}_\epsilon, B) := (\sqrt{\epsilon}A, B)$ and observe that System (60) takes the alternative form

$$\begin{cases} \partial_t \widetilde{A}_\epsilon + v\partial_x \widetilde{A}_\epsilon + (\sqrt{\epsilon}\widetilde{\mathfrak{a}} + \epsilon^{\frac{3}{2}}\delta')\partial_x B &= 0, \\ \partial_t B + \sqrt{\epsilon}\delta\partial_x \widetilde{A}_\epsilon + v\partial_x B &= 0. \end{cases} \tag{62}$$

This suggests as a limiting system

$$\begin{cases} \partial_t \widetilde{A} + v\partial_x \widetilde{A} &= 0, \\ \partial_t B + v\partial_x B &= 0. \end{cases} \tag{63}$$

Yet an inspection of initial data shows that if we were interested in solving System (60) with $(A, B)(0, \cdot) = (A^{(0)}, B^{(0)})$, then we are ending up solving (63) with $(\widetilde{A}, B)(0, \cdot) = (0, B^{(0)})$, hence essentially

$$\begin{cases} 0 &= 0, \\ \partial_t B + v\partial_x B &= 0. \end{cases}$$

In disguise we have simply multiplied the first equation of System (60) by $\sqrt{\epsilon}$ so as to drop the equation at the limit.

There are more subtle ways to arrive at similarly deceptive systems. Let us examine what happens if one insists in diagonalizing System (60) either exactly or at main order. To do so, consider for $\epsilon > 0$,

$$(A_\epsilon^+, A_\epsilon^-) := \left(\sqrt{\epsilon\frac{\delta}{\mathfrak{a}}}A + B, -\sqrt{\epsilon\frac{\delta}{\mathfrak{a}}}A + B\right)$$

and observe that System (60) takes the alternative form

$$\begin{cases} \partial_t A_\epsilon^+ + \left(v + \sqrt{\epsilon\delta\widetilde{\mathfrak{a}}} + \frac{1}{2}\epsilon^{\frac{3}{2}}\sqrt{\frac{\delta}{\mathfrak{a}}}\delta'\right)\partial_x A_\epsilon^+ + \frac{1}{2}\epsilon^{\frac{3}{2}}\sqrt{\frac{\delta}{\mathfrak{a}}}\delta'\partial_x A_\epsilon^- &= 0, \\ \partial_t A_\epsilon^- + \left(v - \sqrt{\epsilon\delta\widetilde{\mathfrak{a}}} - \frac{1}{2}\epsilon^{\frac{3}{2}}\sqrt{\frac{\delta}{\mathfrak{a}}}\delta'\right)\partial_x A_\epsilon^- - \frac{1}{2}\epsilon^{\frac{3}{2}}\sqrt{\frac{\delta}{\mathfrak{a}}}\delta'\partial_x A_\epsilon^+ &= 0. \end{cases} \tag{64}$$

This suggests as a limiting system

$$\begin{cases} \partial_t A^+ + v\partial_x A^+ &= 0, \\ \partial_t A^- + v\partial_x A^- &= 0. \end{cases} \tag{65}$$

Yet an inspection of initial data shows that this is really

$$\begin{cases} \partial_t B + v\partial_x B &= 0, \\ \partial_t B + v\partial_x B &= 0, \end{cases}$$

in disguise.

We stress again that in contrast with what happens with the foregoing sets of variables, Theorem 3 ensures that the set of variables $(k, \alpha, \mathbf{M})$ does not suffer from an hidden reduction of dimension.

## C  Symbolic index

### C.1  Table of symbols

## C.2 Specialization to gKdV

To help the reader gain some slightly more concrete acquaintance with various general symbols, we specialize here some of the general formulas to the case of the generalized Korteweg–de Vries equation

$$\partial_t v - \partial_x(f'(v)) + \partial_x^3 v = 0\,.$$

In this case

$$\mathcal{H} = \frac{1}{2}v_x^2 + f(v)\,, \qquad \kappa \equiv 1\,, \qquad \mathbf{B} = b = 1\,, \qquad \mathcal{Q}(v) = \frac{1}{2}v^2\,,$$

$$\mathcal{W} = -f(v) - c\frac{v^2}{2} - \lambda v\,, \qquad \mathbf{M} = \langle \underline{v} \rangle\,, \qquad \mathbb{B} = \begin{pmatrix} 0 & 1 & 0 \\ 1 & 0 & 0 \\ 0 & 0 & 1 \end{pmatrix}\,,$$

$$\alpha = \frac{1}{2}\int_0^{\Xi} \left(\underline{v}(\xi)^2 - \langle \underline{v} \rangle^2\right) \mathrm{d}\xi = \frac{1}{2}\int_0^{\Xi} \left(\underline{v}(\xi) - \langle \underline{v} \rangle\right)^2 \mathrm{d}\xi\,,$$

$$\Theta = \int_0^{\Xi} \left(\frac{1}{2}\underline{v}_x^2 + f(\underline{v}) + c\frac{\underline{v}^2}{2} + \lambda \underline{v} + \mu\right) \mathrm{d}\xi\,,$$

$$\mathcal{M} = \int_{-\infty}^{+\infty} \left(\frac{1}{2}(\underline{v}_x^s(\xi))^2 + f(\underline{v}^s(\xi)) - f(v_s) - f'(v_s)(\underline{v}^s(\xi) - v_s) + c\frac{(\underline{v}^s(\xi) - v_s)^2}{2}\right) \mathrm{d}\xi\,,$$

$$\partial_c \mathcal{M} = \frac{1}{2}\int_{-\infty}^{+\infty} (\underline{v}^s(\xi) - v_s)^2 \mathrm{d}\xi\,,$$

$$\mathrm{H} = \frac{1}{\Xi}\int_0^{\Xi} \left(\frac{1}{2}\underline{v}_x^2 + f(\underline{v})\right) \mathrm{d}\xi\,,$$

$$\mathbf{V}_i = \begin{pmatrix} 1 \\ \frac{1}{2}v_i^2 \\ v_i \end{pmatrix}\,, \qquad \mathbf{W}_i = \begin{pmatrix} 0 \\ v_i \\ 1 \end{pmatrix}\,, \qquad \mathbf{Z}_i = \begin{pmatrix} 0 \\ 1 \\ 0 \end{pmatrix}\,, \qquad \mathbf{E} = \begin{pmatrix} 1 \\ 0 \\ 0 \end{pmatrix}\,, \qquad \mathbf{F} = \begin{pmatrix} 0 \\ -1 \\ 0 \end{pmatrix}\,,$$

$$\mathbb{S} = \begin{pmatrix} 0 & -1 & 0 \\ -1 & 0 & 0 \\ 0 & 0 & 1 \end{pmatrix}\,, \qquad \mathbb{D}_i = \begin{pmatrix} 0 & 1 & 0 \\ 1 & 0 & 0 \\ 0 & 0 & 1 \end{pmatrix}\,,$$

$$\mathbb{P}_i = \begin{pmatrix} 1 & -\frac{v_i^2}{2} & -v_i \\ 0 & -1 & 0 \\ 0 & v_i & 1 \end{pmatrix}\,, \qquad \mathbb{A} = \begin{pmatrix} -\frac{1}{k} & 0 & 0 \\ -\frac{1}{k}\frac{\langle \underline{v} \rangle^2}{2} & k & \langle \underline{v} \rangle \\ -\frac{1}{k}\langle \underline{v} \rangle & 0 & 1 \end{pmatrix}\,,$$

$$w_i = 1\,, \qquad \Sigma_i = 2\mathfrak{c}_i\,, \qquad \mathbf{x}_i = \mathfrak{b}_i\,, \qquad \mathbf{D}_i = 1\,,$$

$$\mathbf{Y}_0 = \frac{\mathfrak{b}_0}{4\mathfrak{c}_0}\,, \qquad -\mathbf{B}\nabla_\mathbf{U}^2 \mathcal{H}(v_0, 0) = -f''(v_0)\,,$$

$$c_0 = -f''(v_0) - (2\pi)^2 k_0^2\,, \qquad v_g = -f''(v_0) - 3(2\pi)^2 k_0^2\,,$$

$$\widetilde{\mathfrak{a}}_0 = -\frac{1}{6}\frac{(f'''(v_0))^2}{(2\pi)^2} + \frac{1}{2}f''''(v_0)\,k_0^2\,, \qquad \Delta_{MI} = k_0\left[(f'''(v_0))^2 - 3(2\pi)^2 f''''(v_0)\,k_0^2\right]\,.$$

Sylvie Benzoni-Gavage
Univ Lyon, Université Claude Bernard Lyon 1, CNRS UMR 5208,
Institut Camille Jordan, F-69622 Villeurbanne cedex, France
benzoni@math.univ-lyon1.fr

Colin Mietka
Univ Lyon, Université Claude Bernard Lyon 1, CNRS UMR 5208,
Institut Camille Jordan, F-69622 Villeurbanne cedex, France
mietka@math.univ-lyon1.fr

L. Miguel Rodrigues
Univ Rennes & IUF, CNRS, IRMAR - UMR 6625, F-35000 Rennes, France
luis-miguel.rodrigues@univ-rennes1.fr